\documentclass[12pt]{report}

\input{mathPreambuleReport}
\usepackage[]{natbib} 

\usepackage{makeidx}
\makeindex

\begin{document}
\title{Lecture Notes on Free Probability}
\author{Vladislav Kargin}
\date{June 2025}
\maketitle

%
%
%
%

\tableofcontents


\chapter{Non-commutative probability spaces and distributions}

\section{Non-commutative probability spaces}

Let \(\mathcal{A}\) be an algebra of bounded linear operators acting on a Hilbert space \(H\). 
We assume that \(\mathcal{A}\) contains the identity operator (such algebras 
are called \emph{unital}) and that it is closed under taking adjoints, i.e.\ 
if \(X \in \mathcal{A}\), then \(X^{\ast} \in \mathcal{A}\). A (uniform) 
operator norm on \(\mathcal{A}\) is given by
\[
\|X\| \;=\; \sup_{\|v\|=1} \|X\,v\|.
\]

It is often convenient to assume further that \(\mathcal{A}\) is closed with 
respect to this norm topology, in which case \(\mathcal{A}\) is called a 
\(C^\ast\)-algebra. One can also assume that \(\mathcal{A}\) is closed in the 
weak operator topology (that is, \(X_i \to X\) if and only if 
\(\langle u, X_i\,v\rangle \to \langle u, X\,v\rangle\) for all \(u,v \in H\)), 
in which case \(\mathcal{A}\) is called a \(W^\ast\)-algebra or 
\emph{von Neumann algebra}.

\paragraph{Example (commutative algebras).}
\begin{itemize}
\item A typical commutative \(C^\ast\)-algebra is \(C([0,1])\), the algebra 
      of continuous functions on \([0,1]\) with the supremum norm.
\item A typical commutative \(W^\ast\)-algebra is \(L^\infty([0,1])\), the 
      algebra of essentially bounded measurable functions on \([0,1]\).
\end{itemize}
Note that a sequence of continuous functions might converge weakly to a 
discontinuous function, so \(C([0,1])\) is not closed in the weak operator 
topology and hence is not a von Neumann algebra.

\vspace{1em}

We define a \emph{state} on an operator algebra \(\mathcal{A}\) to be a 
continuous linear functional \(\phi : \mathcal{A} \to \mathbb{C}\) satisfying 
the positivity property:
\[
\phi\bigl(X^*\,X\bigr) \;\ge\; 0 
\quad\text{for all } X \in \mathcal{A}.
\]
A typical example of a state is given by
\[
\phi(X) \;=\; \langle u,\, X\,u\rangle
\]
for a unit vector \(u \in H\).  (Indeed, if \(\phi\) is a state on a 
\(C^\ast\)-algebra, one can build a suitable representation of \(\mathcal{A}\) 
on a Hilbert space so that \(\phi\) appears in this form; this is the 
well-known \emph{GNS construction}.)

The name ``state'' reflects the connection to quantum mechanics: if a quantum 
system is in the unit vector state \(u\), then each self-adjoint operator \(X\) 
represents an observable, and its expected measurement outcome is 
\(\langle u,\, X\,u\rangle\).

In what follows, we often use the words ``state'' and ``expectation'' 
interchangeably.

States can have extra properties:
\begin{itemize}
\item \(\phi\) is called \emph{faithful} if \(\phi(A^*\,A)=0\) implies \(A=0\).
\item \(\phi\) is called \emph{normal} if \(X_n \to X\) weakly implies 
      \(\phi(X_n)\to \phi(X)\). 
\item \(\phi\) is called \emph{tracial} if \(\phi(XY) = \phi(YX)\) for all 
      \(X,Y\in \mathcal{A}\). A tracial state is also called a \emph{trace}.
\end{itemize}

\index{non-commutative probability space|textbf}
\begin{defi}
A \emph{non-commutative probability space} \(\bigl(\mathcal{A}, \phi\bigr)\) 
consists of a unital operator algebra \(\mathcal{A}\) and a state \(\phi\) 
such that \(\phi(I)=1\), where \(I\) is the identity operator in 
\(\mathcal{A}\).
\end{defi}

When \(\phi\) is tracial, we call \(\bigl(\mathcal{A},\phi\bigr)\) a 
\emph{tracial non-commutative probability space}. If \(\mathcal{A}\) is a 
\(C^\ast\)-algebra and \(\phi\) is faithful, then we say 
\(\bigl(\mathcal{A},\phi\bigr)\) is a \emph{\(C^\ast\)-probability space}. 
If \(\mathcal{A}\) is a von Neumann algebra and \(\phi\) is normal, we call 
\(\bigl(\mathcal{A},\phi\bigr)\) a \emph{\(W^\ast\)-probability space}. 
As in classical probability, the possibility of taking limits (in either norm 
or weak topology) can make the \(C^\ast\)- or \(W^\ast\)-setting more 
complicated and more interesting.

\subsection*{Examples}

\begin{exa}[A classical probability space]
Let \(\bigl(\Omega,\mathfrak{A},\mu\bigr)\) be a classical probability space. 
Take \(H = L^2(\Omega,\mu)\) and let \(\mathcal{A}\) be the algebra of bounded 
measurable functions on \(\Omega\), acting by pointwise multiplication on \(H\). 
Define
\[
\phi(f) \;=\; \int_{\Omega} f \,\mu(d\omega).
\]
This coincides with the usual expectation of \(f\).
\end{exa}

\begin{exa}[The algebra of \(N \times N\) matrices]
Here \(\mathcal{A}\) is the algebra of all \(N\times N\) complex matrices. 
A natural state on \(\mathcal{A}\) is the \emph{normalized trace}:
\[
\phi(X) \;=\; \frac{1}{N} \sum_{i=1}^N X_{ii}.
\]
This is clearly tracial, since \(\mathrm{tr}(XY) = \mathrm{tr}(YX)\).
\end{exa}

\begin{exa}[The algebra of random \(N \times N\) matrices]
\label{example_random_matrices}
Suppose we have random \(N\times N\) matrices \(X\). As long as all joint 
moments of the entries are finite, one can define a state by
\[
\phi(X) \;=\; \mathbb{E}\bigl[\mathrm{tr}(X)\bigr].
\]
In physics literature, the expectation of a random variable \(a\) is often 
denoted by \(\langle a \rangle\). Thus \(\phi\) acts by 
``take the trace, then take the ensemble average.''
\end{exa}

\begin{exa}[The group algebra of a finitely generated group]
\label{example_group_algebra}
Let \(G\) be a finitely generated group, and let \(\mathbb{C}G\) be its group 
algebra. One embeds \(\mathbb{C}G\) into the bounded operators on \(\ell^2(G)\) 
via the \emph{left regular representation} and takes either norm closure or weak closure in $B(\ell^2(G))$ to obtain $C^\ast$ or $W^\ast$ algebra, respectively. Let \(e \in G\) be the identity 
element. The state on \(\mathbb{C}G\) is given by
\[
\phi(X) \;=\; \langle e,\,X\,e\rangle,
\]
which is just the coefficient of the identity element in \(X\). This state is 
both tracial and faithful.
\end{exa}

\index{Toeplitz--Cuntz algebra|textbf}
\index{Fock space|textbf}
\index{Cuntz algebra|textbf}
\begin{exa}[The Toeplitz--Cuntz algebra]
\label{exa_Fock_space}
Let \(H\) be a finite-dimensional Hilbert space. Define the \emph{full Fock space} \(F(H)\) by
\[
F(H) \;=\; \mathbb{C}\,\Omega 
          \;\oplus\; H 
          \;\oplus\; (H \otimes H)
          \;\oplus\; (H \otimes H \otimes H) 
          \;\oplus\;\cdots,
\]
where \(\Omega\) is the \emph{vacuum vector}. For each \(v \in H\), define the 
\emph{left creation} and \emph{annihilation} operators \(a^*(v)\) and \(a(v)\) 
on \(F(H)\) by linearity and the rules
\[
a^*(v)\, \Omega = v,  \quad 
a^*(v)\,(x_1 \otimes x_2 \otimes \cdots \otimes x_n)
  \;=\; v \otimes x_1 \otimes x_2 \otimes \cdots \otimes x_n,
\]
\[
a(v)\,\Omega \;=\; 0,
\quad
a(v)\,(x_1 \otimes x_2 \otimes \cdots \otimes x_n)
  \;=\; \langle v,x_1\rangle\, (x_2 \otimes \cdots \otimes x_n).
\]
The algebra \(\mathcal{A}(F(H))\) is generated by these creation and 
annihilation operators. A natural state on this algebra is
\[
\phi(X) \;=\; \langle \Omega,\,X\,\Omega\rangle.
\]
This state is \emph{not} faithful (for instance, 
\(\phi(a^*(v)\,a(v))=0\) but \(a^*(v)\,a(v)\neq 0\)), and it is not tracial 
(\(\phi(a(v)\,a^*(v))=1\neq 0=\phi(a^*(v)\,a(v))\)).
\end{exa}

Overall, these definitions provide a framework to treat ``expectations'' of 
non-commuting operators in a way that parallels classical probability theory, 
but with genuinely new phenomena arising from non-commutativity.

\section{Distributions}

Suppose that \(X_1, X_2, \ldots, X_n\) are elements of a non-commutative probability space \((\mathcal{A},\phi)\). We will refer to these as random variables. Their \emph{distribution}\index{distribution|textbf} is defined as the linear map from the algebra of polynomials in non-commuting variables \(x_1, \ldots, x_n\) to \(\mathbb{C}\) given by
\[
f(x_1,\ldots,x_n) \mapsto \phi\bigl[f(X_1,\ldots,X_n)\bigr].
\]
Similarly, the \(\ast\)\emph{-distribution} is defined for polynomials in the non-commuting variables \(x_1,\ldots,x_n,y_1,\ldots,y_n\) by
\[
f(x_1,\ldots,x_n,y_1,\ldots,y_n) \mapsto \phi\bigl[f(X_1,\ldots,X_n,X_1^\ast,\ldots,X_n^\ast)\bigr].
\]

In other words, the distribution of a family of random variables is entirely determined by their joint moments.

We write \(X\cong Y\) and say that \(X\) is \emph{equivalent}\index{equivalent random variables} to \(Y\) if these two families of random variables have the same \(\ast\)-distribution.

\index{convergence in distribution|textbf}
\begin{defi}\label{definition_convergence_distribution}
A sequence of \(n\)-tuples of random variables \(\bigl(X_1^{(i)}, \ldots, X_n^{(i)}\bigr)\), \(i=1,2,\ldots,\) is said to \emph{converge in distribution} to \((X_1,\ldots,X_n)\) if for every non-commutative polynomial \(f\) in \(n\) variables one has
\[
\phi\bigl[f\bigl(X_1^{(i)},\ldots,X_n^{(i)}\bigr)\bigr] \longrightarrow \phi\bigl[f(X_1,\ldots,X_n)\bigr]
\]
as \(i\to\infty\).
\end{defi}

Convergence in \(\ast\)-distribution is defined similarly.

\bigskip

If we have a single self-adjoint random variable, then its distribution can be described by a probability measure.

\begin{propo}\label{proposition_spectral_distribution}
Suppose that \(X\) is a bounded self-adjoint element of a non-commutative probability space \((\mathcal{A},\phi)\). Then there exists a probability measure \(\mu\) on \(\mathbb{R}\) such that
\[
\phi\bigl(X^k\bigr) = \int_{\mathbb{R}} x^k\, \mu(dx)
\]
for all \(k\ge 0\).
\end{propo}

\textbf{Proof:} By the spectral theorem, one may write
\[
X = \int_{\mathbb{R}} x\, P(dx),
\]
where \(P\) is the projection-valued measure associated with \(X\). Define \(\mu(A) := \phi(P(A))\) for any Borel set \(A\subset \mathbb{R}\). It is straightforward to check that \(\mu\) is a probability measure. Then,
\[
X^k = \int_{\mathbb{R}} x^k\, P(dx),
\]
and taking the expectation yields
\[
\phi\bigl(X^k\bigr) = \int_{\mathbb{R}} x^k\, \mu(dx).
\]
\(\Box\)

\bigskip

This proposition extends to normal elements \(X\) (i.e., operators satisfying \(X^\ast X = XX^\ast\)). In that case, the measure \(\mu\) is defined on \(\mathbb{C}\), and for every polynomial \(P\) in two variables we have
\[
\phi\bigl(P(X,X^\ast)\bigr) = \int_{\mathbb{C}} P(z,\overline{z})\, \mu(dz).
\]

The measure \(\mu\) defined in Proposition~\ref{proposition_spectral_distribution} is called the \emph{probability distribution}\index{probability distribution of non-commutative r.v|textbf} of non-commutative random variable \(X\).

If \(X\) is not normal, then a different notion---the Brown measure---can be used to associate it with a probability measure. However, if we have a collection \(X_1,\ldots,X_n\) with \(n\ge 2\) and the operators do not commute, there is in general no way to relate their joint distribution to a single probability measure.

If one knows the moments of a self-adjoint random variable and wishes to recover the corresponding probability distribution, one can use the Cauchy transform method.

The \emph{Cauchy transform} of a probability distribution \(\mu\)\index{Cauchy transform!of a probability distribution} is defined by
\begin{equation}
\label{defi_Cauchy_transform}
G(z) = \int_{\mathbb{R}} \frac{\mu(dt)}{z-t},
\end{equation}
where \(z\) is a complex variable in the upper half-plane \(\mathbb{C}^+ = \{z\in\mathbb{C}\mid \operatorname{Im}z>0\}\). Another frequently used name for this object (or sometimes for its negation) is the \emph{Stieltjes transform}\index{Stieltjes transform|textbf}. We will use these names as synonyms.

If \(\mu\) is the spectral probability distribution of a random variable \(X\), then
\begin{equation}\label{formula_Cauchy_transform}
G(z) = \phi\Bigl(\frac{1}{z-X}\Bigr).
\end{equation}
Thus, we call \(G(z)\) defined in \eqref{formula_Cauchy_transform} the \emph{Cauchy transform of the random variable} \(X\)\index{Cauchy transform!of a random variable}. Since \(X\) is bounded, one may expand \((z-X)^{-1}\) as a series in \(z^{-1}\):
\[
G(z) = z^{-1} + \sum_{k=1}^\infty \phi\bigl(X^k\bigr) z^{-k-1},
\]
and this series converges for \(|z|\ge \|X\|\).

Once the Cauchy transform is known, the probability distribution can be recovered via the \emph{Stieltjes inversion formula}\index{Stieltjes inversion formula|textbf}. Namely,
\begin{equation}\label{formula_Stieltjes_inversion}
\mu(B) = -\frac{1}{\pi} \lim_{\varepsilon\downarrow0} \int_B \operatorname{Im} G(x+i\varepsilon)\,dx,
\end{equation}
provided that \(B\) is a Borel set with \(\mu(\partial B)=0\).

In particular, if \(G(z)\) admits an analytic continuation to a point \(x\in\mathbb{R}\), then \(\mu\) is absolutely continuous at \(x\) with density
\[
p(x) = -\frac{1}{\pi}\operatorname{Im}G(x).
\]

\bigskip

\textbf{Examples.}

\begin{exa}
If \(X\) is a Hermitian matrix and the state is given by the normalized trace, then the spectral probability distribution is
\[
\mu_X = \frac{1}{n}\sum_{i=1}^{n}\delta_{\lambda_i},
\]
where \(\lambda_i\) are the eigenvalues of \(X\) (counted with multiplicity).
\end{exa}

\begin{exe}
Suppose that \(X\) is a random Hermitian matrix considered as an element of the non-commutative probability space from Example~\ref{example_random_matrices} on page~\pageref{example_random_matrices}. What is its spectral probability distribution?
\end{exe}

\index{Haar unitary|textbf}
\begin{exa}[Haar unitary r.v.]
Consider the probability space given by the group algebra \(\mathbb{C}G\) as defined in Example~\ref{example_group_algebra} on page~\pageref{example_group_algebra}. In particular, let \(G\) be isomorphic to \(\mathbb{Z}\) and let \(g\) be its generator. Then, \(g^\ast = g^{-1}\) and hence \(g\) is unitary. One easily verifies that
\begin{equation}\label{equation_Haar_unitary}
\phi\bigl(g^k\bigr) =
\begin{cases}
1, & \text{if } k=0,\\[1mm]
0, & \text{if } k\neq 0.
\end{cases}
\end{equation}
Any unitary random variable with moments as in \eqref{equation_Haar_unitary} is called a \emph{Haar unitary}.

The probability distribution of \(g\) is the uniform (i.e., Haar) measure on the unit circle. In fact,
\[
\phi\bigl(g^k (g^\ast)^l\bigr) = \delta_{kl},
\]
and one has
\[
\int_{\mathbb{C}} z^k\overline{z}^l\,\mu(dz) = \frac{1}{2\pi}\int_0^{2\pi} e^{i(k-l)\theta}\,d\theta = \delta_{kl}.
\]
\end{exa}

\index{arcsine distribution|textbf}
\begin{exa}[Arcsine distribution]
Now consider the same probability space as in the previous example but let
\[
X = g + g^{-1}.
\]
This is a self-adjoint operator, and its moments are given by
\begin{equation}\label{moments_arcsine}
\phi\bigl(X^k\bigr) =
\begin{cases}
0, & \text{if } k \text{ is odd},\\[1mm]
\binom{k}{k/2}, & \text{if } k \text{ is even}.
\end{cases}
\end{equation}
\end{exa}

What is its probability distribution?

\begin{figure}[tbph]
\begin{center}
\includegraphics[width=6cm]{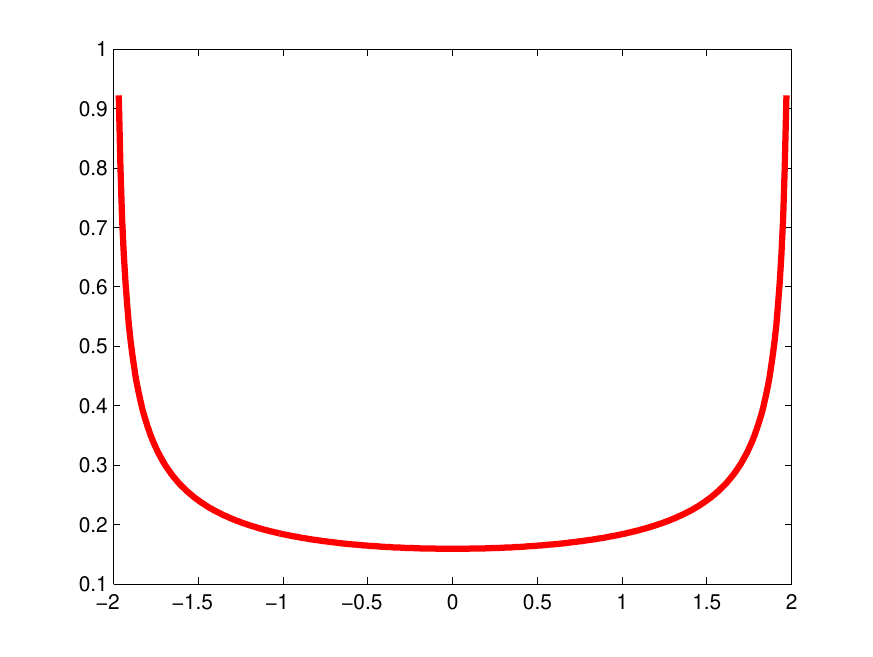}
\end{center}
\caption{The density of the arcsine distribution.}
\label{figure_arcsine}
\end{figure}

\begin{exe}\label{arcsine_distribution}
Check that if \(X = g+g^{-1}\), then the Cauchy transform is given by
\[
G(z) = \frac{1}{\sqrt{z^2-4}},
\]
and conclude that the spectral distribution of \(X\) has density
\[
p(x) =
\begin{cases}
\displaystyle \frac{1}{\pi}\frac{1}{\sqrt{4-x^2}}, & \text{if } |x|\le2,\\[1mm]
0, & \text{if } |x|>2.
\end{cases}
\]
\end{exe}

\index{distribution!arcsine}
This probability distribution is called the \emph{arcsine distribution}. Its density is illustrated in Figure~\ref{figure_arcsine}.

\index{semicircle r.v.|textbf}

\index{Fock space}
\index{Toeplitz--Cuntz algebra}
\index{Cuntz algebra}
\begin{exa}[Semicircle r.v.]\label{exa_semicircle}
Now consider the full Fock space \(F(H)\) where \(H\) is a one-dimensional Hilbert space with a unit vector \(e\). For brevity, denote by \(a\) the creation operator \(a(e)\) and by \(a^\ast\) the annihilation operator \(a^\ast(e)\). Also, define
\[
e_n := \begin{cases}
\underbrace{e\otimes\cdots\otimes e}_{n\text{ times}}, & n\ge 1,\\[1mm]
\Omega, & n=0.
\end{cases}
\]
The algebra generated by \(a\) and \(a^\ast\) is called the Toeplitz algebra. Recall that we defined the expectation by
\[
\phi(X) = \langle e_0, X e_0 \rangle,
\]
i.e., the expectation of \(X\) is the \((0,0)\)-entry in the matrix representation of \(X\) with respect to the basis \(\{e_n\}\).

What is the spectral distribution of the operator \(a+a^\ast\) with respect to this expectation?
\end{exa}

\begin{figure}[tbph]
\begin{center}
\includegraphics[width=6cm]{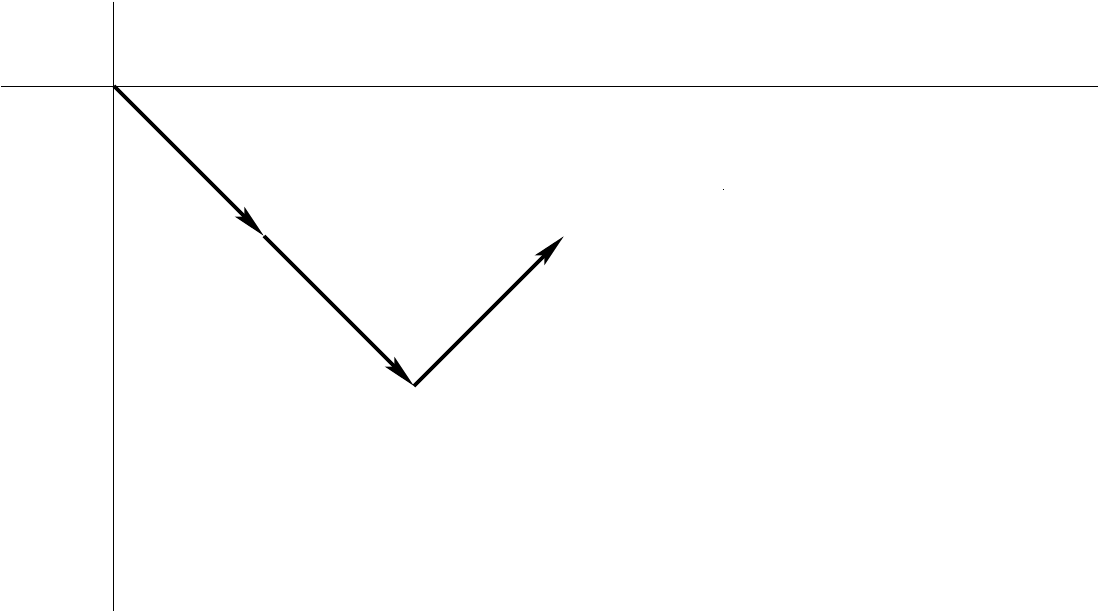}
\end{center}
\caption{Path corresponding to the term \(aa^\ast a^\ast\).}
\label{figure_walk1}
\end{figure}

As a first step, we calculate \(\phi\bigl((a+a^\ast)^n\bigr)\). To this end, we code the terms in the expansion of \((a+a^\ast)^n\) by certain paths on the \(\mathbb{Z}^2\) lattice. Starting at the origin (i.e., the point \((0,0)\)), at each step we move one unit to the right and either up or down, depending on whether the term \(a\) or \(a^\ast\) is chosen. Note that the terms in the expansion are read from right to left. For example, the term
\[
aa^\ast a^\ast
\]
in the expansion of \((a+a^\ast)^3\) corresponds to the path shown in Figure~\ref{figure_walk1}.

\bigskip

\begin{figure}[tbph]
\begin{center}
\includegraphics[width=6cm]{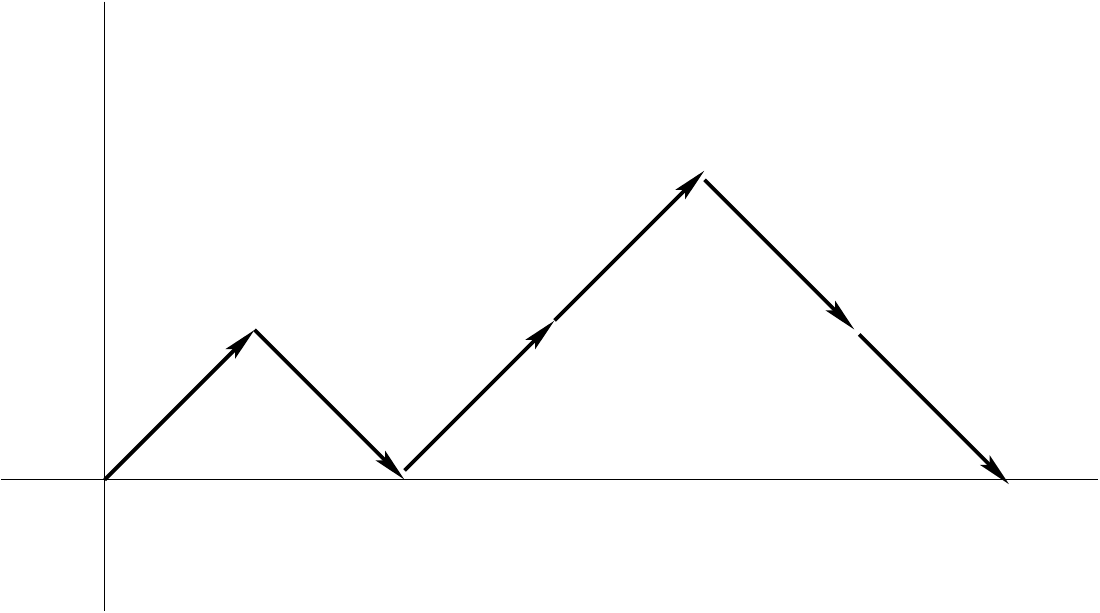}
\end{center}
\caption{Path corresponding to the term \(a^\ast a^\ast a a a^\ast a\).}
\label{figure_walk2}
\end{figure}

A term in the expansion of \((a+a^\ast)^n\) has zero expectation if and only if the corresponding path either goes below the horizontal axis or does not return to zero at the end (or both). The paths that remain on or above the horizontal axis and return to zero at the end are known as \emph{Dyck paths}\index{Dyck paths}. Hence, we must count the Dyck paths.

Clearly, the number of Dyck paths is zero when \(n\) is odd. When \(n=2k\) is even, the number of Dyck paths is given by the \emph{Catalan numbers}\index{Catalan numbers}. Specifically,
\[
C_k = \frac{1}{k+1}\binom{2k}{k} = \frac{(2k)!}{k!(k+1)!}.
\]

Indeed, the total number of paths from the origin to \((2k,0)\) is \(\binom{2k}{k}\) (since one must choose \(k\) upward steps among \(2k\) steps). To count those that drop below the horizontal axis, one uses the reflection principle to obtain a bijection with the paths from \((0,0)\) to \((2k,-2)\). The number of such paths is \(\binom{2k}{k-1}\). Therefore, the number of Dyck paths is
\[
\binom{2k}{k} - \binom{2k}{k-1} = \frac{1}{k+1}\binom{2k}{k} = C_k.
\]

The first few Catalan numbers are:
\[
\begin{array}{ccccccc}
k: & 0 & 1 & 2 & 3 & 4 & 5 \\
C_k: & 1 & 1 & 2 & 5 & 14 & 42
\end{array}
\]

A useful recursive formula for the Catalan numbers is:
\[
C_n = \sum_{k=0}^{n-1} C_{n-1-k}\, C_k.
\]

It is worth noting that Catalan numbers occur very often in enumerative combinatorics (see, e.g., \cite{stanley2015}). Two enumeration problems that naturally appear in this context are the counting of planar rooted trees and non-crossing pairings. Both are closely related to large random matrices and free probability.

\bigskip

We have now shown that the moments of \(a+a^\ast\) are given by the Catalan numbers. This allows us to recover its spectral probability distribution.

\bigskip
\begin{figure}[tbph]
\begin{center}
\includegraphics[width=6cm]{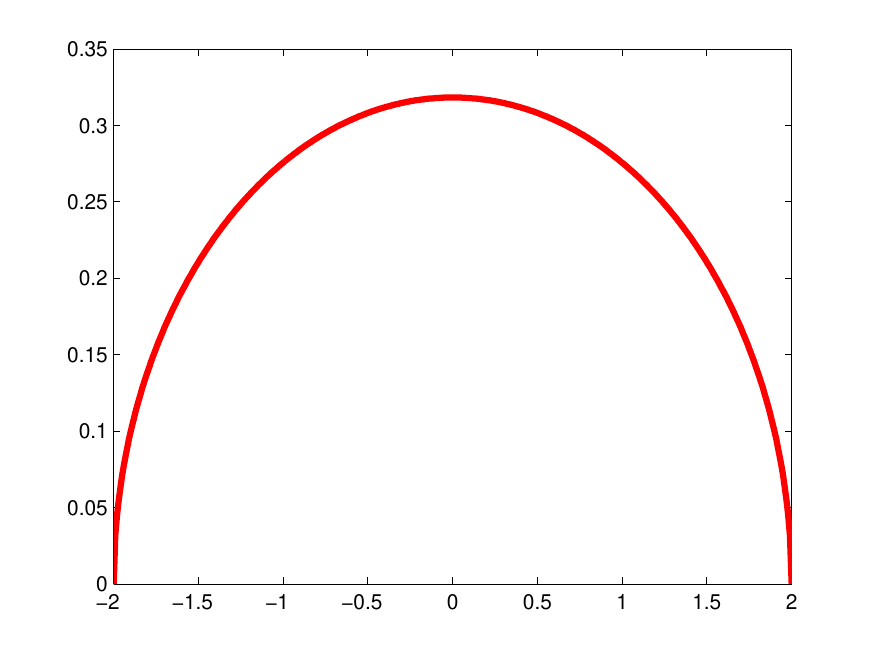}
\end{center}
\caption{The density of the semicircle distribution.}
\label{figure_semicircle}
\end{figure}

\begin{exe}\label{semicircle_distribution}
\label{exe_semicircle_distr}
Let \(a\) and \(a^\ast\) be the creation and annihilation operators defined above. Check that the Cauchy transform of \(a+a^\ast\) is given by
\[
G(z) = \frac{1}{2}\Bigl(z - \sqrt{z^2-4}\Bigr),
\]
and conclude that the spectral distribution of \(a+a^\ast\) has density
\[
p(x) =
\begin{cases}
\displaystyle \frac{1}{2\pi}\sqrt{4-x^2}, & \text{if } |x|\le2,\\[1mm]
0, & \text{if } |x|>2.
\end{cases}
\]
\end{exe}

This distribution is called the \emph{semicircle distribution}\index{distribution!semicircle}. Its density is illustrated in Figure~\ref{figure_semicircle}.

\section{Freeness}

We are often interested in computing expectations of products of several random variables. In classical probability, such computations are typically possible only if the terms in the product are independent. In the non-commutative setting, the condition of independence is replaced by the notion of \emph{freeness} (also called \emph{free independence}).

Recall first the definition of classical independence.

\index{independence!classical}
\begin{defi}
\label{definition_independence}
Let $\mathcal{A}_{1}$ and $\mathcal{A}_{2}$ be sub-\(\sigma\)-algebras of a (commutative) $\sigma$-algebra $\mathcal{A},$ and $E$ is the expectated value on $\mathcal{A}$. Sub-algebras $\mathcal{A}_{1}$ and $\mathcal{A}_{2}$ are called \emph{independent} if 
\[
E\bigl(a_{1}a_{2}\bigr) \;=\; E(a_{1}) \,E(a_{2})
\]
for all $a_{1}\in \mathcal{A}_{1}$ and $a_{2}\in \mathcal{A}_{2}.$
\end{defi}

Clearly this extends to finitely or countably many subalgebras in the usual way.

Here is the definition of free independence.

\index{independence!free}
\index{free independence|textbf}
\begin{defi}
\label{definition_freeness}
Let $(\mathcal{A}, \phi)$ be a noncommutative probability space, and let $\mathcal{A}_{1}, \mathcal{A}_{2}, \dots, \mathcal{A}_{n}$ be unital subalgebras of $\mathcal{A}$. They are called \emph{free} if 
\[
\phi\bigl(a_{1} \,a_{2}\,\cdots\, a_{s}\bigr) \;=\; 0
\]
for all choices $a_{1}\in \mathcal{A}_{i(1)}, \dots, a_{s}\in \mathcal{A}_{i(s)}$ such that 
\[
\phi(a_{k}) \;=\; 0 \quad\text{for each }k
\quad\text{and}\quad
i(k) \neq i(k+1)\ \text{for }1 \le k < s.
\]
\end{defi}

Random variables \(x_{1},\dots,x_{m}\) are called \emph{free} if the unital subalgebras they generate (together with their adjoints) are free.

\subsection*{Computing expectation of products of free variables.}  The freeness condition in Definition \ref{definition_freeness} in fact allows one to compute the expectation of any product of free random variables. Suppose \(a_{1},\dots,a_{s}\) is an \emph{alternating} sequence (i.e.\ the index of the subalgebra to which \(a_{k}\) belongs differs from that of \(a_{k+1}\)), but we do \emph{not} assume that \(\phi(a_{k})=0\). Then the expansion
\[
\phi\Bigl[(\,a_{1}-\phi(a_{1})\,)\,\cdots\,(a_{s}-\phi(a_{s}))\Bigr]
\;=\;0
\]
(by freeness) can be rearranged inductively to express 
\(\phi(a_{1}a_{2}\cdots a_{s})\) in terms of expectations of shorter products.

\begin{exe}\label{exa_prod_rule1}
For free elements \(x,y\), expanding via the above idea shows
\[
\phi(x\,y) \;=\; \phi(x)\,\phi(y)\,.
\]
\end{exe}

\begin{exe}\label{exa_prod_rule1a}
Similarly, if \((x_{1}, x_{2})\) is free from \(y\), then
\[
\phi(x_{1}\,y\,x_{2}) \;=\; \phi(x_{1}x_{2}) \,\phi(y)\,.
\]
\end{exe}

\begin{exe}\label{exe_prod_rule1b}
Let $a_1, \ldots, a_s$ be an alternating sequence of variables (i.e., $a_{i}$ and $a_{i + 1}$ are always from different sub-algebras), and suppose that $\phi(a_i) = 0$ for all $i \ne k$. Then 
\[
\phi(a_1 \ldots a_s ) \;=\; \phi(a_k) \phi(a_1 \ldots \wh {a_k} \ldots a_s),
\]
where $a_1 \ldots \wh {a_k} \ldots a_s$ denote the production $a_1 \ldots a_s$ with excluded term $a_k$. 
\end{exe}

These particular identities look identical to what we get for classically independent variables. However, the notion of freeness can differ dramatically from classical independence in more complicated products:

\begin{exa}
\label{exa_prod_rule2}
For free \(x,y\), expanding $\phi\big((x - \phi(x))\,(y - \phi(y))\,(x - \phi(x)) \,( y - \phi (y) \big)$ yields
\[
\phi(x\,y\,x\,y) \;=\;
\phi(x^2)\,\phi(y)^2 \,+\, \phi(y^2)\,\phi(x)^2 \;-\; \phi(x)^2\,\phi(y)^2,
\]
which can also be rewritten as
\begin{align*}
\phi(x\,y\,x\,y)
\;=\; &\phi\bigl(x^{2}\bigr)\,\phi\bigl(y^{2}\bigr)
\;-\;\Bigl[\phi\bigl(x^{2}\bigr)\;-\;\bigl(\phi(x)\bigr)^{2}\Bigr]
      \,\Bigl[\phi\bigl(y^{2}\bigr)\;-\;\bigl(\phi(y)\bigr)^{2}\Bigr] 
\\[6pt]
=\;  &\phi\bigl(x^{2}\bigr)\,\phi\bigl(y^{2}\bigr) \;-\;\mathrm{Var}(x)\,\mathrm{Var}(y)\,,
\end{align*}
where we define \(\mathrm{Var}(x) := \phi(x^{2}) - (\phi(x))^{2}\).
\end{exa}

\subsubsection*{Difference from classical independence.}\,
If \(x\) and \(y\) are classically independent random variables in a \emph{commutative} probability space, then
\(\phi(x\,y\,x\,y) = \phi(x^2)\,\phi(y^2)\).  Comparing with the free formula from Example \ref{exa_prod_rule2} shows that \(x\) and \(y\) can be both classically independent \emph{and} free only if
\(\mathrm{Var}(x) = 0\) or \(\mathrm{Var}(y)=0\). 
If \(\phi\) is faithful in a \(\mathrm{C}^*\)-algebraic sense, a zero variance implies the variable in question is a scalar multiple of the identity, so one of them is constant.

\subsubsection{Group algebra example of freeness.}\,
Suppose \(G_1\) and \(G_2\) are subgroups of a group \(G\).  Recall that \(G_1\) and \(G_2\) are called \emph{free} (as subgroups of \(G\)) if for every sequence \(g_1,\dots,g_s\) with \(g_k\neq e\), \(g_k\in G_{r(k)}\), and \(r(k)\neq r(k+1)\) for \(1\le k < s\), we have 
\[
g_1\,g_2\,\cdots\,g_s \;\neq\; e.
\]
Consider the group algebra \(\mathcal{A} = \mathbb{C}G\) with its standard state \(\phi\), given by \(\phi(e)=1\) on the identity and \(\phi(g)=0\) for \(g\neq e\).  Then \(\mathcal{A}_1=\mathbb{C}G_{1}\) and \(\mathcal{A}_2=\mathbb{C}G_{2}\) embed as unital subalgebras of \(\mathcal{A}\).  One shows:

\begin{propo}
\label{prop_groups_free_equiv}
The following are equivalent:
\begin{enumerate}
\item $G_1$ and $G_2$ are free subgroups of $G$.
\item $\mathcal{A}_1=\mathbb{C}G_{1}$ and $\mathcal{A}_2=\mathbb{C}G_{2}$ are free subalgebras of $\mathcal{A}=\mathbb{C}G$ with respect to the state \(\phi\).
\end{enumerate}
\end{propo}

\begin{proof}[Sketch of proof]
\quad
\textbf{\((1)\Rightarrow (2)\).}\,
Let \(a_i = \sum_j \alpha_{i,j}\,g_{i,j}\in \mathbb{C}G_{r(i)}\) (finite sums), with \(r(1)\neq r(2)\neq\cdots\neq r(s)\).  If \(\phi(a_{i})=0\) for each \(i\), this means none of the \(g_{i,j}\) appearing is the identity.  Then 
\(\phi(a_1\cdots a_s)\) is the coefficient of the identity in the product
\(\bigl(\sum_j \alpha_{1,j}g_{1,j}\bigr)\,\cdots\,\bigl(\sum_j \alpha_{s,j}g_{s,j}\bigr).\)
But by the freeness assumption on \(G_1\) and \(G_2\), there is no way to multiply a reduced word of alternating group elements to get \(e\), so that coefficient is zero.

\smallskip

\noindent
\textbf{\((2)\Rightarrow (1)\).}\,
If $G_1$ and $G_2$ were \emph{not} free subgroups of $G$, there is a reduced product \(g_1g_2\cdots g_s=e\) with each \(g_k\in G_{r(k)}\ne \{e\}\).  Then \(\phi(g_k)=0\) but \(\phi(g_1\cdots g_s)=1\), which would violate freeness in \(\mathbb{C}G\). 
\end{proof}

\noindent
In particular, if \(G=G_1 \ast G_2\) is the free product of two groups, then \(\mathbb{C}G_1\) and \(\mathbb{C}G_2\) are free subalgebras inside \(\mathbb{C}G\).

\section{Multivariate free CLT}

The goal of this section is to show that freeness is a very natural generalization of the independence concept for non-commutative variables. We do it by demonstrating that it implies an analogue of the classical Central Limit Theorem (``CLT'') for the sums of free random variables. Fortunately, for bounded non-commutative random variables this theorem can be proved without any advanced machinery,  by using simple combinatorics. In addition, the proof  will clearly show the special role played by non-crossing pairings. Since the proof is essentially the same both in univariate and multivariate setting, we prove the more general multivariate version.

\begin{figure}[tbph]
\begin{center}
\includegraphics{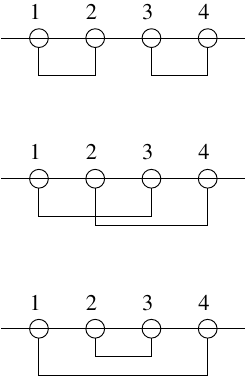}
\end{center}
\caption{Pairings of \(\{1,2,3,4\}\).}
\label{figure_pairings}
\end{figure}

Let \(\mathcal{P}_{2}(n)\) denote the set of pairings of the set \(\{1,2,\ldots,n\}\). (If \(n\) is odd, then this set is empty.) For example, \(\mathcal{P}_{2}(4)\) consists of three elements: \((12)(34)\), \((13)(24)\), and \((14)(23)\). These pairings are represented graphically in Figure~\ref{figure_pairings}.

\bigskip

\index{non-crossing pairing|textbf}
A pairing has a \emph{crossing} if there exist four elements \(i,j,k,l\) in \(\{1,2,\ldots,n\}\) with \(i<j<k<l\), such that \(i\) is paired with \(k\) and \(j\) is paired with \(l\). For example, the pairing \((13)(24)\) of \(\{1,2,3,4\}\) has a crossing. Let \(\mathcal{NCP}_{2}(n)\) be the set of \emph{non-crossing} pairings of \(\{1,2,\ldots,n\}\), that is, pairings that do not have a crossing. For example, \(\mathcal{NCP}_{2}(4)\) has two elements: \((12)(34)\) and \((14)(23)\).

\index{semicircle family|textbf}
\index{Wick's formula|textbf}
\begin{defi}
\label{defi_semicircular_family}
Let \(s_{i}\), \(i=1,\ldots,n\), be self-adjoint (centered) random variables such that \(E(s_{i}s_{j}) = c_{ij}\), and suppose higher moments are given by
\begin{equation}
\label{equ_free_Wick}
\phi \bigl(s_{i_{1}} s_{i_{2}} \cdots s_{i_{n}}\bigr)
\;=\;\sum_{\pi \in \mathcal{NCP}_2(n)} \;\prod_{(p,q)\in \pi} c_{\,i_{p} i_{q}}.
\end{equation}
Then the family \(\{s_{i}\}\) is called the \emph{semicircle family} with covariance \((c_{ij})\).
\end{defi}

Do such families exist? In other words, can we find variables \(s_{1},\ldots,s_{n}\) whose joint moments are given by the above formula? It turns out that the semicircle family exists if and only if the matrix \((c_{ij})\) is positive semidefinite. That is, 

\[
\sum_{i,j} c_{ij}\,z_i\,z_j \;\ge\; 0 \quad \text{for all vectors }(z_1,\dots,z_n).
\]

This is analogous to the classical case, where a multivariate Gaussian with covariance matrix \((c_{ij})\) exists if and only if \((c_{ij})\) is positive semidefinite. Formula \eqref{equ_free_Wick} is the free probability analogue of Wick’s formula for the moments of multivariate Gaussian r.v.s.

Let us postpone the proof of the existence of semicircle families and first prove the multivariate CLT.

\index{CLT!multivariate}
\begin{theo}\label{theorem_multivariate_CLT}
Let \(\bigl\{\,X_{i}^{(1)}, \ldots, X_{i}^{(r)}\bigr\}\) be a sequence of \(r\)-tuples of self-adjoint bounded random variables. Assume that:
\begin{enumerate}
\item the joint distribution of each \(r\)-tuple does not depend on \(i\),
\item \(\phi\bigl(X_{i}^{(\alpha)}\bigr) = 0\),
\item \(\phi\bigl(X_{i}^{(\alpha)} X_{i}^{(\beta)}\bigr) = c_{\alpha\beta}\),
\item the \(r\)-tuples are free for different \(i\).
\end{enumerate}
Define
\[
S_{n}^{(\alpha)} \;=\;\sum_{i=1}^{n} X_{i}^{(\alpha)},\quad \alpha=1,\dots,r.
\]
Then the \(r\)-tuple
\[
\Bigl\{\,S_{n}^{(1)}/\sqrt{n}, \ldots, S_{n}^{(r)}/\sqrt{n}\Bigr\}
\]
converges in distribution to the semicircle family with covariance \(\{c_{\alpha\beta}\}\).
\end{theo}

\begin{proof}
For simplicity, let us prove the theorem in the case \(r=2\), i.e.\ we have two components \(Y_i = X_i^{(1)}\) and \(Z_i = X_i^{(2)}\). Note that $Y_i$ and $Z_i$ are not necessarily free, however $\{Y_i, Z_i\}$ is free of $\{Y_j, Z_j\}$ if $i \ne j$. The general case follows by a similar argument. We examine a typical moment such as
\[
\phi\Bigl(S_{n}^{(1)} \,S_{n}^{(2)} \,S_{n}^{(1)} \,S_{n}^{(2)}\Bigr).
\]
Expanding the sums, we obtain
\[
\sum_{\,i_{1},\,i_{2},\,i_{3},\,i_{4}=1}^{n} \phi\bigl(Y_{i_{1}}\,Z_{i_{2}}\,Y_{i_{3}}\,Z_{i_{4}}\bigr).
\]
The value of each term depends only on how the indices \(i_{1},i_{2},i_{3},i_{4}\) are grouped among distinct tuples, rather than on their exact numeric values. For instance,
\[
\phi\bigl(Y_{100}\,Z_{50}\,Y_{50}\,Z_{100}\bigr)
\;=\;\phi\bigl(Y_{1}\,Z_{2}\,Y_{2}\,Z_{1}\bigr),
\]
but this might differ from
\[
\phi\bigl(Y_{1}\,Z_{1}\,Y_{1}\,Z_{1}\bigr).
\]
It is convenient to encode each pattern of indices by a partition. In our example, the partition corresponding to \(\phi(Y_{100}Z_{50}Y_{50}Z_{100})\) and \(\phi(Y_{1}Z_{2}Y_{2}Z_{1})\) is \([1,4],[2,3]\), whereas the partition for \(\phi(Y_{1}Z_{1}Y_{1}Z_{1})\) is \([1,2,3,4]\).

Any part of a partition that is a singleton (i.e.\ an index that appears exactly once) gives zero contribution, because we assumed that the variables have mean zero and are free from the remaining variables. Hence partitions with singletons do not contribute.

Next, we note that for a given moment, the leading contribution as \(n\to\infty\) comes from partitions with the largest number of blocks (i.e.\ pairings). Indeed, each block corresponds to a distinct index value in the product, so for large $n$ there are roughly \(n\) choices of index per block. Therefore, a partition with \(k\) blocks corresponds to roughly \(n^{k}\) summation terms, so pair partitions (which have \(r/2\) blocks for \(r\) factors) dominate in the large-\(n\) limit.

Now we evaluate contributions from pairings. Suppose a particular pairing \(\pi\) of \(\{1,\dots,r\}\) is given. If \emph{no} pair in \(\pi\) connects neighboring factors in the product, then the sequence of elements in the product is alternating meaning that $i_{p+1} \ne i_{p+1}$ for any $p$. Since $\phi(Y_i) = \phi(Z_i) = 0$ by assumption  then this product has zero expectation by freeness. (For example, $\phi(Y_1 Z_2 Y_1 Z_2) = 0$.)

On the other hand, if a pair \((p,q)\) connects neighboring positions, then we can factor out the corresponding expectation by one of the properties of $\phi$. For instance,
\[
\phi\bigl(Z_{3}\,Y_{1}\,\underline{Y_{2}\,Z_{2}}\,Z_{1}\,Z_{3}\bigr)
\;=\; \phi \bigl(Y_{2}\,Z_{2}\bigr)\, E\bigl(Z_{3}\,Y_{1}\,Z_{1}\,Z_{3}\bigr)
\;=\; c_{YZ} \;E\bigl(Z_{3}\,Y_{1}\,Z_{1}\,Z_{3}\bigr).
\]
We continue removing such pairs step by step. One sees that the only way this process \emph{fails} to connect neighboring factors at some point is if pairing \(\pi\) has a crossing: in that scenario, once certain pairs are removed, eventually no remaining pair connects neighbors, making that product’s expectation zero. Conversely, if \(\pi\) is non-crossing, we always can remove a neighboring pair, factoring out its covariance, until the entire product is determined by
\[
\prod_{(p,q)\in \pi} c_{\alpha(p)\,\alpha(q)},
\]
where \(\alpha(p)\in \{Y,Z\}\) indicates which variable appears at position \(p\). Moreover, each non-crossing pairing accounts for on the order of \(n^{r/2}\) terms in the original sum. Thus for a sequence $(i_1, \ldots, i_r)$ with $i_k \in \{1, 2\}$, we have 
\[
\phi \bigl(S_{n}^{(i_1)}\,S_{n}^{(i_2)}\,\dots\bigr)
\;\sim\; n^{r/2}\;\sum_{\pi \in \mathcal{NCP}_{2}(r)}
\prod_{(p,q)\,\in\,\pi} c_{\alpha(p)\,\alpha(q)}.
\]
After dividing by \(\sqrt{n}\) for each factor, we see that
\[
\bigl\{\,S_{n}^{(X)}/\sqrt{n},\;S_{n}^{(Y)}/\sqrt{n}\bigr\}
\]
converges to the semicircle family with covariances \(c_{XX}, c_{XY}, c_{YX}, c_{YY}\). This completes the proof of Theorem~\ref{theorem_multivariate_CLT}.
\end{proof}

Now let us address the question of existence of semicircle families. First, if $n = 1$, one can observe that the distribution from Exercise \ref{exe_semicircle_distr} has the moments defined in \eqref{equ_free_Wick} with $c_{11} = 1$. We can find a corresponding self-adjoint random variable, which we call a standard (univariate) semicircle random variable. This establishes existence of the semicircle family for the dimension $n = 1$. 

Then, we can apply Theorem ~\ref{theorem_multivariate_CLT} in the univariate setting and since the limit in Theorem ~\ref{theorem_multivariate_CLT} must be stable, we establish the following property of free standard semicircle random variables $s_1, \ldots, s_n$:
\begin{equation}
\label{semicircle_variable_stability}
s = (s_1 + \ldots + s_n)/\sqrt{n}
\end{equation}
is a standard semicircle. 

Now, let \(s=(s_{1},\ldots,s_{r})\) be an \(r\)-tuple of random variables, each of which has the standard semicircle distribution, and assume \(s_{1},\ldots,s_{r}\) are free. Such an \(r\)-tuple is called a \emph{semicircle system}. (Note the difference from the semicircle family: here we require that the components of the $r$-tuple are free from each other.) If \(s^{(1)},\ldots,s^{(n)}\) are \(n\) mutually free semicircle systems, then by using \eqref{semicircle_variable_stability} component-wise, we have that 
\begin{equation}
\label{semicircle_system_stability}
s \;=\;\frac{s^{(1)}+\cdots+s^{(n)}}{\sqrt{n}}
\end{equation}
is itself a semicircle system.

Now let \(C=(c_{ij})\) be a non-negative definite \(r\times r\) matrix. We can factor it as \(C = A^{\prime}A\), where \(A^{\prime}\) is the transpose of \(A\). Let \(x = A\,s\), where \(s\) is a semicircle system. It is then straightforward to verify that each component of \(x\) has the semicircle distribution, and the covariance of \(x_{i}\) and \(x_{j}\) is \(C_{ij}\).

(The linear transformation \(x = A\,s\) is the free-probability counterpart of the classical Gaussian construction \(x = A\,g\), where $g$ is the standard multivariate Gaussian r.v.)

If \(s^{(1)},\ldots,s^{(n)}\) are \(n\) free semicircle systems, define \(x^{(i)} = A\,s^{(i)}\). Then the \(r\)-tuples \(x^{(i)}\) are mutually free, all have the same joint distribution, and by \eqref{semicircle_system_stability} that distribution coincides with that of
\[
X^{(n)} \;=\;\frac{x^{(1)} + \cdots + x^{(n)}}{\sqrt{n}}.
\]
By the same counting argument as in Theorem~\ref{theorem_multivariate_CLT}, the moments of the limit of \(X^{(n)}\) are given by
\[
\sum_{\pi \in \mathcal{NCP}_{2}(n)} \,\prod_{(p,q)\in \pi} C_{\,i_{p} i_{q}}.
\]
Since the distribution of \(X^{(n)}\) does not depend on \(n\), each $X^{(n)}$ and therefore each \(x^{(i)}\) must share these same moments; hence each \(x^{(i)}\) is a semicircle family with covariance matrix \(C\). This settles the existence of semicircle families.

\section{Exercises}

\begin{exe}
Show that the algebra of constants \(\{\,cI: c\in \mathbb{C}\}\) is free from \emph{any} other unital subalgebra of \(\mathcal{A}\). 
\end{exe}

\begin{exe}
Let \(x\) and \(y\) be commuting operators in a noncommutative probability space and suppose they are free.  Use the rules in Examples~\ref{exa_prod_rule1} and \ref{exa_prod_rule2} to show that at least one of \(x\) or \(y\) must be a constant multiple of the identity.
\end{exe}

\begin{exe}
\label{exe_conjugation_by_unitary}
Let $(\AA,\phi)$ be a $\ast$-probability space. Consider a unital subalgebra $\BB \subset \AA$ and a Haar unitary $u \in \AA$ such that $\{u, u^\ast\}$ and $\BB$ are free. Show that then also $\BB$ and $u^\ast \BB u$ are free. Here the algebra $u^\ast B u$ is 
\[
u^\ast \BB u := \{u^\ast b u | b \in  \BB\} \subset A.
\]
\end{exe}

\index{Fock space}
\begin{exe}
\label{exe_Fock_space}
Let \(\xi\) and \(\eta\) be orthogonal vectors in a Hilbert space \(H\).  Consider the algebra
\(\mathcal{A}=\mathbb{C}\langle a(\xi),a^{\ast}(\xi),a(\eta),a^{\ast}(\eta)\rangle\)
of creation and annihilation operators on the free Fock space over \(H\).  Suppose \(\bigl(\mathcal{A},\phi\bigr)\) is endowed with the standard vacuum state \(\phi\). Define
\(\mathcal{A}_{1}=\mathbb{C}\langle a(\xi),a^{\ast}(\xi)\rangle\)
and
\(\mathcal{A}_{2}=\mathbb{C}\langle a(\eta),a^{\ast}(\eta)\rangle\).
Show that \(\mathcal{A}_{1}\) and \(\mathcal{A}_{2}\) are free.
\end{exe}

\section{Notes}
The concept of freeness in operator algebras, along with its connection to independence in probability theory, was introduced in \cite{voiculescu83, voiculescu86}. A standard textbook and reference in free probability theory is \cite{nica_speicher06}. In particular, Exercise \ref{exe_conjugation_by_unitary} corresponds to Exercise 5.24 in this reference.

\chapter{Free Cumulants}

\section{Motivation and definition}

Cumulants were introduced by the Danish scientist T.\,N.\,Thiele at the end of the
19th century. They proved very useful for analyzing nonlinear transformations
of random variables. In particular, Brillinger \cite{brillinger80} used them
as a main technical tool in the spectral analysis of time series.

Why are cumulants useful? Let us fix $n$ random variables
$X_{1},X_{2},\ldots ,X_{n}.$ Their dependence structure can be explored
by examining the moments $\E\bigl(X_{i_{1}}X_{i_{2}}\cdots X_{i_{n}}\bigr).$
However, raw moments alone can be somewhat opaque. As a simple example, even if
$X_{1}$ and $X_{2}$ are independent, $\E\bigl(X_{1}X_{2}\bigr)$ is not necessarily
zero (unless at least one of them has mean zero). A better strategy is to isolate the
\emph{new} or \emph{higher-order} interaction among $X_{1},\dots,X_{n}$ by
subtracting contributions that arise from lower-order interactions.

For instance, in the classical theory, the first nontrivial interaction term
between $X_{1}$ and $X_{2}$ is the covariance:
\[
\Cov(X_1, X_2) \;=\; \E(X_1 X_2)\;-\;\E(X_1)\,\E(X_2).
\]

\index{cumulants!classical}
Classical cumulants generalize this principle inductively. One sets the
first-order cumulants to be $c_1(X_i)=\E(X_i)$ and then writes
\begin{equation}
\label{classical_cumulants}
\E \bigl(X_1 \cdots X_n\bigr) \;=\; \sum_{\pi \in \PP(n)} c_\pi(X_1, \ldots, X_n),
\end{equation}
where the sum is over all partitions $\pi$ of the set $[n] = \{1, \ldots, n\}$,
and $c_\pi(X_1,\ldots,X_n)$ is the product of the corresponding elementary cumulants
$c_k(\cdot)$ over the blocks of the partition. For example,
\[
c_{\{1,3\},\{2\}}(X_1,X_2,X_3) \;=\; c_2(X_1, X_3)\,\cdot\, c_1(X_2).
\]
In this manner, one can solve recursively for the elementary cumulants
$c_k$. For $n=2$, \eqref{classical_cumulants} implies
\[
\E\bigl(X_1 X_2\bigr)
\;=\;
c_2(X_1,X_2)
\;+\;
c_1(X_1)\,c_1(X_2).
\]
Hence
\[
c_2(X_1, X_2) \;=\; \E(X_1X_2)\;-\; \E(X_1)\,\E(X_2)
\;=\; \Cov(X_1,X_2).
\]

A convenient lattice-theoretic perspective on cumulants was given by
T.~Speed \cite{speed83}, building on Rota's theory of M\"{o}bius inversion
on lattices in \cite{rota64}.

\index{cumulants!free}
\index{free cumulants}
Speed's theory assumes commutative random variables, so it employs set
partitions of $\{1,2,\ldots,n\}$ (an \emph{unordered} set) in the above
formula. It does not apply \emph{as is} to noncommutative random variables.
However, Roland Speicher extended this idea to \emph{free} random variables
by changing the relevant lattice from all partitions to the
\emph{noncrossing} partitions of the \emph{ordered} set $[n]$. Concretely,
the free version of \eqref{classical_cumulants} becomes
\begin{equation}
\label{free_cumulants}
\phi\bigl(X_1 \cdots X_n\bigr)
\;=\;
\sum_{\pi \in NCP(n)} \kappa_\pi(X_1, \ldots, X_n),
\end{equation}
where $NCP(n)$ denotes the set of \emph{noncrossing} partitions $\pi$ of
the ordered set $[n]$, and where $\kappa_\pi$ is the product of the
elementary free cumulants $\kappa_k$ over the blocks of $\pi$.

One must check that the recipe \eqref{free_cumulants} is consistent and can be used to calculate free cumulants. Fortunately, M\"{o}bius inversion on the noncrossing partition lattice ensures that it is.

A key and beautiful result due to Speicher states that if any two of the
variables in $\kappa_n(X_1,\ldots,X_n)$ belong to different \emph{free}
subalgebras, then that free cumulant vanishes. Thus, just as in the
classical case (where mixed cumulants vanish under classical independence),
here mixed free cumulants vanish under free independence. This vanishing of
mixed cumulants enables a variety of explicit calculations of moments and
leads, for instance, to the straightforward proof of the Free Central Limit Theorem.

\section{Lattices and the Möbius inversion}

\index{poset|textbf}

Let us say some words about the Möbius inversion. 
Recall that a \emph{poset} (a partially ordered set) is a set with an
order relation $\leq$ which is defined only for some pairs of elements of the set. The
order relation must satisfy:

\begin{enumerate}
\item \textbf{Reflexivity:} $a \leq a$ for all $a$.  
\item \textbf{Transitivity:} if $a \leq b$ and $b \leq c$, then $a \leq c$.  
\item \textbf{Antisymmetry:} if $a \leq b$ and $b \leq a,$ then $a = b$.
\end{enumerate}

(We will write $a < b$ if $a \leq b$ and $a \neq b$.)

\medskip

\noindent
\textbf{Examples of posets:}
\begin{enumerate}
\item Positive integers with respect to the divisibility relation.

\item Subsets of a finite set with respect to inclusion. (For instance, $\{1\} < \{1,2\}$.)

\item Partitions of a finite set with respect to refinement. (For instance, $\{\{1\}, \{2,3\}\} < \{\{1,2,3\}\}$.)

\item Linear subspaces of a finite-dimensional vector space with respect to
inclusion.

\item Hermitian matrices with the Loewner (semidefinite) order: $A \leq B$ if and only
if $B - A$ is positive semidefinite. 
\end{enumerate}

(In the case of the Loewner order, it is usually written that $A < B$ if $B - A$ is positive definite, which is different from our notation, since in our notation it means that $B - A$ is positive semidefinite and $A \ne B$.)

\subsection*{Lattices}
\index{lattice|textbf}
A \emph{lattice} is a poset with two additional properties. First,
for every two elements $a$ and $b$, there is a least \emph{upper bound} $c$
such that:
\begin{enumerate}
\item $a \leq c$ and $b \leq c$,
\item if $a \leq c'$ and $b \leq c'$, then $c \leq c'$.
\end{enumerate}
This element is called the \emph{join} of $a$ and $b$ and is denoted $a \vee b$.

Similarly, for every $a$ and $b$, there is a greatest \emph{lower bound} $c$
such that:
\begin{enumerate}
\item $c \leq a$ and $c \leq b$,
\item if $c' \leq a$ and $c' \leq b$, then $c' \leq c$.
\end{enumerate}
This element is called the \emph{meet} of $a$ and $b$ and is denoted $a \wedge b$.

\medskip

Note that all of the above examples except the last one are \emph{lattices}. For the set of all subsets of a finite set, where $\vee$ and $\wedge$ correspond to union and intersection, respectively. 

In the Hermitian matrix example, for two matrices $A$ and $B$, there may be many matrices 
$C$ such that $C - A$ and $C - B$ are positive semidefinite, and these $C$'s need not be comparable. In particular, the join $A \vee B$ does not exist in general.

As another example, note that the set of algebraic integers is partially ordered by the divisibility relation but is often not a lattice. One typically considers the poset of ideals instead, which does form a lattice.

\medskip

Now, let us consider a finite lattice which has a maximum and a minimum
element. The maximum element, denoted by $1$, is greater than any other element
in the lattice, and the minimum element, $0$, is smaller than any other
element. 

Let us take a function $f(x)$ defined on the elements of the lattice and define
\begin{equation}\label{directFormula}
g(y) \;=\; \sum_{x \leq y} f(x).
\end{equation}%
We wish to invert this relation, that is, to find a formula expressing $f(x)$ in terms of $g(y)$. We look for a function $\mu(x,y)$ such that  
\begin{equation}\label{inversionFormula}
f(x)\;=\;\sum_{y \leq x} g(y)\,\mu(y,x).
\end{equation}%
The function $\mu(x,y)$ is called the \emph{Möbius function} of the lattice.

\index{M\"obius function|textbf}

\subsection*{Calculation of the Möbius function}  
The Möbius function can be calculated by the following recursion:
\begin{enumerate}
\item $\mu(x,x) = 1$ for all $x$.
\item For $x < y$, define
\[
  \mu(x,y) \;=\; - \sum_{\,x \leq z < y\,}\mu(x,z).
\]
\end{enumerate}
Equivalently, for every fixed pair $x \leq y$,
\begin{equation}
\label{mu-sum-identity}
\sum_{x \leq z \leq y} \mu(x,z)
\;=\;
\begin{cases}
1, & \text{if } x = y, \\[6pt]
0, & \text{if } x \neq y.
\end{cases}
\end{equation}

\subsection*{Rota's inversion formula}

\index{Rota's inversion formula|textbf}
\begin{theo}[Rota's Möbius Inversion]
If $\mu(x,y)$ is defined as above, and $f$ and $g$ satisfy 
\[ g(y) \;=\; \sum_{x \leq y} f(x), \]
then the inversion formula
\[ f(x) \;=\; \sum_{y \leq x} g(y)\,\mu(y,x) \]
holds.
\end{theo}

\begin{proof}
For a fixed $x$, we compute:
\begin{align*}
\sum_{y \leq x} g(y)\,\mu(y,x) 
&= \sum_{y \leq x} \Biggl(\sum_{z \leq y} f(z)\Biggr)\,\mu(y,x) \\[6pt]
&= \sum_{z \leq x} f(z)\,\sum_{\,z \leq y \leq x\,}\mu(y,x) \\[6pt]
&= f(x),
\end{align*}
where the last step uses the identity \eqref{mu-sum-identity}.
\end{proof}

The Möbius Inversion can be generalized to some infinite posets provided that at least the assumption of local finiteness holds.

One also has a generalized version of Rota's formula, which will be useful later. 

\index{Rota's inversion formula!generalized}
\begin{propo}[Rota-type inversion]
\label{inf_formula2}
Suppose that in a finite poset we have
\[
   g(y) \;=\; \sum_{\,x \,\le\, y}\; f(x).
\]
Then for every \(a \le z\) we have
\begin{equation}
\label{equ_Rota_inversion3}
   \sum_{\,x \,\vee\, a \;=\; z}\; f(x)
   \;=\;
   \sum_{\,a \,\le\, y \,\le\, z}\; g(y)\,\mu(y,z).
\end{equation}
\end{propo}

\begin{proof}
We compute:
\[
\begin{aligned}
   \sum_{\,a \,\le\, y \,\le\, z}\; g(y)\,\mu(y,z)
   &=\;
   \sum_{\,a \,\le\, y \,\le\, z}
   \,\Bigl(\sum_{\,x \,\le\, y}\; f(x)\Bigr)\,\mu(y,z)
\\
   &=\;
   \sum_{x} \; f(x)
   \,\sum_{\substack{a \,\le\, y \,\le\, z\\ x \,\le\, y}}
   \;\mu(y,z)
   \;=\;
   \sum_{x} \; f(x)
   \,\sum_{\substack{a \,\vee\, x \,\le\, y \,\le\, z}}
   \;\mu(y,z).
\end{aligned}
\]
Since \(\sum_{\,w \,\le\, y \,\le\, z}\mu(y,z)=0\) unless \(w=z\), this
inner sum is nonzero precisely when \(a\vee x = z\).  Thus the right
side becomes
\[
   \sum_{\,x \,\vee\, a \;=\; z}\; f(x),
\]
as claimed.
\end{proof}

\begin{coro}
\label{coro_vee_a_equals_one}
For every \(a \neq 0\) in our lattice,
\[
   \sum_{\,x \,\vee\, a = 1}\; \mu(0,x)
   \;=\;
   0,
\]
and hence
\begin{equation}
\label{mobius01}
   \mu(0,1)
   \;=\;
   -\sum_{\substack{x \,\vee\, a = 1\\ x \neq 1}}
   \mu(0,x)\,.
\end{equation}
\end{coro}

\begin{proof}
Apply Proposition~\ref{inf_formula2} with \(f(x)=\mu(0,x)\).  We have
\[
g(y) = \sum_{\,x \,\le\, y}\; \mu(0, x) = \delta_0(y)
\]
by \eqref{mu-sum-identity}, and therefore,
\[
   \sum_{x \vee a = 1} \mu(0,x) =
   \sum_{a \le  y  \le  z} \delta_0 (y)\,\mu(y, z) = 0,
 \] 
 for $a \ne 0$.  And \eqref{mobius01} is just a rearrangement of this sum.
\end{proof}

\section{The lattice of non-crossing partitions}
Terry Speed’s analysis of cumulants is based on the M\"{o}bius inversion
for the lattice of all partitions of the (unordered) set $\{1, \ldots, n\}$.
Roland Speicher’s construction of free cumulants uses instead the lattice of
\emph{non-crossing} partitions of the (ordered) set $\{1, \ldots, n\}$.
Most properties of this lattice were first discovered by Kreweras.

\index{non-crossing partitions|textbf}
\begin{defi}[Non-crossing partitions]
Given the ordered set $\Omega = \{1,2,\ldots,n\}$, a partition $\pi$ of
$\Omega$ is said to have a \emph{crossing} if there exist indices
$i<j<k<l$ such that
\begin{enumerate}
\item $i$ and $k$ are in the same block of $\pi$, and
\item $j$ and $l$ are in another block of $\pi$.
\end{enumerate}
When no such quadruple $(i,j,k,l)$ exists, we call $\pi$ a \emph{non-crossing
partition}.
\end{defi}
 In other words, if you plot the points $1,2,\ldots,n$ on a
horizontal line and connect the points belonging to each block of $\pi$
by curves above the line, those curves should never intersect.

For example, the partition $\bigl(\{1,3\}, \{2,4\}\bigr)$ of $\{1,2,3,4\}$
\emph{does} have a crossing (since $1<2<3<4$, the block $\{1,3\}$ crosses
the block $\{2,4\}$), whereas
$\bigl(\{1,4\},\{2,3\}\bigr)$ does \emph{not} have a crossing.

\index{non-crossing partitions!lattice structure}
\subsection*{Lattice structure.}
The set of all non-crossing partitions of $\{1,\ldots,n\}$ forms a lattice
under the \emph{refinement} order. Recall that $\alpha \le \beta$
means that $\alpha$ is a refinement of $\beta$, i.e.\ each block of $\alpha$
is contained in a block of $\beta$. We denote this lattice by $NC(n)$.

\begin{figure}[tbph]
\begin{center}
\includegraphics[width=12cm]{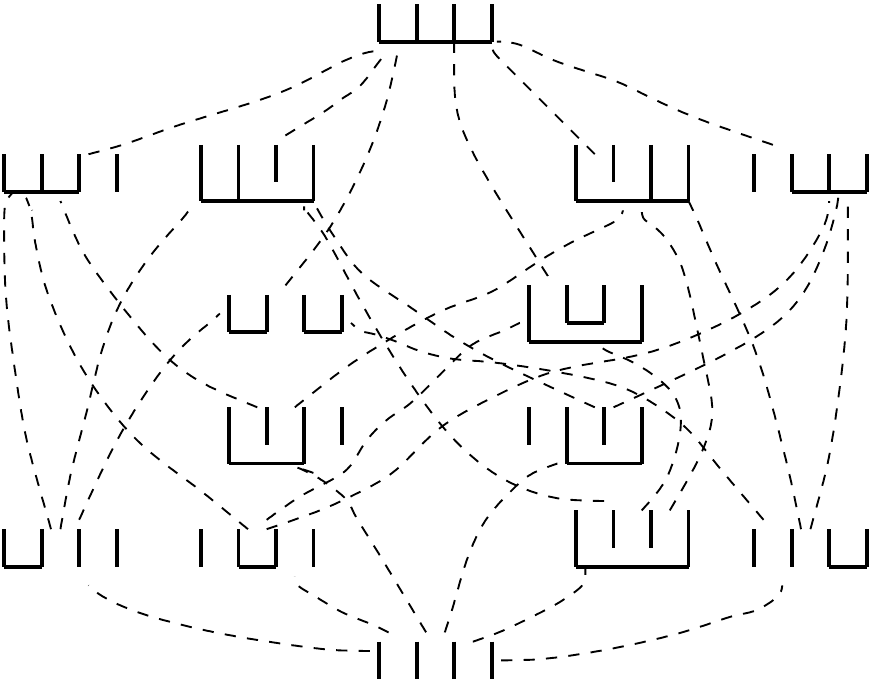}
\end{center}
\caption{$NC(4)$, the non-crossing partitions of the ordered set
$\{1,2,3,4\}$. Smaller elements appear lower in the diagram. The solid lines
indicate which points belong to each block, and the dashed lines indicate
the partial order.}
\label{figure_noncrossing_partitions}
\end{figure}

\subsection*{Kreweras complement}
\index{Kreweras complement|textbf}
An important property of $NC(n)$ is the existence of the \emph{Kreweras
complement} operation. Let $\pi$ be a non-crossing partition of
$\{1,\ldots,n\}$. We introduce an enlarged ordered set
$\{1,\overline{1},2,\overline{2},\ldots,n,\overline{n}\}$.
The Kreweras complement $K(\pi)$ is defined to be the \emph{largest}
non-crossing partition (in the refinement sense) of
$\{\overline{1},\ldots,\overline{n}\}$ such that the union of the blocks
of $\pi$ (on the unbarred elements) and the blocks of $K(\pi)$
(on the barred elements) forms a non-crossing partition of the entire
double set $\{1,\overline{1},\ldots,n,\overline{n}\}$.

\begin{figure}[tbph]
\begin{center}
\includegraphics[width=8cm]{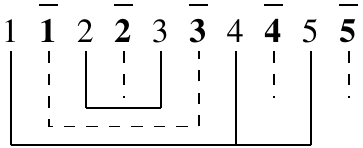}
\end{center}
\caption{An example of the Kreweras complement. The original partition $\pi$
is shown by solid lines, and its complement $K(\pi)$ by dashed lines.}
\label{figure_Kreweras_complement}
\end{figure}

Figure~\ref{figure_Kreweras_complement} shows an example. The solid arcs
depict $\pi = (145)(23)$, and the dashed arcs depict its Kreweras complement
$K(\pi) = (\overline{1}\,\overline{3})(\overline{2})(\overline{4})(\overline{5})$.

The Kreweras complement induces an \emph{anti-isomorphism} of $NC(n)$:
if $\lambda \le \pi$, then $K(\lambda) \ge K(\pi)$. In particular,
$K(0) = 1$ and $K(1) = 0$, where $0$ and $1$ are the minimal and maximal
partitions, respectively. Also, $K\bigl(K(\pi)\bigr)$ is essentially
a “shift” of $\pi$ rather than $\pi$ itself, due to the labeling
$\{\overline{1},\ldots,\overline{n}\}$. One can fix this by considering
the points $1,\ldots,n$ on a circle and suitably modifying
the definition.

\section{Definition of free cumulants}

Fix $n$ noncommutative variables $X_1,X_2,\ldots,X_n$. Let
$\lambda =(\lambda^1,\ldots,\lambda^k)$ be a partition of the ordered set
$\{1,\ldots,n\}$. This means that each $\lambda^i$ is a (consecutive)
ordered subsequence of $\{1,\ldots,n\}$. Concretely, we can write
\[
\lambda^i \;=\; \bigl[\lambda_1^i,\,\lambda_2^i,\ldots,\lambda_{|\lambda^i|}^i\bigr]
\]
with
\(\lambda_1^i < \lambda_2^i < \cdots < \lambda_{|\lambda^i|}^i\)
for each $i$.

\index{moments!generalized|textbf}
For such a partition $\lambda$, define the
\emph{generalized moment} of $(X_1,\ldots,X_n)$ by
\begin{equation}
\label{defi_gen_moments}
\phi_\lambda(X_1,\ldots,X_n)
\;:=\;
\prod_{i=1}^k
\phi\!\bigl(X_{\lambda_1^i}\,X_{\lambda_2^i}\cdots X_{\lambda_{|\lambda^i|}^i}\bigr).
\end{equation}
That is, for each block of the partition (in order), you take the product
of the corresponding variables and evaluate $\phi(\,\cdot\,)$ (the linear
functional), then multiply these blockwise results together. For instance,
if $\lambda = \bigl(\{1,4\},\{2,3\}\bigr)$, then
\[
\phi_\lambda(X_1,X_2,X_3,X_4)
\;=\;
\phi(X_1 X_4)\;\phi(X_2 X_3).
\]
Since the variables $X_1,\ldots,X_n$ will be understood from context, we
often abbreviate this as $\phi_\lambda$.

\index{free cumulants|textbf}

\begin{defi}
The \emph{free cumulants} of the variables $X_1,\ldots,X_n$ are obtained
by applying Rota’s M\"{o}bius inversion to the generalized moments on
the lattice $NC(n)$. 
\end{defi}
Concretely, we have
\[
\phi_\sigma \;=\; \sum_{\lambda \,\le\, \sigma} k_\lambda
\qquad\text{for all }\sigma \in NC(n),
\]
where $k_\lambda$ are the free cumulants we want to define. Equivalently,
\[
k_\sigma
\;=\;
\sum_{\lambda \,\le\, \sigma}
\mu(\lambda,\sigma)\,\phi_\lambda,
\]
where $\mu$ is the M\"{o}bius function of $NC(n)$. 
It is also common to write $k_n$ for the cumulant $k_{1_n}$, where
$1_n$ denotes the one-block partition $\bigl(\{1,\ldots,n\}\bigr)$.

The generalized moments $\phi_\lambda$ were defined in \eqref{defi_gen_moments} as products of usual moments over the blocks over partition $\lambda$. It turns out that the free cumulants inherit this important \emph{multiplicativity} property of generalized moments.  For example, 
if $\lambda = \bigl(\{1,4\},\{2,3\}\bigr)$, then
\[
\kappa_\lambda(X_1,X_2,X_3,X_4)
\;=\;
\kappa_2(X_1, X_4)\;\kappa_2(X_2, X_3).
\]
We will not prove this property but refer instead to Lectures 10 and 11 in \cite{nica_speicher06}.

In practice, one seldom needs to work directly with the M\"{o}bius function $\mu$ to compute free
cumulants. Far more useful is the key property that \emph{mixed free
cumulants} vanish whenever their arguments come from different free
subalgebras. This result underlies many concrete calculations, including
the free central limit theorem.

\section{Vanishing of mixed free cumulants}

The main property of free cumulants is that mixed cumulants of free variables 
vanish. In other words, if we have a set of variables 
\(X_{1},X_{2},\ldots ,X_{n}\) which can be split so that some of them lie in a 
subalgebra \(\mathcal{A}_{1}\) and the rest lie in a subalgebra \(\mathcal{A}_{2}\) 
that is free from \(\mathcal{A}_{1}\), then the mixed free cumulant 
\[
\kappa_{n}\bigl(X_{1},X_{2},\ldots ,X_{n}\bigr)
\]
is zero whenever at least one variable comes from \(\mathcal{A}_{1}\) and at 
least one variable comes from \(\mathcal{A}_{2}\).

\index{free cumulants!vanishing of mixed cumulants}

\begin{theo}\label{theorem_fundamental_cumulants}
Consider variables \(x_{1},\ldots ,x_{n}\) with \(n>1.\) Suppose that some of 
these variables belong to a subalgebra \(\mathcal{A}_{1}\) and the remaining 
variables belong to a subalgebra \(\mathcal{A}_{2}\), where 
\(\mathcal{A}_{1}\) and \(\mathcal{A}_{2}\) are free.  
If at least one \(x_{i}\) lies in \(\mathcal{A}_{1}\) and at least one \(x_{j}\) 
lies in \(\mathcal{A}_{2}\), then
\[
\kappa_{n}\bigl(x_{1},\ldots ,x_{n}\bigr) \;=\; 0.
\]
\end{theo}

\begin{proof}[Proof sketch]
We break the argument into three main steps.

\medskip

\noindent
\textbf{(1) One of the variables is constant.}

\begin{lemma}\label{lemma_freeness_and_constant}
Suppose \(x_{1},\ldots ,x_{n}\) (with \(n>1\)) includes at least one constant 
variable, say \(x_{i} = c \cdot 1\), where \(1\) is the unit.  Then
\[
\kappa_{n}\bigl(x_{1},\ldots ,x_{n}\bigr) \;=\; 0.
\]
\end{lemma}

\begin{proof} Suppose without loss of generality that $x_{i}=1.$
It is straightforward to 
check that \(\kappa_{2}(x,1) = \kappa_2(1, x) = 0\) for any \(x\), since 
\(\kappa_{2}(x,1) = \phi(x1) - \phi(x)\,\phi(1)\) and \(\phi(1) = 1\).
In order to prove the statement for $n>2,$ let us write 
\begin{equation*}
\phi\left( x_{1}\ldots x_{i-1}1x_{i+1}\ldots x_{n}\right) =\sum_{\lambda \in NC\left(
n\right) }\kappa_{\lambda }\left( x_{1},\ldots ,x_{i-1},1,x_{i+1},\ldots
,x_{n}\right) .
\end{equation*}%
By induction on $n$, all cumulants must equal zero except, perhaps, in two cases: when $%
\lambda =1_{n}$ and when $i$ is a one-element block of partition $\lambda .$
Hence, we obtain: 
\begin{eqnarray*}
\phi\left( x_{1}\ldots 1x_{i+1}\ldots x_{n}\right) &=&\kappa_{n}\left( x_{1},\ldots
,x_{i-1},1,x_{i+1},\ldots ,x_{n}\right) \\
&&+\sum_{\lambda \in NC\left( n-1\right) }\kappa_{\lambda }\left( x_{1},\ldots
,x_{i-1},x_{i+1},\ldots ,x_{n}\right) \\
&=&\kappa_{n}\left( x_{1},\ldots ,x_{i-1},1,x_{i+1},\ldots ,x_{n}\right) +\phi\left(
x_{1}\ldots x_{i-1}x_{i+1}\ldots x_{n}\right) .
\end{eqnarray*}%
This implies that 
\begin{equation*}
\kappa_{n}\left( x_{1},\ldots ,x_{i-1},1,x_{i+1},\ldots ,x_{n}\right) =0.
\end{equation*}
\end{proof}

\medskip

\noindent
\textbf{(2) The variables form an alternating sequence from two free subalgebras.}

Assume \(x_{1},x_{3},x_{5},\ldots\in \mathcal{A}_{1}\) and 
\(x_{2},x_{4},x_{6},\ldots\in \mathcal{A}_{2}\).  Moreover, at first suppose each \(x_{i}\) 
is centered, i.e.\ \(\phi(x_{i})=0\).  We will relax this assumption later. We claim that 
\(\kappa_{n}(x_{1},\ldots,x_{n})=0\).  It is clear for \(n=2\) because then
\(\kappa_{2}(x_{1},x_{2}) = \phi(x_{1}x_{2}) - \phi(x_{1})\phi(x_{2})=0.\)

For \(n>2\),
\[
\phi\bigl(x_{1}\cdots x_{n}\bigr)
 \;=\;
 \kappa_{n}(x_{1},\ldots,x_{n})
 + \sum_{\substack{\lambda\in NC(n)\\ \lambda \neq 1_{n}}}
     \kappa_{\lambda}(x_{1},\ldots,x_{n}).
\]
If a nontrivial block of \(\lambda\neq 1_{n}\) connects variables from 
different algebras, the induction hypothesis forces \(\kappa_{\lambda}=0.\)  
If each block of \(\lambda\) stays within a single algebra, then the alternating structure of the sequence $x_{1},\ldots , x_{n}$ forces this partition to have single-element blocks. (Note that it is crucial here that the partitions are non-crossing!).  But each 
\(x_i\) is centered, so each singleton block contributes zero factor. By multiplicativity, $\kappa_\lambda = 0$.   
Thus \(\phi(x_{1}\cdots x_{n}) = \kappa_{n}(x_{1},\ldots,x_{n}).\)

By freeness, an alternating product of centered variables in 
\(\mathcal{A}_{1}\cup \mathcal{A}_{2}\) has \(\phi(x_{1}\cdots x_{n})=0,\) 
so \(\kappa_{n}(x_{1},\ldots,x_{n})=0.\)

Removing the ``centered'' assumption can be done by using multi-linearity of cumulants, replacing 
\(x_{i}\) by \(x_{i} - \phi(x_{i})\) and using 
Lemma~\ref{lemma_freeness_and_constant} for the constant parts.

Indeed, in order to show that $\kappa_{n}\left( x_{1},\ldots ,x_{n}\right) =0$ holds for alternating
sequences $x_{1},\ldots ,x_{n}$ even if variables $x_{i}$ are not centered, 
we can write equations like 
\begin{eqnarray*}
\kappa_{n}\left( x_{1},\ldots ,x_{n}\right) &=&\kappa_{n}\left( x_{1}-\phi\left(
x_{1}\right) ,x_{2},\ldots ,x_{n}\right) +\kappa_{n}\left( \phi\left( x_{1}\right)
,x_{2},\ldots ,x_{n}\right) \\
&=&\kappa_{n}\left( x_{1}-\phi\left( x_{1}\right) ,x_{2},\ldots ,x_{n}\right)
\end{eqnarray*}%
and apply them several times.

\medskip

\noindent
\textbf{(3) The general case.}

Let \(x_{1},\ldots,x_{n}\) be an arbitrary sequence in 
\(\mathcal{A}_{1}\cup\mathcal{A}_{2}\), possibly with consecutive runs in 
the same algebra.  Partition the indices into blocks:
\[
X_{1} = x_{1}\cdots x_{i_{1}},
\quad
X_{2} = x_{i_{1}+1}\cdots x_{i_{2}},
\quad
\dots,
\quad
X_{m} = x_{i_{m-1}+1}\cdots x_{n},
\]
so that \(X_{1},\dots,X_{m}\) \emph{alternate} between 
\(\mathcal{A}_{1}\) and \(\mathcal{A}_{2}\).  From the step~(2) argument, 
\(\kappa_{m}(X_{1},\dots,X_{m})=0\).

On the other hand, by the general product formula 
(Theorem~\ref{theorem_cumulants_of_products} below),
\[
\kappa_{m}(X_{1},\dots,X_{m})
\;=\;
\sum_{\substack{\pi\in NC(n)\\ \pi \,\vee\, \widehat{0}_{m}=1_{n}}}
 \kappa_{\pi}(x_{1},\dots,x_{n}),
\]
where \(\widehat{0}_{m}\) is the partition 
\(\{\{1,\dots,i_{1}\},\dots,\{\,i_{m-1}+1,\dots,n\}\}.\)

Since \(\kappa_{m}(X_{1},\dots,X_{m})=0\), we get
\[
0
\;=\;
\kappa_{n}(x_{1},\ldots,x_{n})
 + \sum_{\substack{\pi\neq 1_{n}\\ \pi \,\vee\, \widehat{0}_{m}=1_{n}}}
     \kappa_{\pi}(x_{1},\ldots,x_{n}).
\]
The condition $\pi \vee \widehat{0}_{m}=1_{n}$ implies that $\pi $ must connect two variables
from different algebras. But any such \(\pi\neq 1_{n}\) 
factors into smaller blocks and at least one of them contains variables from different algebras, so by induction \(\kappa_{\pi}=0\).  
Hence \(\kappa_{n}(x_{1},\ldots,x_{n})=0\).

\end{proof}

\medskip

\index{free cumulants!formula for products}

\begin{theo}[Cumulants of products]
\label{theorem_cumulants_of_products}
Let \(x_{1},\dots,x_{n}\) be given, and fix 
\(0=i_{0}<i_{1}<\dots<i_{m}=n\).  Define
\[
X_{1}\;=\;x_{1}\cdots x_{i_{1}},\quad
X_{2}\;=\;x_{i_{1}+1}\cdots x_{i_{2}},\quad
\dots,\quad
X_{m}\;=\;x_{i_{m-1}+1}\cdots x_{n}.
\]
Then
\[
\kappa_{m}(X_{1},\ldots,X_{m})
\;=\;
\sum_{\substack{\pi\in NC(n)\\ \pi \,\vee\, \widehat{0}_{m}=1_{n}}}
  \kappa_{\pi}(x_{1},\ldots,x_{n}),
\]
where \(\widehat{0}_{m}\) is the partition 
\(\bigl\{\{1,\dots,i_{1}\},\{\,i_{1}+1,\dots,i_{2}\},\dots,\{\,i_{m-1}+1,\dots,n\}\bigr\}.\)
\end{theo}

\index{partition lift|textbf}
\begin{proof}[Proof of Theorem~\ref{theorem_cumulants_of_products}]
Recall that we have a sequence \(0=i_{0}<i_{1}<\dots<i_{m}=n\) such that
\[
X_{1}\;=\;x_{1}\cdots x_{i_{1}},\quad
X_{2}\;=\;x_{i_{1}+1}\cdots x_{i_{2}},\quad
\dots,\quad
X_{m}\;=\;x_{i_{m-1}+1}\cdots x_{n}.
\]

\begin{figure}[tbph]
\begin{center}
\includegraphics[width=10cm]{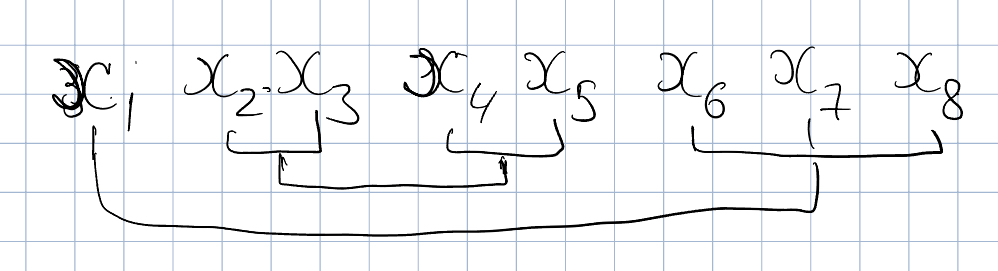}
\end{center}
\caption{Lift of partition $\protect\pi =  ((1, 4), (2,3))$ based on the sequence $i_1 = 1 < i_2 = 3 < i_3 = 5 < i_4 = 8$}
\label{figure_partition_lift}
\end{figure}

First, we need a lemma on the M\"obius function. To formulate it, let us define \emph{a lift of a partition}
of $\left\{ 1,\ldots ,m\right\} $ to a partition of $\left\{ 1,\ldots
,n\right\},$ based on a sequence   $(1=i_{0}<i_{1}<\dots<i_{m}=n)$. We simply substitute every element $k$ in $\left\{ 1,\ldots ,m\right\} $ with connected elements $i_{k-1}+1, i_{k-1} + 2, \ldots ,i_{k}$ in  $\left\{ 1,\ldots
,n\right\}$ (by convention $i_0 = 0$). If elements $l$ and $k$ were connected in $\pi$, then the elements $i_{l-1}+1,\ldots ,i_{l}$ and $i_{k-1}+1,\ldots ,i_{k}$ belong to the same block in the lift of $\pi$.  It is easy to see that a lift of a
non-crossing partition is non-crossing. We will denote a lift of partition $\pi $ by $\widehat{\pi }.$

For example if $n = 4$, $m = 8$, and the sequence is $i_1 = 1 < i_2 = 3 < i_3 = 5 < i_4 = 8$, then the partition $\pi = ((1, 4), (2,3))$ is lifted to $((1, 6, 7, 8), (2, 3, 4, 5)$.  See Figure \ref{figure_partition_lift}

\begin{lemma}\label{lemma_Mobius_lift}
\[
\mu\bigl(\pi,1_{m}\bigr)
\;=\;
\mu\bigl(\,\widehat{\pi},1_{n}\bigr)
\quad\text{for all } \pi\in NC(m).
\]
\end{lemma}

\begin{proof}
The map \(\pi\mapsto\widehat{\pi}\) is an isomorphism of poset intervals:
\(\,[\pi,1_{m}]\subset NC(m)\,\) onto
\(\,[\widehat{\pi},1_{n}]\subset NC(n)\,\).  
(Indeed, it is clear that if  if $\pi \leq \lambda $ then $\widehat{\pi }\leq 
\widehat{\lambda }$. In the other direction suppose that $\widehat{\pi } \leq \tau$ for some $\tau$. That is, $\tau$ is a coarsening of $\widehat{\pi }$. Then, it is easy to see that $\tau $ is a lift of a corresponding coarsening of  $\pi$ in $NC\left( m\right) ,$ which we can
call $\lambda ,$ so that $\tau =\widehat{\lambda }$, and $\pi \leq
\lambda .$)

Hence \(\mu\) values coincide on those intervals by definition of the M\"obius function.
\end{proof}

Using the definition of free cumulants,
\[
\kappa_{m}(X_{1},\dots,X_{m})
\;=\;
\sum_{\lambda \in NC(m)} \mu(\lambda,1_{m}) \,\phi_{\lambda}(X_{1},\dots,X_{m}).
\]
But \(\phi_{\lambda}(X_{1},\dots,X_{m})\) can be rewritten as 
\(\phi_{\widehat{\lambda}}(x_{1},\dots,x_{n})\).  Applying 
Lemma~\ref{lemma_Mobius_lift} to replace \(\mu(\lambda,1_{m})\) with 
\(\mu(\widehat{\lambda},1_{n})\), we find 

\begin{eqnarray*}
k_{m}\left( X_{1},\ldots ,X_{m}\right) &=&\sum_{0_{m} \leq \lambda }\mu
( \lambda ,1_{n}) \phi_{\lambda }\left( x_{1},\ldots ,x_{n}\right) \\
&=&\sum_{\pi \vee \widehat{0}_{m}=1_{n}}k_{\pi }\left( x_{1},\ldots
,x_{n}\right) ,
\end{eqnarray*}%
where the second line follows from a more general variant of Rota's inversion formula -- formula \eqref{equ_Rota_inversion3} on page \pageref{equ_Rota_inversion3}.

\end{proof}

\section{Additivity and $R$-transform}


The vanishing of mixed free cumulants leads to an additivity property. Let us
denote $\kappa_{n}(X,\ldots,X)$ by $\kappa_{n}(X)$.

\index{free cumulants!additivity}
\begin{theo}
\label{theo_additivity_cumulants}
Let $X$ and $Y$ be free. Then 
\begin{equation*}
\kappa_{n}(X+Y) = \kappa_{n}(X) + \kappa_{n}(Y).
\end{equation*}
\end{theo}

\begin{proof}
By definition,
\[
\kappa_{n}(X+Y) \;=\;
   \kappa_{n}\bigl(
     \underbrace{X+Y,\ldots,X+Y}_{n\text{-times}}
   \bigr).
\]
We expand this expression via multilinearity and note that all mixed cumulants vanish by Theorem~\ref{theorem_fundamental_cumulants}. Hence,
\[
\kappa_{n}(X+Y)
 \;=\;
 \kappa_{n}(X,\ldots,X)
 \;+\;
 \kappa_{n}(Y,\ldots,Y)
 \;=\;
 \kappa_{n}(X) + \kappa_{n}(Y).
\]
\end{proof}

In order to apply this result, we need to be able to compute the cumulants of individual random variables. It is natural to form a generating function
\[
R_X(z) \;=\; \sum_{n=1}^{\infty} \kappa_{n}(X)\,z^{n-1} = \phi(X) + \Var(X) z + \ldots
\]

\index{R-transform|textbf}
By using the explicit formula for the M\"obius function in Section \ref{section_mobius_function_NC}, it is not difficult to show that for bounded $X$ these series converge for sufficiently small $z$ and define an analytic function which is called the \emph{$R$-transform} of the r.v. $X$.

 For bounded $X$, the Cauchy transform $G_X$ defined in \eqref{defi_Cauchy_transform}can be expanded in series: 
\[
G_X(u) \;=\; \phi\Bigl[\frac{1}{u - X}\Bigr]
 \;=\; \frac{1}{u} \;+\; \frac{\phi(X)}{u^2} \;+\; \frac{\phi(X^2)}{u^3} \;+\;\dots,
\]
convergent for sufficiently large $u \in \bC$.

\index{R-transform!relation to Cauchy transform}
\begin{theo}[Relation between $G_X$ and $R_X$]
\label{theo_R_transform_identity}
Let $X$ be a bounded random variable. For all sufficiently large $u \in \bC$
\[
R_X\bigl(G_X(u)\bigr) \;+\; \frac{1}{G_X(u)} \;=\; u.
\]
\end{theo}
One can check that for the map $u \to G_X(u)$ is invertible on the domain $|u| > C$ if $C$ is sufficiently large.   If $G_X^{(-1)}(z)$ denotes the functional inverse of $G_X(u)$, then the theorem shows that 
\begin{equation}
\label{equ_defi_Rtransform}
R_X(z) \;+\; \frac{1}{z} \;=\; G_X^{(-1)}(z)
\end{equation}
for small $z \in \bC$. 
Equivalently,
\[
G_X\!\Bigl(\frac{1}{z} \;+\; R_X(z)\Bigr) \;=\; z.
\]

For unbounded r.v.s. $X$, formula \eqref{equ_defi_Rtransform} serves as the definition of the $R$-transform $R_X(z)$, however, one has to be careful in the choice of the domain on which the Cauchy transform $G_X(z)$ is invertible. Typically this is a domain in the upper complex half-plane such that $\Im z > C$ and $\Im z > \alpha \, \Re z$ for some positive $C$ and $\alpha$.

\begin{proof}
Define
\begin{equation}
\label{defi_M}
M_X(z)
 \;:=\;
 1 \;+\; \sum_{n=1}^{\infty} \phi\bigl(X^n\bigr)\,z^n.
\end{equation}
By the definition of free cumulants,
\[
\phi\bigl(X^n\bigr)
 \;=\;
 \sum_{\pi \in NC(n)}\kappa_{\pi}(X,\ldots,X).
\]

\begin{figure}[h]
\begin{center}
\includegraphics[width=10cm]{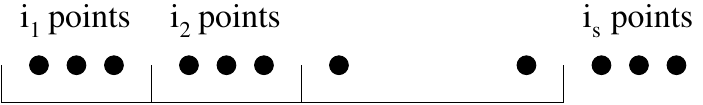}
\end{center}
\caption{Decomposition of a partition by its first block.}
\label{figure_R_cumulants}
\end{figure}

To organize these sums, fix the first block of $\pi$. By non-crossingness, the remaining parts of $\pi$ lie between the elements of that block.  Then we can use the multiplicativity of free
cumulants to write $\kappa_{\pi }$ as a product. Finally, when we add up all
these products and sum over all possibilities for the first block,  we  obtain
the following formula: 
\begin{eqnarray*}
\phi\left( X^{n}\right) &=&\sum_{s=1}^{n}
\sum_{\substack{ i_{1}+\ldots
i_{s}+s=n  \\ i_{r}\geq 0}}\kappa_{s}
\left( \sum_{\pi _{1}\in NC\left(
i_{1}\right) }\kappa_{\pi _{1}}\right) \ldots 
\left( \sum_{\pi _{s}\in NC\left(
i_{s}\right) }\kappa_{\pi _{s}}\right) \\
&=&\sum_{s=1}^{n}\sum_{\substack{ i_{1}+\ldots i_{s}+s=n  \\ i_{r}\geq 0}}%
\kappa_{s}\phi\left( X^{i_{1}}\right) \ldots \phi\left( X^{i_{s}}\right) .
\end{eqnarray*}
Here, we used $\kappa_\pi$ for $\kappa_\pi(X, \ldots, X)$ for shortness. 

Symbolically, this argument is illustrated in Figure \ref{figure_R_cumulants}.

Hence,
\[
M_X(z)
 \;=\;
 1
 \;+\;
 \sum_{n=1}^\infty
   \sum_{s=1}^n
     \kappa_{s}\,z^s
     \Bigl[\sum_{i=0}^\infty \phi(X^i)\,z^i\Bigr]^s
 \;=\;
 1 \;+\; \sum_{s=1}^\infty \kappa_s\,z^s \bigl(M_X(z)\bigr)^s.
\]
Extracting the factor $z\,M_X(z)$ leads to the functional equation:
\begin{equation}
\label{equ_C_and_M}
M_X(z)
 \;=\;
 1 \;+\; z\,M_X(z)\,R_X\!\bigl(z\,M_X(z)\bigr),
\quad
\text{i.e.}
\quad
\frac{1}{z}
 \;=\;
 \frac{1}{z\,M_X(z)}
 \;+\;
 R_X\!\bigl(z\,M_X(z)\bigr),
\end{equation}
valid for small $z \ne 0$. 

One also checks that $z\,M_X(z) = G_X\!\bigl(z^{-1}\bigr)$, so setting $u=z^{-1}$ gives
\[
u
 \;=\;
 \frac{1}{G_X(u)}
 \;+\;
 R_X\!\bigl(G_X(u)\bigr),
\]
valid for all sufficiently large $u \in \bC$.
\end{proof}


\index{R-transform!scaling}
Another important property of the $R$-transform is its behavior under scaling. Indeed,
\begin{equation}
\label{scaling_R}
R_{aX}(z)
 \;=\;
 \sum_{n=1}^\infty \kappa_n(aX)\,z^{n-1}
 \;=\;
 a\,\sum_{n=1}^\infty \kappa_n(X)\,(a\,z)^{\,n-1}
 \;=\;
 a\,R_X(a\,z).
\end{equation}


\index{R-transform!additivity}
\begin{theo}[Additivity of the $R$-transform]
\label{theo_additivity_R}
Let $X$ and $Y$ be free. Then 
\[
R_{X+Y}(z)
 \;=\;
 R_X(z) \;+\; R_Y(z).
\]
\end{theo}

\begin{proof}
Immediate from Theorem~\ref{theo_additivity_cumulants} and the definition of $R_X(z)$.
\end{proof}


\index{semicircle r.v.!free cumulants}
\index{semicircle r.v.!R-transform}
\begin{exa}[Semicircle distribution]
\label{exa_semicircle_cumulants}
Let $s$ have the semicircle distribution as in Example~\ref{exa_semicircle}. 
Its Cauchy transform is
\[
G_s(z)
 \;=\;
 \frac{1}{2}\Bigl(z \;-\;\sqrt{z^2 - 4}\Bigr),
\]
which satisfies $G_s(z) + 1/G_s(z) = z$. By the preceding theorem, we see that
\[
R_s(u)
 \;=\;
 u.
\]
In other words, the only non-vanishing free cumulant of $s$ is $\kappa_2(s)=1$. This perfectly parallels the classical Gaussian, which is characterized by having its only non-zero \emph{classical} cumulant equal to 1 in second order.

Moreover, by the scaling property~\eqref{scaling_R}, the sum of $n$ free semicircle variables with the same distribution is again a semicircle, scaled by $\sqrt{n}$. This is the exact analogue of how sums of independent Gaussians scale by $\sqrt{n}$. 
\end{exa}

\index{arcsine distribution!free cumulants}
\index{arcsine distribution!R-transform}
\begin{exe}[Arcsine distribution]
\label{exe_arcsine}
Let $X$ be a self-adjoint random variable with the arcsine distribution
(see p.~\pageref{arcsine_distribution}). Show that
\[
R_X(u)
 \;=\;
 \frac{1}{u}\bigl(\sqrt{1 + 4u^2}\;-\;1\bigr)
 \;=\;
 2\,u\,\sum_{k=0}^{\infty} C_k\,u^{2k},
\]
where $C_k = \tfrac{1}{k+1}\,\binom{2k}{k}$ are the Catalan numbers.
\end{exe}

\index{bernoulli distribution!R-transform}
\begin{exe}[Free Bernoulli]
\label{exe_bernoulli}
Let $X$ be a self-adjoint random variable with the distribution 
\[
\mu_{\pm 1}
 \;=\;
 \tfrac12\,\delta_{-1} \;+\; \tfrac12\,\delta_{1}.
\]
Show that
\[
R_X(u)
 \;=\;
 \frac{1}{2\,u}\,\bigl(\sqrt{1 + 4\,u^2}\;-\;1\bigr).
\]
In particular, this shows that if $X$ and $Y$ are free and both have 
$\mu_{\pm 1}$, then $X+Y$ has the free arcsine law (Exercise~\ref{exe_arcsine}).
\end{exe}


For an unbounded self-adjoint operator $X$ (affiliated with $(\mathcal{A},\phi)$),
the concepts of freeness still make sense via truncations (defined by functional calculus). Two unbounded variables are free if all their finite-rank truncations are free. One can still define the probability distribution of $X$ via spectral theory, and its Cauchy transform 
\[
G_X(z)
 \;=\;
 \phi\bigl((z - X)^{-1}\bigr)
\]
is well-defined in appropriate regions of the upper half-plane. For sufficiently large imaginary part, $G_X(z)$ is invertible, so one can define $R_X$ using equation \eqref{equ_defi_Rtransform}. One may show that $R_{X+Y} = R_X + R_Y$ still holds. This fact is crucial in studying free infinitely-divisible distributions.


\index{Cauchy distribution!R-transform}
\begin{exe}[Cauchy distribution]
Let $X$ have the standard Cauchy distribution, i.e.\ its spectral distribution has the density
\[
p(x)
 \;=\;
 \frac{1}{\pi}\,\frac{1}{1 + x^2}.
\]
Then $\phi\bigl(f(X)\bigr)$ can be defined suitably for bounded $f$, even though $X$ has infinite variance and not even a well-defined mean. Show that
\[
G_X(z)
 \;=\;
 \frac{1}{\,z + i\,}
\quad\text{for } \mathrm{Im}(z) > 0,
\]
and that the associated $R$-transform is
\[
R_X(z)
 \;=\;
 -\,i.
\]
Thus, even though no combinatorial moment-cumulant calculations are possible (the distribution lacks moments), the $R$-transform is still well-defined.
\end{exe}

\section{Expectation of products of free elements}

The following results are very useful in calculations. 

\begin{theo}\label{theorem_expectation_products}
Let $\{ a_{1},\ldots ,a_{n}\}$ and $\{ b_{1},\ldots ,b_{n}\}$ be free. Then
\begin{align}
\phi\bigl( a_{1}b_{1}a_{2}b_{2}\cdots a_{n}b_{n}\bigr)
&=\sum_{\pi \in NC(n)} \kappa_{\pi}( a_{1},\ldots ,a_{n}) \,\phi_{K(\pi)}( b_{1},\ldots ,b_{n})
\label{eq:phi-alternating-product}\\
&=\sum_{\pi \in NC(n)} \sum_{\lambda \leq \pi }
   \,\mu ( \lambda ,\pi)\, \phi_{\lambda} ( a_{1},\ldots ,a_{n})
  \phi_{K(\pi)}( b_{1},\ldots ,b_{n}),
\label{eq:phi-alternating-product2}
\end{align}
where $K(\pi)$ denotes the Kreweras complement of $\pi$, and  $\mu$
is the Möbius function on $NC(n)$.
\end{theo}

\begin{proof}[Proof of Theorem \ref{theorem_expectation_products}]
By the definition of free cumulants, their multiplicativity,  and by vanishing of all mixed cumulants (due to freeness of $\{a_i\}$ and $\{b_i\}$), we have
\[
   \phi\bigl(a_{1}b_{1}\cdots a_{n}b_{n}\bigr)
   \;=\;
   \sum_{\substack{\pi_{a} \in NC(n)\\ \pi_{b}\in NC(n) \\ \pi_a \cup \pi_b \in NC(2n)}} 
        k_{\pi_{a}}(a_{1},\ldots ,a_{n})\;k_{\pi_{b}}(b_{1},\ldots ,b_{n}),
\]
where $\pi_{a}$ is a noncrossing partition of the indices of the $a_{i}$'s, $\pi_{b}$ is a noncrossing partition of the indices of the $b_{i}$'s, and the combined partition $\pi_{a}\cup \pi_{b}$ is noncrossing on the full alternating string $a_{1},b_{1},\ldots,a_{n},b_{n}$.  

Recalling the definition of the Kreweras complement, one sees that each valid pair $(\pi_{a},\pi_{b})$ is characterized by $\pi_{b} \leq K(\pi_{a})$; hence,
\begin{align}
   \phi\bigl(a_{1}b_{1}\cdots a_{n}b_{n}\bigr)
   &\;=\;
   \sum_{\substack{\pi_{a}\in NC(n)\\ \pi_{b} \leq K(\pi_{a})}}
       k_{\pi_{a}}(a_{1},\ldots,a_{n})\,k_{\pi_{b}}(b_{1},\ldots,b_{n})
   \\
   &\;=\;
   \sum_{\pi\in NC(n)}
       k_{\pi}(a_{1},\ldots,a_{n})\,\phi_{K(\pi)}(b_{1},\ldots,b_{n}).
\end{align}
This proves the first equality in~\eqref{eq:phi-alternating-product}.

For the second equality, one applies Möbius inversion to write
\[
 \kappa_{\pi}( a_{1},\ldots ,a_{n}) =  \sum_{\lambda \leq \pi }
   \,\mu ( \lambda ,\pi)\, \phi_{\lambda} ( a_{1},\ldots ,a_{n})
  \] 
\end{proof}

\medskip

\begin{theo}\label{theorem_expectation_products2}
Let $\{ a_{1},\ldots ,a_{n}\}$ and $\{ b_{1},\ldots ,b_{n}\}$ be free. Then 
\[
  \kappa_n\bigl(a_{1}b_{1},\, a_{2}b_{2},\,\ldots,\,a_{n}b_{n}\bigr)
  \;=\;
  \sum_{\pi \in NC(n)}
       k_{\pi}(a_{1},\ldots,a_{n})
       \,\kappa_{K(\pi)}(b_{1},\ldots,b_{n}),
\]
where $K(\pi)$ is again the Kreweras complement of $\pi$.
\end{theo}

\begin{proof} (Sketch)

By Theorem \ref{theorem_cumulants_of_products}We have 
\[
\kappa_n\bigl(a_{1}b_{1},\, a_{2}b_{2},\,\ldots,\,a_{n}b_{n}\bigr) =
\sum_{\substack{\pi \in NC(2n) \\ \pi \vee \hat 0_n = 1_{2n}}} \kappa  \bigl(a_{1}, b_{1}, a_{2}, b_{2},\ldots,\,a_{n}, b_{n}\bigr),
\]
where $\hat 0_n$ is the partition $\{(1, 2), (3, 4), \ldots, (2 n - 1, 2n)\}$. By vanishing of free cumulants $\pi$ should not couple $a$'s and $b$'s, therefore, we can write $\pi = \pi_a \cup \pi_b$, where $\pi_a \in NC(n)$ and $\pi_b \in NC(n)$. Moreover, since $\pi$ is non-crossing we must have $\pi_b \leq K(\pi_a)$ and if strict inequality holds then $(\pi_a \cup \pi_b) \vee \hat 0_n \ne 1_{2n}$. For example if $\pi_a = \{(1), (3), (5)\}$, then $\pi_b$ must be $\{(2, 4, 6)\} = K(\pi_a)$ so that $\pi$ would connect all blocks in $\{(1, 2), (3, 4), (5, 6)\}$.  After some effort this can be proved in general. 

Then, by multiplicativity of free cumulants, one can write:
\[
\kappa_n\bigl(a_{1}b_{1},\, a_{2}b_{2},\,\ldots,\,a_{n}b_{n}\bigr) = 
\sum_{\substack{\pi_a \in NC(n) \\ \pi_b = K(\pi_a) }} \kappa_{\pi_a}  \bigl(a_{1}, a_{2},\ldots,\,a_{n}\bigr) \kappa_{\pi_b}  \bigl(b_{1}, b_{2},\ldots,\,b_{n}\bigr).
\]
\end{proof}

\section{The Möbius function for non-crossing partitions}
\label{section_mobius_function_NC}

\subsection*{Multiplicativity and Kreweras complement}

\index{Mobius function!for non-crossing partitions}

The M\"obius function of the lattice of non-crossing partitions is
multiplicative in the following sense. Suppose \(\lambda \le \sigma\) (meaning that
\(\lambda\) is a refinement of \(\sigma\)) and that
\(\sigma = \sigma_1 \sigma_2\). Here we think about blocks of $\sigma$, $\sigma_1$ and $\sigma_2$ as the permutation cycles and  assume \(\sigma_1\) and \(\sigma_2\)
act on disjoint subsets of \(\{1,\dots,n\}\).  Then we can factorize
\(\lambda\) as \(\lambda_1 \lambda_2\) with \(\lambda_i \le \sigma_i\).
In this situation,
\[
   \mu(\lambda,\sigma) \;=\; \mu(\lambda_1,\sigma_1)\,\mu(\lambda_2,\sigma_2).
\]

For more information on multiplicativity of the Mobius function, see Lecture 10 in \cite{nica_speicher06}.

 This multiplicativity property
immediately reduces the problem of computing
\(\mu(\lambda,\sigma)\) to the case where \(\sigma = 1_n\) (the
\emph{maximal} partition with only one block).  Indeed, if we can
evaluate \(\mu(\lambda, 1_n)\) for all \(\lambda\) and \(n\), then
we can handle the general case by factoring.

Next, recall that the Kreweras complement \(K\) gives an anti-isomorphism
of \(NC(n)\).  In particular,
\[
   \mu(\lambda,\sigma)
   \;=\;
   \mu\bigl(K(\sigma),\,K(\lambda)\bigr).
\]
Setting \(\sigma = 1_n\) in the above gives
\[
   \mu(\lambda, 1_n)
   \;=\;
   \mu\bigl(0_n,\,K(\lambda)\bigr),
\]
where \(0_n\) is the \emph{minimal} partition \(\{\{1\},\dots,\{n\}\}\). 
Hence, knowing \(\mu(0_k,1_k)\) for all \(k\) also controls
\(\mu(\lambda,1_n)\) by passing through the Kreweras complement.  

\subsection*{Speicher’s formula for \(\mu(0_n,1_n)\) via Catalan recursion}

\begin{theo}[Speicher]
\label{theo_Speicher}
The Möbius function of the lattice of non-crossing partitions
\(NC(n)\) satisfies
\[
   \mu(0_n,1_n)
   \;=\;
   (-1)^{\,n-1}\,C_{\,n-1},
\]
where \(C_{\,k}\) denotes the \(k\)-th Catalan number.
\end{theo}

\begin{figure}[tbp]
\centering
\includegraphics[width=8cm]{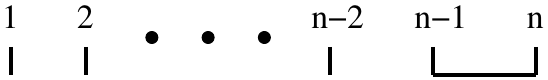}
\caption{A partition \(a\) with last two points blocked together.}
\label{figure_partition_w}
\end{figure}

\begin{figure}[tbp]
\centering
\includegraphics[width=8cm]{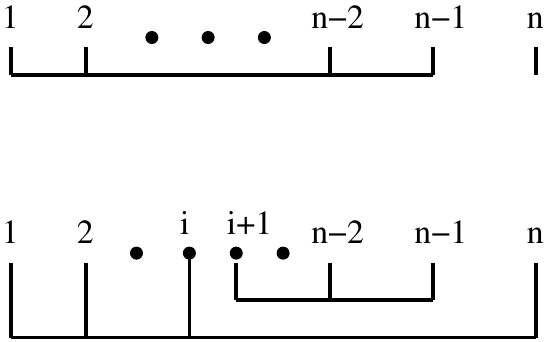}
\caption{Two possible types of partition  \(x\).}
\label{figure_partition_type1and2}
\end{figure}

\begin{proof}
We use \eqref{mobius01} with the particular partition
\[
   a \;=\; \bigl\{\{1\},\,\{2\},\,\dots,\,\{n-2\},\,\{n-1,n\}\bigr\},
\]
as in Figure~\ref{figure_partition_w}.  

It is clear that any partition
\(x \neq 1_n\) satisfying \(x \vee a = 1_n\) must have one of the two forms shown in
Figure \ref{figure_partition_type1and2}. In the first case, $x$ factors into two blocks of length $n - 1$ and $1$, and in the second case, it factors into two blocks of length $n - k$ and $k$ for $ 2 \leq k \leq n - 2$

 By using \eqref{mobius01}  and multiplicativity of $\mu$, we
 thus get a recursion
\[
   \mu(0_n,1_n)
   \;=\;
   - \sum_{k=1}^{n-1}
   \mu(0_k,1_k)\,\mu(0_{n-k},1_{n-k}),
\]
and this implies that $\left( -1\right) ^{n-1}\mu \left( 0_{n},1_{n}\right) $
satisfy the usual recursion for the Catalan numbers
\(C_{n-1}\).  One checks initial values to see that
\(\mu(0_n,1_n)\) coincides with \((-1)^{n-1}C_{n-1}\).  
\end{proof}

\subsection*{Examples of block-by-block computation}

Let us illustrate this multiplicative approach to computing Möbius
values.  

\begin{itemize}
\item 
Consider \(\lambda = \{\{1,4\},\{2,3\}\}\).  By
the Kreweras complement,
\[
   \mu(\lambda,1_4)
   \;=\;
   \mu\bigl(0_4,\,K(\lambda)\bigr),
\]
and one checks that
\(K(\lambda) = \{\{1,3\},\,\{2\},\,\{4\}\}\).  Since those three
blocks \(\{1,3\},\{2\},\{4\}\) are on disjoint subsets, 
\[
   \mu\bigl(0_4,\,K(\lambda)\bigr)
   \;=\;
   \mu(0_2,1_2)\,\mu(0_1,1_1)\,\mu(0_1,1_1)
   \;=\;
   (-1)\,\cdot 1 \,\cdot 1
   \;=\;
   -1.
\]
Hence \(\mu(\lambda,1_4)=-1\).

\item
As another example, let \(\lambda = (\{1\},\,\{2\},\,\{3,4\})\).
Then the Kreweras complement is
\[
   K(\lambda) \;=\; (\{1,2,4\},\,\{3\}).
\]
Hence
\[
   \mu(\lambda,1_4)
   \;=\;
   \mu\bigl(0_4, K(\lambda)\bigr)
   \;=\;
   \mu(0_3,1_3)\,\mu(0_1,1_1)
   \;=\;
   \bigl((-1)^{2}C_2\bigr)\cdot 1
   \;=\;
   2,
\]
since \(C_2 = 2\). 
\end{itemize}

\section{NC partitions and permutations}
\label{section_partitions_permutations}

\index{geodesic permuations}
This section deviates from the topic of free cumulants and discusses a surprising relationship between the lattice
of non-crossing partitions and the group of permutations. Let $S_{n}$ be the group of permutations of the $n$-element set $\{1,\dots,n\}$. The \emph{length} of a permutation%
\index{length of permutation} $\sigma$ is the minimal number of transpositions whose product is $\sigma.$ We
denote this quantity by $|\sigma|.$ If the permutation $\sigma$ decomposes into $c(\sigma)$ disjoint cycles, then
\[
|\sigma| = n - c(\sigma).
\]
We define the distance between two permutations $\sigma$ and $\tau$ by
\[
d(\sigma, \tau) = \bigl|\sigma^{-1}\tau\bigr|.
\]
It is straightforward to check that this indeed defines a metric on $S_n$. 

\begin{exe}
Check the triangle inequality for the metric $d(\sigma, \tau)$.
\end{exe}

We say that a permutation $\sigma$ \emph{belongs to a geodesic}%
\index{geodesic}
between $\rho$ and $\tau$ if 
\[
d(\rho,\tau) \;=\; d(\rho,\sigma)\;+\;d(\sigma,\tau).
\]
Equivalently, $\sigma$ is on a shortest path from $\rho$ to $\tau$ with respect to the distance $d$.

Let $\tau$ be the $n$-cycle $(\,1\,2\,\dots\,n\,)$. We consider all permutations $\sigma$ that lie on some geodesic between the identity permutation $e$ and $\tau$. Concretely, $\sigma$ satisfies
\[
d(e, \tau)\;=\;d(e,\sigma)\;+\;d(\sigma,\tau)\quad\Longleftrightarrow\quad
|\tau|\;=\;|\sigma|+|\sigma^{-1}\tau|.
\]
Since $|\tau| = n-1$, this becomes
\[
|\sigma| + |\tau^{-1}\sigma| = n-1,
\]
or equivalently,
\[
c(\sigma) \,+\, c\bigl(\tau^{-1}\sigma\bigr) = n \;+\; 1.
\]
We give this set of permutations a partial order: we say $\sigma \prec \pi$ if 
\[
d(e,\pi) \;=\; d(e,\sigma)\;+\;d(\sigma,\pi),
\]
i.e.
\[
|\pi| \;=\;|\sigma|\;+\;|\sigma^{-1}\pi|.
\]
One can check that this is indeed a partial order on the set of all permutations lying on a geodesic from $e$ to $\tau$, and that this partially ordered set is in fact a lattice. Even more remarkably, this lattice is isomorphic to the lattice of non-crossing partitions of $\{1,2,\ldots,n\}$.

\subsection*{Bijection between geodesic permutations and non-crossing partitions.}

Recall that each permutation can be decomposed into cycles with disjoint supports. From these cycles, we get a set-partition of $\{1,\dots,n\}$ by taking as blocks the supports of each cycle. We will think about elements of the set $\{1,\dots,n\}$ as points on the circle put in counterclockwise order. Then we can restrict our attention to permutations $\sigma$ for which cycles are in the counter-clockwise order. They correspond to partitions of the ordered set $\{1, \ldots, n\}$.

\begin{lemma}
\label{lemma:geodesic-nc}
Let $\tau = (1\,2\,\dots\,n)$. A permutation $\sigma$ lies on a geodesic between $e$ and $\tau$ if and only if its associated partition (via disjoint cycle decomposition) is a non-crossing partition of the ordered set $\{1, \ldots, n\}$. Moreover, every non-crossing partition of $\{1,2,\dots,n\}$ arises from some $\sigma$ lying on a geodesic from $e$ to $\tau$.
\end{lemma}

\begin{proof}
First suppose $\sigma$ is on a geodesic from $e$ to $\tau$. Since $d(e,\tau) = n-1$, there is a chain
\[
e = \sigma_0,\;\sigma_1,\;\dots,\;\sigma_{n-1}=\tau
\]
where each $\sigma_{k+1}$ differs from $\sigma_k$ by a transposition and $c(\sigma_k) = n - k$.

Now, it is easy to check that if a transposition $s=\left( ij\right) $ is
applied to a permutation $\lambda ,$ then the number of cycles increases if
and only if $i$ and $j$ belong to the same cycle of $\lambda .$

Consider for example, what happens if we apply $s = (ij)$ to $\tau =\left( 12\ldots
n\right) .$ Then 
\begin{equation*}
s\tau =\left( 12\ldots i-1,j,j+1,\ldots n\right) \left( i,i+1,\ldots
,j-1\right) .
\end{equation*}%
This means that the only block of the partition corresponding to $\tau ,$
was split in two blocks which form the partition corresponding to $s\tau .$ The order of the elements (say, counterclockwise) in the cycles of the corresponding permutation is preserved
Most importantly, these two blocks are non-crossing.

It is easy to see that as we continue applying the transpositions, every transposition must break a cycle, since the sequence of $\sigma_k$ form a geodesic. Then, a
block of the partition that corresponds to a cycle of $\sigma_k$ will be split by a transposition into two blocks that correspond to two new cycles in $\sigma_{k + 1}$ and the non-crossing property of the blocks is preserved.  Hence, by induction the partitions corresponding to every of permutations $
\sigma_i$ are non-crossing.

\begin{figure}[h]
\begin{center}
\includegraphics[width=8cm]{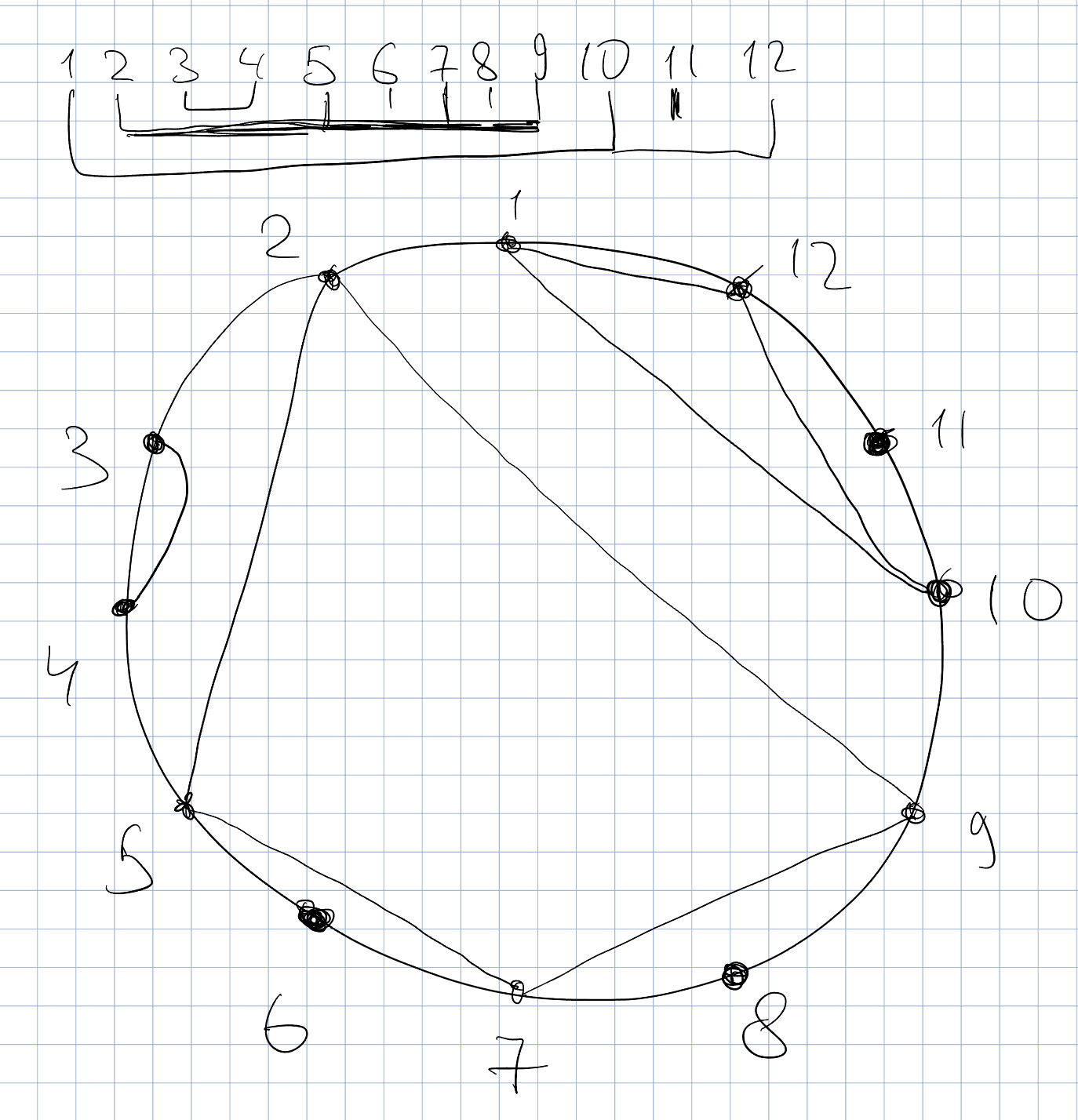}
\end{center}
\caption{Non-crossing partition for permutation $\sigma = (1, 10, 12) (2, 5, 7, 9) (3, 4)$}
\label{fig:non-crossing_on_circle}
\end{figure}

The process can be illustrated using Figure \ref{fig:non-crossing_on_circle}. We can obtain permutation $\sigma$ shown in this figure by applying the sequence of transpositions $(2, 10)$, $(11, 12)$, $(3, 5)$, $(6, 7)$, $(8, 9)$ to $\tau = (1, 2, \ldots, 12)$, so that $\sigma = (8,9) \circ (6, 7) \circ (3, 5) \circ (11, 12) \circ (2, 10) \circ \tau$.


\begin{figure}[h]
\begin{center}
\includegraphics[width=8cm]{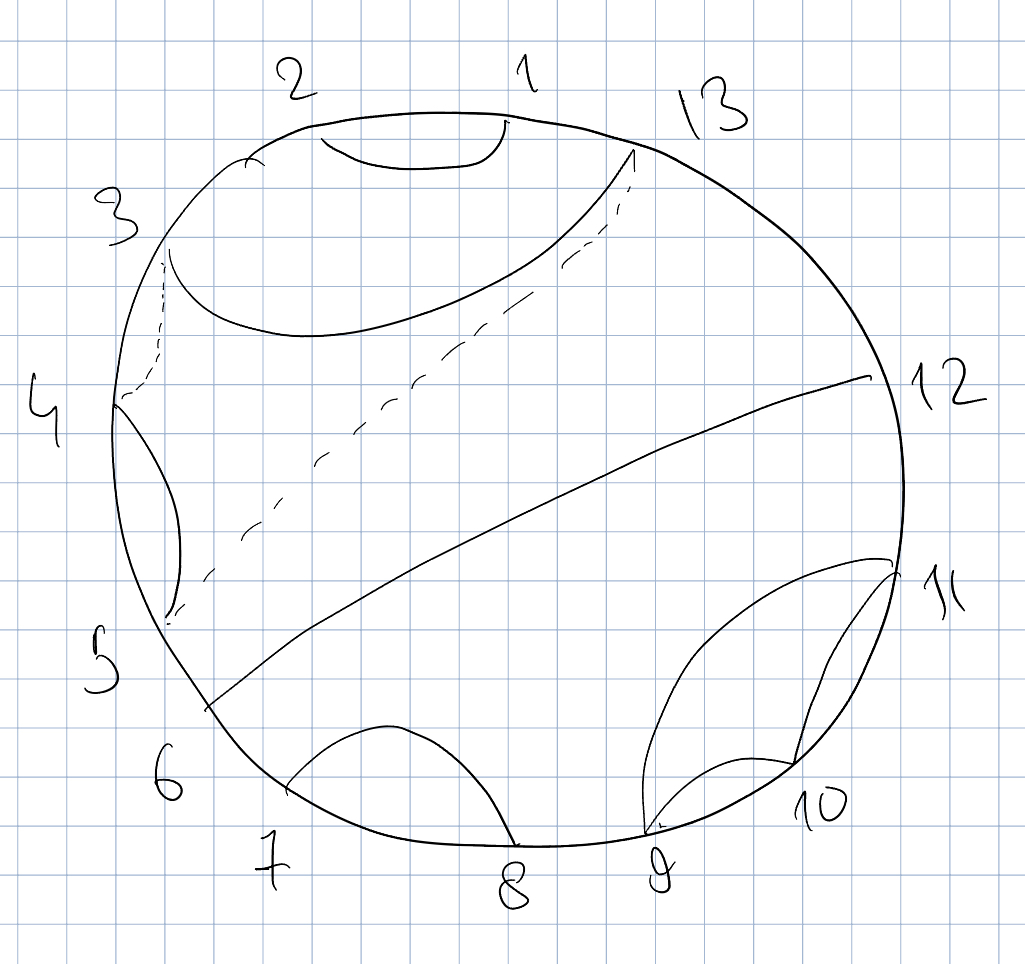}
\end{center}
\caption{}
\label{fig:non-crossing_on_circle2}
\end{figure}

Conversely, take a non-crossing partition $\pi$ of $\{1,\dots,n\}$ with corresponding permutation $\sigma$. We claim that we can join two cycles of $\sigma$ together by multiplying $\sigma$ by a transposition on the right so that the non-crossing property is preserved.   We can start with an \emph{innermost} block $b$—that is, a block forming a contiguous interval around the circle with no other blocks inside it. For example, let this block be $i, i+1, \ldots, j$ (if the block covers $n$, then the argument is similar). Then take the transposition $(i - 1, j)$. Then $\sigma \circ (i - 1, j)$ sends $i - 1$ to $i$ and $j$ to $\sigma(i - 1)$. It preserves correct order of vertices in the cycle, and the connection $[i - 1, i]$ obviously does not create a crossing. In addition, the connection $j, \sigma(i - 1)$ also does not create a crossing because otherwise it would be a crossing with either a pair $i - 1, \sigma(i -1)$ or a crossing with a pair $(i, j)$ which is impossible by our assumption that $\pi$ is not crossing. Hence the new permutation $\sigma \circ (i - 1, j)$ has smaller number of cycles and corresponds to a non-crossing partition. 
This process is illustrated in Figure \ref{fig:non-crossing_on_circle2} with innermost block $b = [4, 5]$.

Repeating this argument, we can connect $\sigma$ to a permutation with one cycle that corresponds to a non-crossing partition, which is $\tau = (1, 2, \ldots, n)$. We can do it by a sequence of transposition that always increase the number of cycles.

 In addition, by an argument above we can connect $\sigma$ to $e$ by breaking the cycles. It follows that $\sigma$ is on the geodesic from $e$ to $\tau$.

\end{proof}

It is equally straightforward to show that for $\lambda$ on a geodesic between some permutations $\sigma$ and $\tau$, the partition corresponding to $\lambda$ refines that of $\sigma$ appropriately. In essence, fewer splits of blocks are required to go from $\sigma$’s partition to $\lambda$’s than from $\sigma$’s partition to $e$’s partition. Consequently, the partial order $\sigma\prec \lambda$ (when $\lambda$ is on a geodesic from $\sigma$ to $\tau$) aligns with the usual refinement order of partitions. 

Altogether, these facts establish a lattice isomorphism between the set of geodesic permutations from $e$ to the $n$-cycle $\tau=(1\,2\,\dots\,n)$ and the lattice of non-crossing partitions on $\{1,2,\dots,n\}$. Through this isomorphism, one can also transfer various operations on non-crossing partitions—such as the Kreweras complement (see Exercise \ref{exercise_Kreweras_complement})—into operations on the corresponding geodesic permutations.

%
\section{Non-cumulant proof of the additivity of the R-transform}

\index{R-transform!additivity}
Let us define $R_X(z)$ of a bounded random variable $X$ using \eqref{equ_defi_Rtransform}: $R_X(z) = G_X^{(-1)}(z) - 1/z$, where $z \in \bC$ is sufficiently small. The Cauchy transform is defined as $G_X(z) = \phi\big((z - X)^{-1}\big)$, for all sufficiently large $z \in \bC$. In particular, $G_X(z)$ is determined by moments of $X$.  

\begin{theo}[\cite{voiculescu86}, \cite{haagerup97}]
\label{theo_additivity_R_2}
Let $X$ and $Y$ be free bounded random variables in a noncommutative probability space. Then the $R$-transform is additive for free convolution, i.e.
\[
R_{X+Y}(z) \;=\; R_X(z) \;+\; R_Y(z).
\]
\end{theo}

\begin{proof}
The functions $R_{X}(z)$ and $R_{Y}(z)$ determine the Cauchy transforms of $X$ and $Y$ and, therefore, they determine the moments of $X$ and $Y$. Since $X$ and $Y$ are assumed to be free, these functions also determine the moments (and hence the Cauchy transform) of $X+Y$. This shows that the function $R_{X+Y}(z)$ is determined by $R_{X}(z)$ and $R_{Y}(z)$. Our strategy is to construct $X$ and $Y$ on a convenient probability space so that they have the prescribed $R$-transforms, and then verify directly that
\[
R_{X+Y}(z) \;=\; R_X(z) + R_Y(z).
\]
Since the choice of $X$ and $Y$ did not rely on any special properties other than these transforms, it follows that the additivity formula holds in general.

\medskip

\index{Fock space}
\noindent
\textbf{Step 1. Construction and basic lemma.} 

Let us work on the full Fock space $F(H)$ with $\dim H \ge 2$. (See Example \ref{exa_Fock_space} and Exercise \ref{exe_Fock_space} on pages  \pageref{exa_Fock_space} and \pageref{exe_Fock_space}, respectively, for details.) Denote by $a^*$ and $a$ the creation and annihilation operators, respectively, corresponding to a unit vector $e \in H$. Define a random variable
\begin{equation}
\label{equ_operator_X}
X \;=\; a^* \;+\; \sum_{k=0}^{\infty} c_k \,a^k
\end{equation}
and set
\[
K_{X}(z) \;:=\; \frac{1}{z} \;+\; \sum_{k=0}^{\infty} c_k \,z^{k}.
\]
We claim that $K_{X}(z)$ is precisely the functional inverse of the Cauchy transform $G_X(z)$, i.e.\ $G_X(K_X(z)) = z$. This is the content of the following lemma.

\begin{lemma}\label{lemma_K_G}
For the above $X$, we have 
\[
G_X\bigl(K_{X}(z)\bigr) \;=\; z.
\]
\end{lemma}

\begin{proof}[Proof of Lemma \ref{lemma_K_G}]
We must show
\[
\phi\Bigl[\bigl(K_{X}(z) - X\bigr)^{-1}\Bigr] \;=\; z.
\]
Consider the vector
\[
w \;=\; \Omega \;+\; \sum_{n=1}^\infty z^n \,e^{\otimes n},
\]
where $\Omega$ is the vacuum vector in $F(H)$. Observe the actions:
\[
a\,w \;=\; z\,w,
\qquad
a^*\,w \;=\;\frac{w - \Omega}{z}.
\]
Hence,
\[
X\,w 
\;=\; \Bigl(a^* + \!\sum_{k=0}^\infty c_k\,a^k\Bigr)w 
\;=\; \frac{w-\Omega}{z} \;+\; \Bigl(\sum_{k=0}^\infty c_k z^k\Bigr)\,w
\;=\; K_X(z)\,w \;-\;\frac{\Omega}{z}.
\]
Rearranging gives 
\[
\bigl(K_{X}(z) - X\bigr)^{-1}\,\Omega 
\;=\; z\,w,
\]
so
\[
\phi\Bigl[\bigl(K_{X}(z) - X\bigr)^{-1}\Bigr]
\;=\;
\langle \Omega,\,z\,w\rangle
\;=\;
z.
\]
This proves the lemma.
\end{proof}

In particular, the lemma shows that for a given function $R(z) = \sum_{k = 0}^\infty c_k z^k$, the random variable $X$ in \eqref{equ_operator_X} has this function as an $R$-transform. 

\medskip

\noindent
\textbf{Step 2. Freely independent operators and additivity.}

Now let $(a^*, a)$ and $(b^*, b)$ be creation and annihilation operators corresponding to two orthogonal vectors $e,f \in H$. Orthogonality of $e$ and $f$ ensures that the subalgebras generated by $(a^*,a)$ and $(b^*,b)$ are free. For given functions $R_X(z) = \sum_{k = 0}^\infty x_k z^k$ and $R_Y(z) = \sum_{k = 0}^\infty y_k z^k$, define random variables
\[
\tilde X \;=\; a^* \;+\; \sum_{k=0}^{\infty} x_k\,a^k,
\quad
\tilde Y \;=\; b^* \;+\; \sum_{k=0}^{\infty} y_k\,b^k.
\]
We know $\tilde X$ and $\tilde Y$ are free. We want to show that
\[
\phi\Bigl[\bigl(K_{\tilde X}(z) + K_{\tilde Y}(z) - \tfrac{1}{z} - (\tilde X+ \tilde Y)\bigr)^{-1}\Bigr]
\;=\;
z.
\]
Consider the vector
\[
w \;=\; \Omega \;+\; \sum_{n=1}^\infty z^n \,(e+f)^{\otimes n}.
\]
By checking the action of $a$ and $b$ on $(e+f)^{\otimes n}$ (in the full Fock space), one finds
\[
a\,w \;=\; z\,w,
\quad
b\,w \;=\; z\,w,
\quad
(a^* + b^*)\,w 
\;=\; \frac{w - \Omega}{z}.
\]
Hence,
\begin{align*}
(\tilde X+\tilde Y)\,w  &= \frac{w - \Omega}{z} + \Big(\sum_{k = 0}^\infty (x_k + y_k) z^k\Big) w \\
& =
\Bigl(K_{\tilde X}(z) + K_{\tilde Y}(z) - \tfrac{1}{z}\Bigr)\,w 
\;-\;\frac{\Omega}{z}.
\end{align*}
It follows that
\[
\bigl[K_{\tilde X}(z) + K_{\tilde Y}(z) - \tfrac{1}{z} - (\tilde X+ \tilde Y)\bigr]^{-1}\,\Omega 
\;=\;
z\,w,
\]
and therefore
\[
\phi\Bigl[\bigl(K_{\tilde X}(z) + K_{\tilde Y}(z) - \tfrac{1}{z} - (\tilde X+ \tilde Y)\bigr)^{-1}\Bigr]
\;=\;
z.
\]
It follows that 
\[
K_{\tilde X}(z) + K_{\tilde Y}(z) \;-\;\frac{1}{z}
\]
is the inverse of $G_{\tilde X+ \tilde Y}(z)$,  and therefore 
\begin{align}
R_{\tilde X+ \tilde Y}(z) &:= G_{\tilde X+ \tilde Y}^{(-1)}(z) - \frac{1}{z} = K_{\tilde X}(z) - \frac{1}{z} + K_{\tilde Y}(z) - \frac{1}{z}
\\
&\;=\;
R_{\tilde X}(z) + R_{\tilde Y}(z),
\end{align}
where we used Lemma \ref{lemma_K_G} in the last step. By the argument outlined in the beginning of the proof this implies that $R_{X + Y}(z) = R_X(z) + R_Y(z)$. 
This completes the proof.
\end{proof}

\section*{Exercises}


\begin{exe}
Let $\AA_1, \ldots, \AA_k$ be free subalgebras of $\AA$ and let $x_1, \ldots, x_n$ be elements of $\AA$, such that $x_i \in A_{r(i)}$. Consider the partition $\pi$ of ordered set $[n]$ with $k$ blocks, in which element $i$ belongs to block $b_j$ if $r(i) = j$. (That is, block $b_j$ consists of indices of all elements that belong to sub-algebra $\AA_j$ according to assignment $r(\cdot)$. Some blocks can be empty.) Use the properties from Example \ref{exa_prod_rule1} to show that \emph{if $\pi$ is non-crossing} then 
\bal{
\phi(x_1 \ldots x_n) = \prod_{j=1}^{|\pi|} \phi\big(x_{b_j(1)} \ldots x_{b_j(|b_j|)}\big),
}
where $b_j(1), \ldots, b_j(|b_j|)$ are elements of the block $b_j$.

For example, if $x_1, x_4, x_6, x_8 \in A_1$, $x_2, x_3 \in A_2$, $x_5 \in A_3$ and $x_7 \in A_4$, then 
\bal{
\phi(x_1 x_2 \ldots x_8) = \phi(x_1 x_4 x_6 x_8) \phi(x_2 x_3) \phi(x_5) \phi(x_7).
} 

\end{exe}

\begin{exe}
\label{multMu}
Suppose that $\lambda \leq
\sigma $ ($\lambda$ is more refined then $\sigma$) and that $\sigma =\sigma _{1}\sigma _{2}.$ Then, we can factorize $%
\lambda $ as $\lambda _{1}\lambda _{2}$ so that $\lambda _{1}\leq \sigma
_{1}$ and $\lambda _{2}\leq \sigma _{2}.$ Show that $\mu \left( \lambda ,\sigma
\right) =\mu \left( \lambda _{1},\sigma _{1}\right) \mu \left( \lambda
_{2},\sigma _{2}\right) .$
\end{exe}

\index{Kreweras complement!for permutations}
\begin{exe}
\label{exercise_Kreweras_complement} Show that if $\pi $ is a geodesic
permutation (i.e., it belongs to a geodesic between $e$ and $\tau =\left(
12\ldots n\right) ,$ then its Kreweras complement can be computed by the
formula $K\left( \pi \right) =\pi ^{-1}\tau .$
\end{exe}

\section*{Notes}
For more information on M\"obius inversion theory, the reader can consult Chapter IV of
\emph{Combinatorial Theory} by \cite{aigner1979} or Sections 3.6--3.10 in \emph{Enumerative Combinatorics, volume 1} by  \cite{stanley2012}. Most of the results in this Lecture are from \cite{nica_speicher06}.


\chapter{Additive free convolution and limit theorems}

\section{Additive free convolution}

\index{free convolution!additive}
Let \(\mu\) and \(\nu\) be two compactly supported probability measures. We can find two \emph{self-adjoint} operators \(X\) and \(Y\) in a tracial $W^\ast$-probability space whose spectral distribution measures are \(\mu\) and \(\nu\), respectively. By choosing an appropriate non-commutative probability space, we can ensure that \(X\) and \(Y\) are free. Then Theorem~\ref{theo_additivity_R} implies that the moments of \(X + Y\) depend only on \(\mu\) and \(\nu\). The probability measure determined by these moments is called the \emph{free additive convolution}\index{free additive convolution}\index{convolution!free additive} of \(\mu\) and \(\nu\), and we denote it by
\[
   \mu \boxplus \nu.
\]

This definition extends to probability measures with unbounded support via a truncation method: one approximates each measure by compactly supported (truncated) versions, computes their free convolutions, and then takes a suitable limit. The details of this construction are somewhat technical, and we omit them here. The interested reader can find a thorough treatment in \cite{bercovici_voiculescu93}.

\smallskip

Recall that the \emph{classical convolution}\index{convolution!classical} of measures \(\mu\) and \(\nu\) is the distribution of \(X + Y\) when \(X\) and \(Y\) are independent; we denote it by
\[
   \mu \ast \nu.
\]

\begin{exe}
Let \(\mu = \delta_x\). Show that \(\mu \boxplus \nu\) is precisely the shift of \(\nu\) by \(x\), i.e.
\[
   \mu \boxplus \nu \,(A) \;=\; \nu\,(A - x).
\]
\end{exe}

In this situation, the free additive convolution coincides with the classical one:
\[
   \delta_x \,\boxplus\, \nu \;=\; \delta_x \,\ast\, \nu.
\]
In general, however, the free additive convolution \(\boxplus\) is very different from the classical convolution \(\ast\). A key point of departure is \emph{non-linearity}. Recall that a linear combination of measures \(\mu\) and \(\nu\) is defined as
\[
   \bigl(t\,\mu \;+\; s\,\nu\bigr)\,(A)
   \;=\;
   t \,\mu(A) \;+\; s \,\nu(A),
\]
for any measurable set \(A\).

\smallskip

When \(t + s = 1\), the classical convolution is linear in these convex combinations:
\[
   \bigl(t\,\mu_1 + s\,\mu_2\bigr) \,\ast\, \nu
   \;=\;
   t \bigl(\mu_1 \ast \nu\bigr)
   \;+\;
   s \bigl(\mu_2 \ast \nu\bigr).
\]
Since we can easily convolve a single atom with another measure, this linearity property makes it straightforward to compute \(\mu \ast \nu\) for discrete measures (those supported on finitely many points), and the result remains discrete and finitely supported.

\smallskip

In contrast, free additive convolution does \emph{not} satisfy such linearity:
\[
   \bigl(t\,\mu_1 + s\,\mu_2\bigr)\,\boxplus\, \nu
   \;\;\neq\;\;
   t\,(\mu_1 \boxplus \nu)
   \;+\;
   s\,(\mu_2 \boxplus \nu).
\]
For example, let
\[
   \mu \;=\; \nu
   \;=\;
   \tfrac12 \,\bigl(\delta_{-1} + \delta_{1}\bigr).
\]
Then (as shown in Exercise~\ref{exe_bernoulli}) the free additive convolution \(\mu \boxplus \nu\) is the arcsine distribution on \((-2,2)\). That is, although \(\mu\) and \(\nu\) are purely atomic, \(\mu \boxplus \nu\) is absolutely continuous.

\smallskip

On the other hand, the free additive convolution does share a continuity property with the classical convolution. Specifically, let \(\mu_n \to \mu\) and \(\nu_n \to \nu\) weakly, meaning that
\[
   \int f \,d\mu_n \;\longrightarrow\; \int f \,d\mu,
   \quad
   \int f \,d\nu_n \;\longrightarrow\; \int f \,d\nu,
\]
for every bounded continuous \(f\). Then
\[
   \mu_n \,\boxplus\, \nu_n
   \;\longrightarrow\;
   \mu \,\boxplus\, \nu
\]
weakly. We refer to~\cite{bercovici_voiculescu93} for a proof. Furthermore, if \(\mu_n\) and \(\nu_n\) have unbounded support, then under suitable assumptions the convergence above is \emph{tight} and $\mu \boxplus \nu$ is a valid probability measure.

\section{Univariate free CLT}

The following result is the analogue of the classical Central Limit Theorem for sums of independent, identically distributed, centered variables in the setting of free probability.

\index{CLT!univariate}
\begin{theo}
Let \(X_{1},X_{2},\ldots\) be a sequence of identically distributed, bounded, self-adjoint random variables. Assume that \(\phi(X_{i})=0,\) \(\phi(X_{i}^{2})=1,\) and that the \(X_i\) are free. Define 
\[
S_{n} \;=\; X_{1} + \cdots + X_{n}.
\]
Then the sequence \(S_{n}/\sqrt{n}\) converges in distribution to the standard semicircle random variable.
\end{theo}

\begin{proof} Since \(\phi(X_i) = 0\) and \(\phi(X_i^2) = 1\), we have
 the \(R\)-transform of each \(X_i\) equal to 
\[
R(z) 
\;=\; 
z \;+\; c_{2} z^{2} \;+\; c_{3} z^{3} \;+\;\dots
\]

By freeness and additivity of the \(R\)-transform, the \(R\)-transform of \\ 
\(\displaystyle S_{n} = \sum_{i=1}^n X_{i}\) is
\[
R_{n}(z) 
\;=\;
n \,R(z).
\]
Then, using the usual scaling property of the \(R\)-transform, we get
\[
R_{S_{n}/\sqrt{n}}\bigl(z\bigr) 
\;=\; 
\frac{1}{\sqrt{n}}\,R_{n}\!\Bigl(\tfrac{z}{\sqrt{n}}\Bigr) 
\;=\;
\sqrt{n}\;R\!\Bigl(\tfrac{z}{\sqrt{n}}\Bigr) 
\;=\;
z \;+\; \frac{c_{2}}{\sqrt{n}}\,z^{2} 
\;+\; \frac{c_{3}}{n}\,z^{3} 
\;+\;\dots
\]
As \(n\to\infty,\) all coefficients beyond the linear term tend to zero, so in the limit,
\[
R_{\mathrm{limit}}(z) 
\;=\; z.
\]
Since a compactly supported distribution is determined by its \(R\)-transform, this means that the moments of \(S_{n}/\sqrt{n}\) converge to the moments of the (unique) distribution whose \(R\)-transform is
$z$, i.e. the standard semicircle distribution. That is, 
\(S_{n}/\sqrt{n}\) converges in distribution to the standard semicircle random variable.
\end{proof}

This theorem prompts several natural questions, parallel to those arising in classical probability:

\begin{enumerate}
\item Can the Central Limit Theorem be extended to unbounded variables, or to variables that are not necessarily identically distributed?
\item Is there a corresponding theory of infinitely divisible distributions in free probability?
\item What can be said about large deviations for the sums \(S_{n}\)?
\end{enumerate}

These questions have been investigated extensively in the literature; we provide only a brief discussion. An excellent overview of results in these areas can be found in \cite{bercovici_pata99}.

\bigskip

Regarding large deviations, the following inequality from \cite{voiculescu86} is particularly relevant:
\[
\Bigl\|\sum_{i=1}^{r} X_{i}\Bigr\|
\;\le\;
\max_{1 \le i \le r} \|X_{i}\|
\;+\;
\sqrt{\sum_{i=1}^{r}\mathrm{Var}\bigl(X_{i}\bigr)}.
\]
If \(X_{i}\) are free, identically distributed, bounded, centered random variables with unit variance, then
\[
\Bigl\|\tfrac{S_{n}}{\sqrt{n}}\Bigr\|
\;\le\;
1 
\;+\;
\frac{\|X_{i}\|}{\sqrt{n}},
\]
showing that for each \(\varepsilon>0,\) the support of \(S_{n}/\sqrt{n}\) is eventually contained in \(\bigl[-1-\varepsilon,\;1+\varepsilon\bigr]\). Together with the high smoothness of the distribution of \(S_{n}/\sqrt{n},\) this phenomenon was termed \emph{superconvergence} in \cite{bercovici_voiculescu95}.

For unbounded variables, the question of large deviations for free sums remains only partially resolved.

%

\section{Free Poisson Limit}

We now discuss a free analogue of another result from classical probability theory, sometimes referred to as the law of small numbers. In the classical setting, this law asserts that counts of rare events follow the Poisson law (the classical example is the distribution of the number of deaths of Prussian cavalry officers from horse kicks).

First, let us define the free analogue of the Poisson law. Let \(\mu\) be a distribution with density
\[
p(x) \;=\; \frac{\sqrt{4x - \bigl(1 - \lambda + x\bigr)^2}}{2\pi\,x}
\quad\text{if}\quad x \in \bigl[(1 - \sqrt{\lambda})^2,\,(1 + \sqrt{\lambda})^2 \bigr],
\]
where \(\lambda\) is a positive parameter. If \(x\) is outside this interval, the density is zero. Moreover, if \(\lambda < 1\), then \(\mu\) has an atom at \(0\) of weight \(1 - \lambda\).

This distribution is called the \emph{free Poisson distribution}%
\index{free Poisson distribution}%
\index{distribution!free Poisson}%
\index{free Poisson r.v.|textbf}
\label{definition_free_Poisson}
with parameter \(\lambda\). It is also known as the \emph{Marchenko-Pastur distribution}%
\index{Marchenko-Pastur distribution}%
\index{distribution!Marchenko-Pastur}
because it was first discovered in~\cite{marchenko_pastur67}.

Recall that the classical Poisson distribution is supported on the set of non-negative integers, given by 
\[
\mu(\{k\}) 
\;=\; e^{-\lambda}\,\frac{\lambda^k}{k!}.
\]
In contrast, the free Poisson distribution is absolutely continuous except for a possible atom at \(0\). The reason it is called the free Poisson distribution is shown by the following theorem. If we replaced “freeness” with “independence” in its hypotheses, we would recover the classical Poisson distribution.

\index{Poisson convergence|textbf}
\begin{theo}
Let \(X_{1,n},\ldots,X_{n,n}\) be self-adjoint random variables each having the Bernoulli distribution
\[
\mu 
\;=\; \Bigl(1 - \frac{\lambda}{n}\Bigr)\delta_{0} 
\;+\; \frac{\lambda}{n}\,\delta_{1}.
\]
Assume that \(X_{1,n},\ldots,X_{n,n}\) are free and define
\[
S_{n} 
\;=\; X_{1,n} + \ldots + X_{n,n}.
\]
Then \(S_{n}\) converges in distribution to the free Poisson distribution with parameter~\(\lambda\).
\end{theo}

\begin{proof} One may compute the Cauchy transform of \(X_{i,n}\) as
\[
G(z) 
\;=\; \frac{z - 1 + \lambda/n}{(z - 1)\,z},
\]
and from there, the \(R\)-transform is
\begin{equation}
\label{Bernoulli_Rtransform}
R(z) 
\;=\; \frac{-1 + z \;-\; \sqrt{\bigl(1 - z\bigr)^2 \;+\; 4\,(\lambda/n)\,z}}{2\,z}.
\end{equation}
Hence, the \(R\)-transform of the sum \(S_{n}\) is
\[
R_{n}(z) 
\;=\; \frac{1 - z}{2\,z}\,n 
\Bigl(-1 + \sqrt{\,1 \;+\; \frac{4\,\lambda}{n}\,\frac{z}{(1 - z)^2}}\Bigr).
\]
It follows that
\[
R_{n}(z) \;\longrightarrow\; \frac{\lambda}{1 - z},
\]
and this convergence is uniform in a sufficiently small disk around \(z = 0\).

\index{free Poisson r.v.!R-transform}
It is straightforward to verify that \(\lambda\,(1 - z)^{-1}\) is the \(R\)-transform of the free Poisson distribution. Since the convergence of \(R\)-transforms implies convergence of the Cauchy transforms (and thus of the moments), we conclude that \(S_{n}\) converges in distribution to the free Poisson law with parameter~\(\lambda\). 
\end{proof}

\chapter{Asymptotic Freeness}

Free probability is closely related to the theory of large random matrices. In this chapter we sketch this connection.

\section{Gaussian matrices}


Let $A^{(N)}=(a_{ij})_{i,j=1}^N$ be an $N\times N$ Hermitian Gaussian random matrix (the GUE).  That is, $a_{ji}=\overline{a_{ij}}$, the entries 
$$\{\Re a_{ij},\,\Im a_{ij}:i<j\}\quad\text{and}\quad\{a_{ii}\}$$ 
are independent Gaussians with zero mean, and
\[
\mathbb{E}[\,a_{ij}\,]=0,\qquad
\mathbb{E}[\,a_{ij}\,a_{kl}\,]=0,\qquad
\mathbb{E}[\,a_{ij}\,\overline{a_{kl}}\,]
=\frac{1}{N}\,\delta_{ik}\,\delta_{jl}.
\]
The last equality can also be written as 
\begin{equation}
\E (a_{ij}a_{kl}) = \E (a_{ij} \overline a_{lk}) = \delta _{il}\delta _{jk}N^{-1}.
\label{formula_af_semicircle_covariance}
\end{equation}%
Equivalently:
\[
\Re a_{ij},\;\Im a_{ij}\sim N\bigl(0,\tfrac1{2N}\bigr)\quad(i<j),
\qquad
a_{ii}\sim N\bigl(0,\tfrac1N\bigr).
\]


Let $D^{( N) }$ be a sequence of $N$-by-$N$ (non-random) Hermitian
matrices. We think about $A^{( N) }$
 $D^{( N) }$ as elements of non-commutative probability space $(M_N(\bC), \tr \otimes \E)$.

Suppose that as $N \to \infty$,  $D^{\left( N\right) }$ converges in distribution to a non-commutative random variable $d$ in a space $\AA_1, \phi_1$. Here the convergence in distribution is
understood as in Definition \ref{definition_convergence_distribution}. Essentially this means that the empirical distribution of eigenvalues of 
$D^{( N) }$ weakly converges to a probability distribution $\mu$. 

It is known from the theory of random matrices that the empirical eigenvalue distribution of the GUE matrices weakly converges to the semicirle distribution. Hence, matrices  $A^{( N) }$ converge in distribution (in the sense of non-commutative r.v.s) to a semicircle random variable $s$ in $(\AA_2, \phi_2$. By taking the free product $(\AA, \phi) = (\AA_1, \phi_1) \ast (\AA_2, \phi_2)$, we can assume that $d$ and $s$ belong to the same non-commutative probability space and free. 

\index{asymptotic freeness|textbf}
\begin{theo}
\label{theorem_af_semicircle}The sequence of $\left( A^{\left( N\right)
},D^{\left( N\right) }\right) $ converges in distribution to $\left(
s,d\right) $ where $s$ has the semicircle distribution and $s$ and $d$ are
free.
\end{theo}

It is not really necessary that $D^{\left( N\right) }$ are non-random. They can be random but in this case they must be independent of the Gaussian matrices $A^{\left( N\right) }.$ The assumption in this case that the (random) empirical eigenvalue distribution of $D^{( N) }$ converges to a probability measure $\mu$ almost surely. 

This phenomenon of the convergence of independent random matrices to free random variables is often called \emph{asymptotic freeness of random matrices}. 

First, let us establish a useful result about the moments of the semicircle variable. 
\begin{propo}
\label{propo_moments_semicircle}
 Let $s, d_1,d_2,\ldots ,d_n$ be elements of a non-commutative probability space $(\AA, \phi)$. Assume that $s$ is a semicircle element which is free from $\{d_1,d_2,\ldots ,d_n\}$. Then, 
\begin{equation}
\phi ( sd_1 s d_2 \ldots s d_n) =\sum_{\pi \in NP_{2}(n)}\phi_{\pi \gamma }( d_1,d_2,\ldots ,d_n),
\label{formula_af_semicircle2}
\end{equation}
\end{propo}
The notation $\pi \in NP_{2}( n) $ means that $\pi $
is a permutation of set $\{ 1,2,\ldots, n\}$ that corresponds to a non-crossing pairing. This means that (i) $\pi $ is a product of
disjoint transpositions with no fixed points, and (ii) the pairing corresponding to these transpositions is
non-crossing. The first requirement implies, of course, that $n$ must be
even. The meaning of the second requirement is that there are no
transpositions $( ij) $ and $( kl) $ in the product
such that $i<k<j<l.$
Then, $\lambda$ is a special permutation defined as
\begin{equation}
\gamma :=( 12\ldots n) ,  \label{definition_af_gamma}
\end{equation}%
and $\pi \gamma $ is the product of permutations $\pi $ and $\gamma .$

\begin{proof}[Proof of Proposition \ref{propo_moments_semicircle}]
We apply Theorem \ref{theorem_expectation_products} and write
\[
\phi ( sd_1 s d_2 \ldots s d_n)  = \sum_{\pi \in NC(n)} \kappa_\pi(s, \ldots, s) \phi_{K(\pi)}(d_1, \ldots, d_n). 
\]
The cumulant $\kappa_\pi(s, \ldots, s)$ is not zero if and only if $\pi$ is a non-crossing pairing, and in this case it equals $1$. Therefore,  
\[
\phi ( sd_1 s d_2 \ldots s d_n)  = \sum_{\pi \in NC_2(n)} \phi_{K(\pi)}(d_1, \ldots, d_n) = \sum_{\pi \in NC_2(n)} \phi_{\pi \gamma}(d_1, \ldots, d_n),
\]
where the last equality follows from Exercise \ref{exercise_Kreweras_complement}. (Indeed, $\pi$ is an involution, $\pi^{-1} = \pi$, and, as a product of disjoint transpositions, it belongs to the geodesic between $e$ and $\lambda$.)  
\end{proof}

\begin{proof}[Proof of Theorem \ref{theorem_af_semicircle}] \ Let $D^{1},D^{2},\ldots ,D^{n}$ denote arbitrary
polynomials of the matrix $D^{\left( N\right) },$ and let us for conciseness omit the superscript in $A^{( N) }.$ We wish to compute $\E \tr (AD^{1}AD^{2}\ldots AD^{n}).$ Note that $D^{i}$ can be equal to the
identity operator and therefore the formula for $\E \tr ( AD^{1}AD^{2}\ldots AD^{n})$ will cover the expectations of the products which include powers of $A.$

We are going to show that 
\begin{equation}
\E \tr ( AD^{1}AD^{2}\ldots AD^{n})  -\sum_{\pi \in NP_{2}\left( n\right) }\tr_{\pi \gamma}( D^{1},D^{2},\ldots ,D^{n}) \to 0
\label{formula_af_semicircle1}
\end{equation}
as $N\,$approaches infinity. Here $\tr_{\pi \gamma }\left( X^{1},\ldots ,X^{n}\right) $ denotes the ``cyclic trace'' of matrices $X^{1},\ldots ,X^{n}.$ Namely, suppose that
permutation $\pi \gamma $ equals the product of $s$ cycles $c_{i},$ $%
i=1,\ldots ,s,$ and let the elements of cycle $c_{i}$ be denoted $%
c_{i,1},c_{i,2},\ldots ,c_{i,k_{i}}.$ Then 
\begin{align*}
\mathrm{tr}_{\pi \gamma }\left( X^{1},\ldots ,X^{n}\right)
&:=\prod_{i=1}^{s}\mathrm{tr}\left( X^{c_{i,1}}X^{c_{i,2}}\ldots
X^{c_{i,k_{i}}}\right) \\
&=N^{-s}\prod_{i=1}^{s}\mathrm{Tr}\left( X^{c_{i,1}}X^{c_{i,2}}\ldots
X^{c_{i,k_{i}}}\right) .
\end{align*}%
For example, if $n=6$ and $\pi =\left( 14\right) \left( 23\right) \left(
56\right) ,$ then $\pi \gamma =\left( 13\right) \left( 2\right) \left(
46\right) \left( 5\right) ,$ and 
\begin{equation*}
\mathrm{tr}_{\pi \gamma }\left( X^{1},\ldots ,X^{6}\right) =N^{-4}\mathrm{Tr}%
\left( X^{1}X^{3}\right) \mathrm{Tr}\left( X^{2}\right) \mathrm{Tr}\left(
X^{4}X^{6}\right) \mathrm{Tr}\left( X^{5}\right) .
\end{equation*}



A comparison of formulas (\ref{formula_af_semicircle1}) and (\ref%
{formula_af_semicircle2}) shows that the joint moments of matrices $%
A^{\left( N\right) }$ and $D^{i) }$ converge to the
corresponding joint moments of random variables $s$ and $d_i,$ and this is exactly what is needed  to establish the validity of
Theorem \ref{theorem_af_semicircle}.

Let us expand the left-hand side of (\ref{formula_af_semicircle1}):%
\begin{equation}
\E \tr ( AD^{1}AD^{2}\ldots AD^{n}) =%
\frac{1}{N}\sum_{i_{1}j_{1}\ldots }\E (a_{i_{1}j_{1}}d_{j_{1}i_{2}}^{1}a_{i_{2}j_{2}}\ldots
d_{j_{n}i_{1}}^{n}).  \label{formula_af_semicircle25}
\end{equation}%
Since $d$ variables are not random, we can take them outside the expectation $\E$. In addition, we can use the Wick formula to compute the
expectations of the products of Gaussian variables. Then, 
\begin{align*}
\frac{1}{N}\sum_{i_{1}j_{1}\ldots }\E(a_{i_{1}j_{1}}a_{i_{2}j_{2}}\ldots a_{i_{n}j_{n}}) &
d_{j_{1}i_{2}}^{1}\ldots d_{j_{n}i_{1}}^{n}  \\
&=\frac{1}{N}\sum_{i_{1}j_{1}\ldots }\sum_{\pi \in \mathcal{P}_{2}\left(
n\right) }\prod_{s=1}^{n}\left\langle a_{i_{s}j_{s}}a_{i_{\pi \left(
s\right) }j_{\pi \left( s\right) }}\right\rangle ^{1/2}d_{j_{s}i_{s+1}}^{s}.
\end{align*}%
Here, $\mathcal{P}_{2}\left( n\right) $ is the set of all possible pairings
of $n$ elements. (In particular, this set is empty if $n$ is odd.) The power 
$1/2$ in this formula is needed to avoid double counting. Also, by
convention $i_{n+1}:=i_{1}.$

For example, for $n = 2$ we have only one pairing $\pi = (12)$ and, for instance,
\[
\E(a_{29} d^1_{93} a_{37} d_{72}^2) = \Big[\E(a_{29} a_{37})\Big]^{1/2} d_{93}^1 \Big[\E(a_{37} a_{29})\Big]^{1/2} d_{72}^2.
\]


Recall our assumption about the covariance structure of entries, (\ref%
{formula_af_semicircle_covariance}). This assumtion implies that for all non-zero terms
in this sum, we have $i_{s}=j_{\pi \left( s\right) }$ and $j_{s}=i_{\pi
\left( s\right) }.$ This allows us to express everything in terms of $j$%
-indices. Namely, we can re-write the previous expression as 
\begin{align}
\frac{1}{N^{n/2+1}}\sum_{j_{1},j_{2},\ldots ,j_{n}}\sum_{\pi \in \PP_2(n) }\prod_{s=1}^{n}d_{j_{s}j_{\pi ( s+1) }}^{s} 
&= \frac{1}{N^{n/2+1}}\sum_{j_{1},j_{2},\ldots
,j_{n}}\sum_{\pi \in \PP_2(n)}\prod_{s=1}^{n}d_{j_{s}j_{\pi \gamma(s) }}^{s} \notag \\
&=\frac{1}{N^{n/2+1}}\sum_{\pi \in \PP_2(n) }\Tr_{\pi \gamma }( D^{1},\ldots ,D^{n}) \notag \\
&= \sum_{\pi \in \PP_2( n) }N^{\#\left( \pi \gamma
\right) -n/2-1}\tr_{\pi \gamma }( D^{1},\ldots ,D^{n}) \label{formula_af_semicircle3} 
\end{align}
Here $\#\left( \pi \gamma \right) $ denotes the number of cycles (including
trivial) in the permutation $\pi \gamma .$

By assumption, the normalized traces of the polynomials of variables $D^{i}$
converge to certain limits. Hence, in order to find the asymptotic behavior
of the trace we need to find those pairings $\pi $ for which $\#\left( \pi
\gamma \right) $ takes the maximal value.

\begin{lemma}
\label{Lemma_number_cycles}Let $\gamma =\left( 1,2,\ldots ,2m\right) $ and
let $\pi $ be a permutation corresponding to a pairing of the set $\left\{
1,2,\ldots ,2m\right\} .$ Then $\#\left( \pi \gamma \right) \leq m+1,$ and
the equality is achieved if and only if $\pi $ is non-crossing.
\end{lemma}

\begin{proof}
Indeed, $\gamma $ has just one cycle, and a multiplication by a
transposition can increase the number of cycles by no more than 1. Hence, $%
\#\left( \pi \gamma \right) \leq m+1.$

For non-crossing pairings $\pi $, we can show that the number of cycles in $%
\pi \gamma $ is $m+1$ by induction. Indeed, let $\left( ij\right) $ is one
of the outer-most pairs of $\pi ,$ that is, there is no pair $\left(
i^{\prime }j^{\prime }\right) $ such that $i^{\prime }<i$ and $j^{\prime
}>j. $ If we compute $\left( ij\right) \gamma ,$ then we will get two
cycles: $\left( 1,\ldots ,i-1,j,\ldots ,2m\right) $ and $\left( i,\ldots
,j-1\right) . $ The pairs that are outside of $\left( ij\right) $ will
operate on the first cycle, and the pairs that are inside of $\left(
ij\right) $ will operate on the second cycle. By induction hypothesis the
multiplication by these transpositions will always increase the number of
cycles by one, hence the total number of cycles will be $m+1.$

In contrast, if there is a crossing and the pairing $\pi $ contains $\left(
ij\right) (kl)$ with $i<k<j<l,$ then it is easy to check that $\left(
ij\right) (kl)\gamma $ is a cycle. Multiplication by remaining
transpositions can increase the number of cycles by $m-2$ at most. Hence, $%
\#\left( \pi \gamma \right) \leq m-1$ in this case. This completes the proof
of Lemma \ref{Lemma_number_cycles}.
\end{proof}

Lemma \ref{Lemma_number_cycles} shows that the right-hand side in (\ref%
{formula_af_semicircle3}) converges to the sum over non-crossing pairings
only, and therefore, 
\begin{equation*}
\E \tr ( AD^{1}AD^{2}\ldots AD^{n}) \to \sum_{\pi \in \mathcal{NP}_{2}( n) }\tr_{\pi \gamma }( D^{1},\ldots ,D^{n}) .  
\end{equation*}

This completes the proof of the theorem.
\end{proof}

\begin{figure}[tbph]
\begin{center}
\includegraphics[width=6cm]{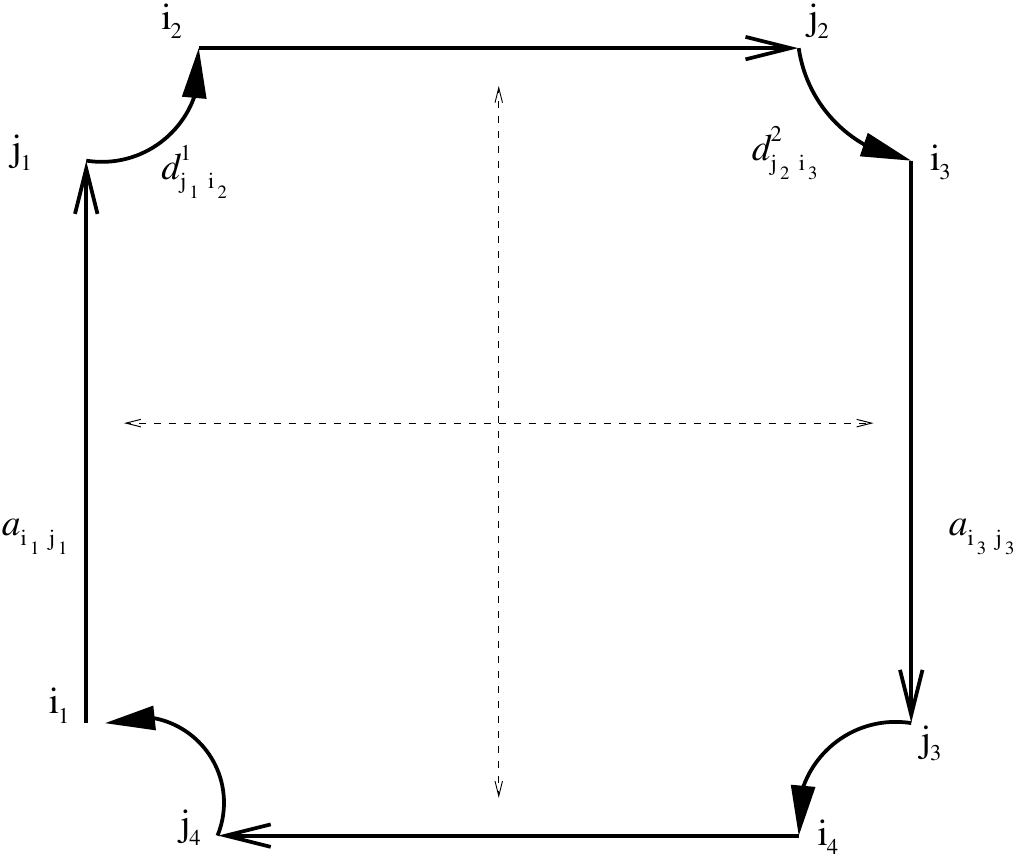}
\end{center}
\caption{Graphical representation of trace. Genus 1.}
\label{figure_genus_expansion1}
\end{figure}

In order to see better what is going on in the proof of the theorem, let us
represent the sum in (\ref{formula_af_semicircle25}) by an oriented polygon
with labeled vertices. See Figure \ref{figure_genus_expansion1}.

The vertex labels can take values from $1$ to $N.$ However, not all
labelings survive after we take expectations.

\index{Wick's formula}
Indeed, Wick's formula allows us to compute the expectation as a sum over
pairings. Each pairing can be thought of as a gluing of the polygon and this
gluing and the assumed covariance structure forces some of the labels to be
identified. For example presented in Figure \ref{figure_genus_expansion1},
we have ,%
\begin{eqnarray*}
i_{1} &=&j_{3},\text{ }i_{2}=j_{4}, \\
i_{3} &=&j_{1},\text{ }i_{4}=j_{2}.
\end{eqnarray*}%
Next, it turns out that arcs that correspond to $D$-variables form cycles
and the number of these cycles equals the number of cycles in the
permuations $\pi \gamma .$

Let us ignore these cycles for a second and think about them as vertices
that remain distinct after we glued the polygon. Then we get a closed
surface and the Euler characteristic of this surface can be computed as 
\begin{equation*}
\chi =F-E+V,
\end{equation*}%
where $F,$ $E,$ and $V$ are the number of faces, edges, and vertices in the
map that we results from the edges and vertices of the polygon after the
gluing. Clearly, $F=1$ and $E=n/2.$ The number of distinct vertices $%
V=\#\left( \pi \gamma \right) .$ Hence, we have 
\begin{equation*}
\chi =1-n/2+\#\left( \pi \gamma \right) .
\end{equation*}%
Recall that the genus of a surface $g$ is related to its Euler
characteristic by the formula $2g=2-\chi .$ Hence, we have%
\begin{equation*}
2g=1+n/2-\#\left( \pi \gamma \right) .
\end{equation*}%
This means that the exact formula for the trace (\ref{formula_af_semicircle3}%
) can be written in the following way:%
\begin{equation*}
\E \tr ( AD^{1}AD^{2}\ldots AD^{n}) =\sum_{\pi \in \PP_2( n) }\frac{1}{%
N^{-2g( \pi ) }}\tr_{\pi \gamma }( D^{1},\ldots,D^{n}) ,  
\end{equation*}%
where $g( \pi ) $ denotes the genus of the closed surface
constructed by gluing the $n$-polygon according to the pairing $\pi .$

\begin{figure}[tbph]
\begin{center}
\includegraphics[width=6cm]{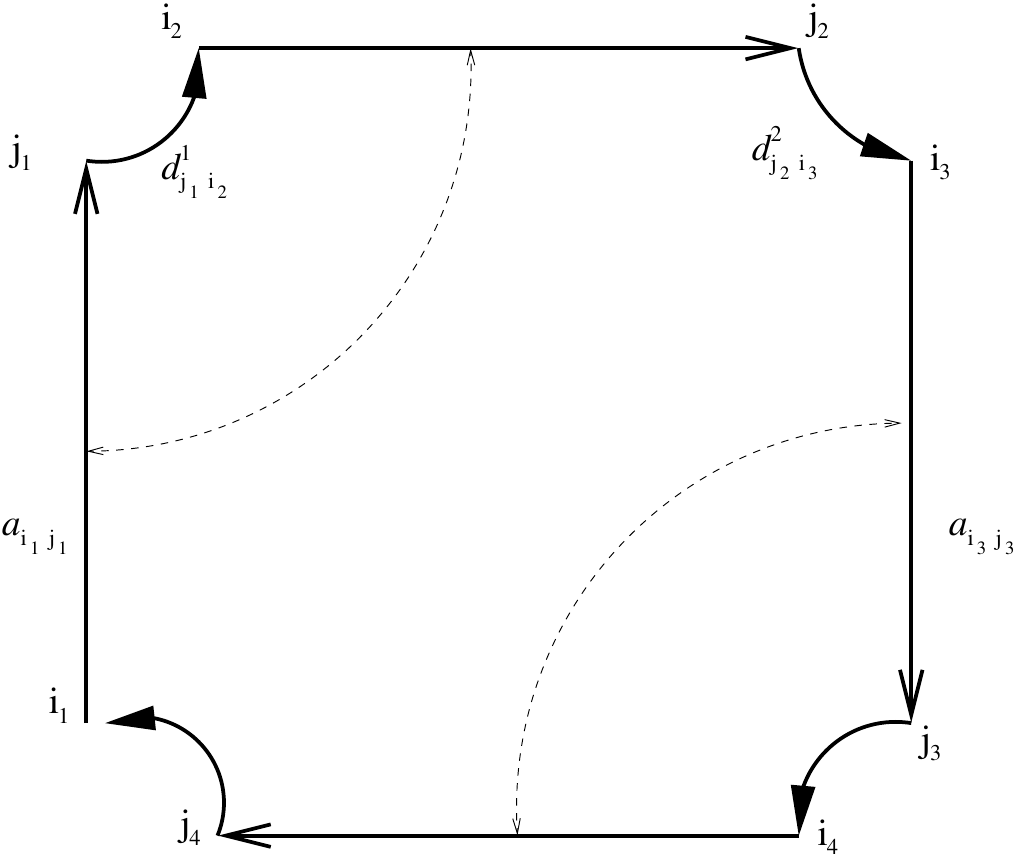}
\end{center}
\caption{Graphical representation of trace. Genus 0.}
\label{figure_genus_expansion2}
\end{figure}

For the gluing in Figure \ref{figure_genus_expansion1}, we have $g=1,$ and
for the gluing in Figure \ref{figure_genus_expansion2}, we have $g=0.$

Formulas of this type are often called the \emph{genus expansion} 
\index{genus expansion} formulas for matrix models.

Let us illustrate this formula with an example in which $D^{i}=I$ for all $%
i. $ Then, 
\begin{equation*}
\E \tr ( A^n) =\sum_{g}k_n( g) N^{-2g},
\end{equation*}%
where $k_n (g) $ is the number of distinct glueings of a
polygon with $n$ sides that induce the surface of genus $g.$ For example,
if $n=2,$ then%
\begin{equation*}
\E \tr ( A^2)  =1;
\end{equation*}%
if $n=4,$ then 
\begin{equation*}
\E \tr ( A^{4})  =2+N^{-2};
\end{equation*}%
if $n=6,$ then 
\begin{equation*}
\E \tr ( A^4)  =5+10 N^{-2},
\end{equation*}%
and so on.

\bigskip

It is clear from the proof of the theorem that instead of one matrix $D_{N},$
we can use several matrices $D_{N}^{\left( 1\right) },\ldots ,D_{N}^{(m)}$,
provided that they converge in distribution to some variables $d_{1},\ldots
,d_{m}$ as $N\rightarrow \infty .$

Moreover, instead of non-random matrices $D^{(i)}_N,$ we can use random matrices
which are independent of Gaussian matrices $A_{N}.$

Finally, we can use independent Gaussian matrices $A_{N}^{\left( i\right) }$
as some of these $D_{N}^{\left( i\right) }.$ In fact, the following theorem is true.

\begin{theo}
Let $A_{N}^{\left( 1\right) },\ldots ,A_{N}^{\left( p\right) }$ be
independent Hermitian Gaussian random matrices and let $D_{N}^{\left(
1\right) },\ldots ,D_{N}^{\left( r\right) }$ be random matrices which are
independent of $A$ matrices. Assume that the $r$-tuple $D_{N}^{\left(
1\right) },\ldots ,D_{N}^{\left( r\right) }$ converge in distribution to an $%
r$-tuple $\left( d_{1},\ldots ,d_{r}\right) .$ Then $A_{N}^{\left( 1\right)
},\ldots ,A_{N}^{\left( p\right) },D_{N}^{\left( 1\right) },\ldots
,D_{N}^{\left( r\right) }$ converge in distribution to $\left( s_{1},\ldots
,s_{p},d_{1},\ldots ,d_{p}\right) $ where $s_{i}$ are semicircle random
variables and the sets $\{s_1\}, \ldots, \{s_p\}$, and  $\{d_{1},\ldots
,d_{p}\}$ are free. 
\end{theo}

\section{Haar unitary matrices}

We are going to prove the following theorem.

\index{asymptotic freeness!Haar unitary conjugation}
\begin{theo}
\label{theorem_af_unitary}Let $\{A_{N} \in (M_N(\bC), \tr \otimes \E)\}$ and $\{B_{N} \in (M_N(\bC), \tr \otimes \E)\}$  be two sequence of $N$-by-$N$ Hermitian random matrices, which individually converge in distribution to non-commutative random variables $a$ and $b \in (\AA, \phi)$. Assume without loss of generality that $a$ and $b$ are free. 

Let $\{U_{N}\}$ be a sequence of $N$-by-$N$ independent random unitary matrices
that have the Haar distribution on the unitary group $\mathcal{U}\left(
N\right) .$ Then $A_{N}$ and $U_{N}B_{N}U_{N}^{\ast }$ jointly converge in
distribution to the pair of free variables $a$ and $b$. 
\end{theo}

The theorem essentially says that if we take two large random matrices $%
A_{N} $ and $B_{N}$ and if we conjugate one of them by a uniformly random
unitary transformation $U_{N},$ then the resulting pair of matrices $A_{N}$
and $U_{N}B_{N}U_{N}^{\ast }$ will be approximately free.

As a slogan, this can be put as follows:

\begin{center}
``Two large random matrices
in general position
are asymptotically free!''
\end{center}

The proof of Theorem \ref{theorem_af_unitary} is similar to the proof of
Theorem \ref{theorem_af_semicircle}. We are going to derive an asymptotic
formula for the following expected trace: 
\begin{equation}
\E \tr \left( A_{N}^{\left( 1\right)
}U_{N}B_{N}^{(1)}U_{N}^{\ast }A_{N}^{\left( 2\right) }\ldots
U_{N}B_{N}^{\left( n\right) }U_{N}^{\ast }\right).
\label{formula_af_unitary0}
\end{equation}%
Here $A_{N}^{\left( i\right) }$ and $B_{N}^{\left( i\right) }$ denote
polynomials of matrices $A_{N}$ and $B_{N},$ respectively.

The main difficulty in the derivation of a formula for (\ref%
{formula_af_unitary0}) is that we need an analogue of  Wick's formula for
the expectations of products of the elements of a Haar-distributed unitary
matrix.

\index{Weingarten function|textbf}
Such formulas are called Weingarten formulas (see \cite{weingarten78}).

Namely, note that the distribution of matrix $U_{N}$ is invariant if we
multiply it by $e^{i\alpha }I.$ This implies that the expection of $%
u_{i_{1}j_{1}}\ldots u_{i_{p}j_{p}}\overline{u}_{i_{1}^{\prime
}j_{1}^{\prime }}\ldots \overline{u}_{i_{q}^{\prime }j_{q}^{\prime }}$ is
zero if $p\neq q.$ If $p=q,$ then we have the following formula: 
\begin{equation}
\E \bigl(u_{i_{1}j_{1}}\ldots u_{i_{q}j_{q}}\overline{u}_{i_{1}^{\prime
}j_{1}^{\prime }}\ldots \overline{u}_{i_{q}^{\prime }j_{q}^{\prime
}}\bigr) =\sum_{\alpha ,\beta }\delta _{i_{1}i_{\beta \left( 1\right)
}^{\prime }}\ldots \delta _{j_{1}j_{\alpha \left( 1\right) }^{\prime
}}\ldots \mathrm{Wg}\left( N,\beta \alpha ^{-1}\right) .
\label{formula_Weingarten}
\end{equation}%
Here the sum is over permutations $\alpha $ and $\beta $ of the set $\left\{
1,2,\ldots ,q\right\} $. The coefficient $\mathrm{Wg}\left( N,\beta \alpha
^{-1}\right) $ is called the Weingarten function. We can define it by the
following equality: 
\begin{equation*}
\mathrm{Wg}\left( N,\alpha \right) :=\E\bigl( u_{11}\ldots u_{qq}\overline{u}%
_{1\alpha \left( 1\right) }\ldots \overline{u}_{q\alpha \left( q\right)
}\bigr) .
\end{equation*}

There is a beautiful explicit formula for $\mathrm{Wg}\left( N,\alpha
\right) $ due to B. Collins (see \cite{collins02}). Namely, let $N\geq q.$
Then, 
\begin{equation*}
\mathrm{Wg}\left( N,\alpha \right) =\frac{1}{\left( q!\right) ^{2}}%
\sum_{\lambda \vdash q}\frac{\chi ^{\lambda }\left( \mathrm{id}\right)
^{2}\chi ^{\lambda }\left( \alpha \right) }{s_{\lambda ,N}\left( \mathrm{id}%
\right) }.
\end{equation*}%
Here, the sum is over all partitions of $q,$ $\chi ^{\lambda }$ is the
character of the irreducible representation of the symmetric group $S_{q}$
that corresponds to the partition $\lambda ,$ and $s_{\lambda ,N}$ is the
character of the irreducible representation of the unitary qroup $\mathcal{U}%
\left( N\right) ,$ that corresponds to the partition $\lambda .$

It is clear from this formula that $\mathrm{Wg}\left( N,\alpha \right) $
depends only on the conjugacy class of the permutation $\alpha .$

\begin{exe}
(from \cite{collins02}). Let us use the notation for the partition class as
the second argument, so, for example, the identity permutation in $S_{q}$
corresponds to $1^{q}.$ Check the following formulas: 
\begin{eqnarray*}
\mathrm{Wg}\left( N,1\right) &=&\frac{1}{N}, \\
\mathrm{Wg}\left( N,1^{2}\right) &=&\frac{1}{N^{2}-1}, \\
\mathrm{Wg}\left( N,1^{3}\right) &=&\frac{N^{2}-2}{N\left( N^{2}-1\right)
\left( N^{2}-4\right) }, \\
\mathrm{Wg}\left( N,2\right) &=&\frac{-1}{N\left( N^{2}-1\right) }, \\
\mathrm{Wg}\left( N,21\right) &=&\frac{-1}{\left( N^{2}-1\right) \left(
N^{2}-4\right) }, \\
\mathrm{Wg}\left( N,3\right) &=&\frac{2}{N\left( N^{2}-1\right) \left(
N^{2}-4\right) }.
\end{eqnarray*}
\end{exe}

\index{Weingarten function!asymptotic}
For our purposes, we are more interested in the asymptotic behavior of the
Weingarten function. It turns out that it is given by the following formula: 
\begin{equation}
\mathrm{Wg}\left( N,\alpha \right) =\mu( \alpha) N^{\#\left(
\alpha \right) -2q}+O\left( N^{\#\left( \alpha \right) -2q-2}\right) ,
\label{formula_asymptotic_Weingarten}
\end{equation}%
where $\mu( \alpha) $ is not zero.

The explicit formula for $\mu (\alpha) $ is as follows. Let
the conjugacy class for $\alpha $ be described by partition $\lambda =(\lambda _{1},\ldots ,\lambda _{s}) .$ Then, 
\begin{equation}
\label{defi_mu}
\mu (\alpha) =\prod_{i=1}^{s}\left( -1\right) ^{\lambda
_{i}-1}C_{\lambda _{i}},
\end{equation}%
where $C_{n}$ denote Catalan numbers.

For example, $\mu ( id) =1,$ $\mu \bigl( ( 12)\bigr) =-1,$ $\mu \bigl( ( 123) \bigr) =2,$ $\mu \bigl(( 12) ( 34) \bigr) =1.$

\begin{proof}[Sketch of the proof of Theorem \ref{theorem_af_unitary}] By using
formulas (\ref{formula_Weingarten}) and (\ref{formula_asymptotic_Weingarten}%
), it is possible to derive the following asymptotic formula for the trace:%
\begin{align}
&\lim_{N\rightarrow \infty }\E \tr \left( A_{N}^{\left(
1\right) }U_{N}B_{N}^{(1)}U_{N}^{\ast }A_{N}^{\left( 2\right) }\ldots
U_{N}B_{N}^{\left( n\right) }U_{N}^{\ast }\right) 
\label{formula_af_unitary_trace} \\
&=\sum_{\substack{ \alpha ,\beta \in S_{n}  \\ \left| \alpha ^{-1}\beta
\right| +\left| \alpha \right| +\left| \beta ^{-1}\gamma \right| =n-1}}%
\phi_{\alpha }\left( a_{1},\ldots ,a_{n}\right) \phi_{\beta ^{-1}\gamma }\left(
b_{1},\ldots ,b_{n}\right) \mu \bigl( \alpha ^{-1}\beta \bigr)  \notag
\end{align}

In this formula, $A_{N}^{\left(
i\right) }$ and $B_{N}^{\left( j\right) }$ are polynomials of matrices $A_{N}$ and $B_{N},$ and variables $a_{i}$ and $b_{j}$  are the corresponding polynomials of variables $a$ and 
$b$. By assumption, for every permutation $\pi$, $\E \tr_\pi (A_N^{(1)}, \ldots, A_N^{(n)}$ converges to $\phi_\pi(a_1, \ldots, a_n)$. Similarly, $\E \tr_\pi (B_N^{(1)}, \ldots, B_N^{(n)}$ converges to $\phi_\pi(b_1, \ldots,b_n)$

The permutation $\gamma $ equals by definition $\left( 1,2,\ldots ,n\right)
, $ and $\left| \sigma \right| $ denotes the \emph{length of permutation}%
\index{length of permutation} $\sigma ,$ that is, the minimal number of
transpositions which is needed to represent $\sigma $ as their product. 

Note that the length function defines a distance on the set of all
permutations: $d\left( \sigma ,\rho \right) =d\left( \sigma ^{-1}\rho
\right) $ and it turns out that this distance is actually a metric. Hence we
can define a notion of a geodesic: a set of permuations is a \emph{geodesic} 
\index{geodesic}if (i) for every triple of permutations from this set, we
can order them in such a way $\left( \rho ,\sigma ,\tau \right) $ that $%
\left| \rho ^{-1}\sigma \right| +\left| \sigma ^{-1}\tau \right| =\left|
\rho ^{-1}\tau \right| ,$ and (ii) it is not possible to add another
permutation without violating property (i).

\begin{figure}[tbph]
\begin{center}
\includegraphics[width=6cm]{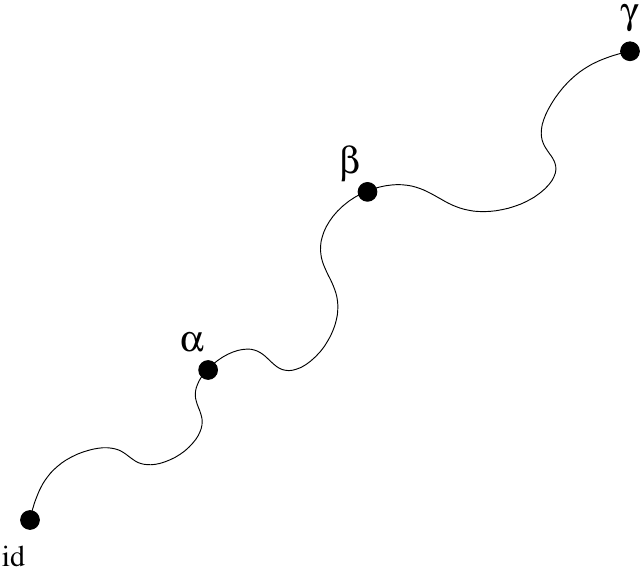}
\end{center}
\caption{A geodesic from $id$ to $\protect\gamma$}
\label{figure_geodesic}
\end{figure}

We are interested in the geodesics from the identity element $id$ to the
cycle $\gamma .$ It is possible to introduce the total ordering on \ each of
these geodesic. The permutation $\alpha $ precedes $\beta $ if $\alpha $ is
between $id$ and $\beta .$ In this ordering $id$ and $\gamma $ are largest
and the smallest element respectively.

In this terms, the condition 
\begin{equation*}
\left| \alpha ^{-1}\beta \right| +\left| \alpha \right| +\left| \beta
^{-1}\gamma \right| =n-1
\end{equation*}%
means that $\alpha $ and $\beta $ lie on a geodesic between $id$ and $\gamma 
$ and that $\alpha $ precedes $\beta $. This situation is shown in Figure %
\ref{figure_geodesic}.

\begin{lemma}
\label{lemma_af_fixed_element}Suppose that 
\begin{equation*}
\left| \alpha ^{-1}\beta \right| +\left| \alpha \right| +\left| \beta
^{-1}\gamma \right| =n-1.
\end{equation*}%
Then either $\alpha $ or $\beta ^{-1}\gamma $ has a fixed point.
\end{lemma}

\begin{proof}[Proof of Lemma \ref{lemma_af_fixed_element}] The assumption implies that either $\left| \alpha \right| \leq \left(
n-1\right) /2,$ or $\left| \beta ^{-1}\gamma \right| \leq \left( n-1\right)
/2.$ This means that one of these permutations can be represented as a
product of less than or equal to $\left( n-1\right) /2$ transposition. Since
a transposition moves only $2$ elements, therefore the product of no more
than $\left( n-1\right) /2$ transpositions moves no more than $n-1$
elements. Hence, at least one element remains unmoved.
\end{proof}

Suppose now that $\lim_{N\rightarrow \infty }\mathrm{tr}\left( A_{N}^{\left(
i\right) }\right) =0$ and $\lim_{n\rightarrow \infty }\mathrm{tr}\left(
B_{N}^{\left( i\right) }\right) =0$ for all $i.$ Then, Lemma \ref%
{lemma_af_fixed_element} and formula (\ref{formula_af_unitary_trace}) imply
that 
\begin{align*}
\lim_{N\to \infty }\E \tr\bigl( A_{N}^{(
1) }U_{N}B_{N}^{(1)}U_{N}^{\ast }A_{N}^{( 2) }\ldots
U_{N}B_{N}^{( n) }U_{N}^{\ast }\bigr)  =0 =\phi ( a_{1}%
b_{1}a_{2}\ldots b_{n}).
\end{align*}%

This means that the fundamental equation of freeness holds approximately for syb-algebra algebras generated by $A_N$ and $B_N$, respectively and concludes the proof of the theorem. 
\end{proof}
%
%
This completes the proof of the theorem.

\begin{coro}
Suppose $\{a_{1},\ldots ,a_{n}\}$ and $\{b_{1},\ldots ,b_{n}\}$ are free families of elements of
 a probability space with
a faithful tracial expectation $\phi$, then 
\begin{align}
&\phi\left( a_{1}b_{1}a_{2}\ldots b_{n}\right)
\label{formula_af_unitary_trace1} \\
&=\sum_{\substack{ \alpha ,\beta \in S_{n}  \\ \left| \alpha ^{-1}\beta
\right| +\left| \alpha \right| +\left| \beta ^{-1}\gamma \right| =n-1}}%
\phi_{\alpha }\left( a_{1},\ldots ,a_{n}\right) \phi_{\beta ^{-1}\gamma }\left(
b_{1},\ldots ,b_{n}\right) \mu \left( \alpha ^{-1}\beta \right)  \notag
\end{align}
\end{coro}

\begin{exa} 
Consider $\phi( abab) ,$ where $a$ and $b$ are free. Set $n=2$ in formula \eqref{formula_af_unitary_trace1}, $a_{1}=a_{2}=a,$ $b_{1}=b_{2}=b$.  Then, $\gamma = (1, 2)$ and we have three possibilities: (i) $\alpha =\beta =id,$ (ii) $\alpha =\beta =\gamma =(1,2) ,$ and (iii) $\alpha =id$, $\beta =\left( 12\right) .$ Note that the case when $%
\alpha =\left( 12\right) $ and $\beta =id$ is impossible, because $\alpha $
must precede $\beta $ on the geodesic from $id$ to $\gamma .$ To these
cases, we have the correspondent summands: (i) $\phi( a) ^2 \phi (b^2)$, (ii) $\phi(a^2)\phi(b)^2$, and (iii) $-\phi(
a) ^2 \phi(b)^2$. Hence,
\begin{equation*}
\phi(abab) =\phi( a) ^{2}\phi(b^{2}) +\phi(a^{2})\phi(b) ^{2}-\phi( a) ^2 \phi(b)^2.
\end{equation*}%
This formula coincides with formula \eqref{exa_prod_rule2} on page \pageref{exa_prod_rule2}.
\end{exa}

Formula \eqref{formula_af_unitary_trace1} is essentially the same as \eqref{eq:phi-alternating-product2} on page \pageref{eq:phi-alternating-product}. Indeed, there is a remarkable isomorphism between the lattice of non-crossing partitions $NC(n)$ and the lattice of permutations that belong to a geodesic $[e, \gamma]$, with $\alpha \leq \beta$ if $\beta$ is between $\alpha$ and $\gamma$. Writing the same symbol for both a permutation and the corresponding non-crossing partition, we can re-write \eqref{formula_af_unitary_trace1} as 
\[
\phi( a_{1}b_{1}a_{2}\ldots b_{n}) = \sum_{\alpha \in NC(n)} \phi_{\alpha}(a_1, \ldots, a_n) \sum_{\alpha \leq \beta} \phi_{K(\beta)}(b_1, \ldots, b_n) \mu(\alpha, \beta),   
\]
after checking that for $\alpha \leq \beta$, $\mu(\alpha^{-1} \beta)$, with $\mu$ defined as in 
\eqref{defi_mu}, equals to $\mu(\alpha, \beta)$ where $\mu$ is the M\"obius function of the lattice.

\section*{Notes} Asymptotic free independence was discovered in \cite{voiculescu91}. See also Chapter 4 in \cite{voiculescu_dykema_nica92}. The theory of Weingerten functions originated in  \cite{collins02}, \cite{collins2003}, where it is also applied to show results on asymptotic freeness of random matrices. For further developments, including results for orthogonal and symplectic groups,  see \cite{collins_sniady2006} and \cite{cmn2022}. Another account of asymptotic freeness for families of random matrices can be found in Chapter 4 of \cite{hiai_petz00}.

\chapter{S-transform}

\label{chapter_Stransform}

We have seen that if \(X\) and \(Y\) are free, the moments of \(X + Y\) can be calculated by using the additivity property of free cumulants and the relation between free cumulants and moments. Analytically, this corresponds to the fact that the \(R\)-transform is additive, \(R_{X + Y}(z) = R_X(z) + R_Y(z)\). Since the \(R\)-transform is closely related to the functional inverse of the Cauchy transform, this also gives us access to the Cauchy transform of \(X + Y\).

For the product \(XY\), the analogous analytic object is the \(S\)-transform, which has the property \(S_{XY}(z) = S_X(z)\,S_Y(z)\). The proof of this property, however, is more involved.

\index{S-transform|textbf}
\begin{defi}
Suppose \(X\) is a bounded non-commutative random variable such that \(\phi(X) \neq 0\). We introduce a moment-generating function 
\[
\psi_X(z) \;:=\; M_X(z) - 1 = \sum_{k = 1}^\infty \phi\bigl(X^k\bigr)\,z^k,
\]
where $M_X(z)$ was defined in \eqref{defi_M} and define \(S_X(z)\) by the formula:
\begin{equation}
\label{defi_Stransform}
S_X(z) \;=\; \frac{1 + z}{z}\,\psi_X^{(-1)}(z).
\end{equation}
\end{defi}


\index{S-transform!relation to R-transform}
Let us also give an equivalent definition of the \(S\)-transform. Let
\[
C_X(u) \;=\; u\,R_X(u) \;=\; \sum_{n = 1}^\infty \kappa_n(X)\,u^n,
\]
where \(\kappa_n(X)\) is the \(n\)-th free cumulant of \(X\). Then formula \eqref{equ_C_and_M} implies the relation
\[
\psi_X(z) \;=\; C_X\bigl[z\,(1 + \psi_X(z))\bigr],
\]
which can be written, by setting \(z = \psi_X^{(-1)}(u)\), as
\[
u \;=\; C_X\!\Bigl(\psi_X^{(-1)}(u) \times (1 + u)\Bigr).
\]
Hence, from \eqref{defi_Stransform}, 
\begin{equation}
\label{defi_Stransform_alt}
S_X(u) \;=\; \frac{1}{u}\;C_X^{(-1)}(u),
\end{equation}
where \(C_X^{(-1)}(u)\) is the functional inverse of \(C_X(u)\). 
(For more details, see pp.\ 335--336 in \cite{haagerup_larsen00}.)

We can also define the \(S\)-transform of probability measures with non-zero first moment (even if they are unbounded) by setting 
\[
\psi_\mu (z) \;=\; \int_0^\infty \frac{z\,t}{1 - z\,t}\, d\mu_x(t),
\]
provided that the integral is well defined. For example, if the probability measure $\mu_x$ has the support in $\R+$, one can define the integral for all $z$ that exclude the positive real axis.  
Then one can define $S_\mu(z)$ using
\[
S_\mu(z) \;=\; \frac{1 + z}{z}\;\psi_\mu^{(-1)}(z).
\]

Even if the first moment of \(\mu\) is zero and the function $\psi_\mu$ is not immediately invertible, one can still define the \(S\)-transform via a suitable generalization of the functional inverse of \(\psi_\mu(z)\); see \cite{rajrao_speicher2007}. Below are some examples.

\begin{exa}[Point mass at \(a\)]
Let \(\mu = \delta_a\). Then 
\[
\psi_\mu(z) \;=\; \frac{z\,a}{1 - z\,a}, 
\qquad
\chi_\mu(u) \;=\; \frac{u}{a\,(1 + u)},
\qquad
S_\mu(u) \;=\; \frac{1}{a}.
\]
\end{exa}

\index{bernoulli distribution!S-transform}
\begin{exa}[Bernoulli]
Consider the random variable \(X\) with the Bernoulli distribution \(\mu_X = p\,\delta_1 + (1-p)\,\delta_0\), where \(p > 0\). From (\ref{Bernoulli_Rtransform}) we know
\[
R_X(z) 
\;=\; 
\frac{-1 + z \;-\; \sqrt{\,(1 - z)^2 \;+\; 4\,p\,z}}{2\,z}.
\]
Then, by applying (\ref{defi_Stransform}) we obtain
\[
S_X(u) 
\;=\; 
\frac{1 + u}{p + u}.
\]
\end{exa}

\index{free Poisson r.v.!S-transform}
\begin{exa}[Free Poisson]
For a free Poisson random variable \(X\) with rate \(\lambda\), we have
\[
C_X(u) \;=\; \frac{\lambda\,u}{1 - u}.
\]
Hence,
\[
S_X(u) 
\;=\; 
\frac{1}{\lambda + u}.
\]
\end{exa}

\medskip

The main property of the \(S\)-transform is that it linearizes the product of free non-commutative random variables.

\begin{theo}[\cite{voiculescu87}]
Let \((\AA, \phi)\) be a non-commutative probability space, and let \(X,Y \in \AA\) with \(\phi(X)\neq 0\), \(\phi(Y)\neq 0\). If \(X\) is free from \(Y\), then 
\[
S_{XY}(z) 
\;=\; 
S_X(z)\,S_Y(z).
\]
\end{theo}

\cite{rajrao_speicher2007} showed that this theorem can be extended even to the case in which one of the variables \(X\) or \(Y\) has zero mean. The following proof follows \cite{rajrao_speicher2007}.

\begin{proof}
Define \(M_X(z) := 1 + \psi_X(z)\). We also define the power series
\[
M_1(z) 
\;:=\; 
\sum_{n = 0}^\infty \phi\bigl(Y\,(X\,Y)^n\bigr)\,z^n,
\qquad
M_2(z) 
\;:=\; 
\sum_{n = 0}^\infty \phi\bigl(X\,(Y\,X)^n\bigr)\,z^n.
\]
We want to relate \(M_{XY}\), \(M_1\), and \(M_2\). Our tool is  Theorem~\ref{theorem_expectation_products}.

\begin{figure}[htbp]
\centering
\includegraphics[width=0.9\textwidth]{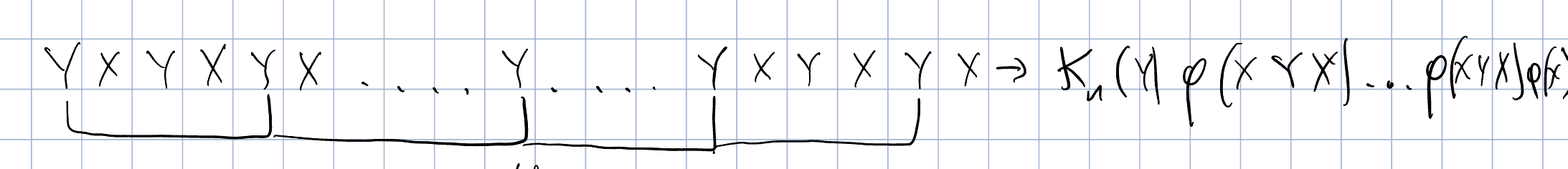}
\caption{One term in the decomposition for \(\phi\bigl((Y\,X)^n\bigr)\).}
\label{fig:S_transform1}
\end{figure}

Figure~\ref{fig:S_transform1} illustrates how Theorem~\ref{theorem_expectation_products} allows us to write \(\phi\bigl((Y\,X)^n\bigr)\) as a sum over non-crossing partitions \(\pi\) that connect the \(Y\)s. By choosing the part \(p\) of \(\pi\) that contains the first \(Y\), the complementary partition \(K(\pi)\) splits into terms located between the elements of \(p\). After summing over all such parts \(p\), \(\phi\bigl((Y\,X)^n\bigr)\) can be expressed as 
\(\sum_p \kappa_{|p|}(Y)\,\phi(a_1)\cdots \phi(a_{|p|}),\)
where \(a_1,\dots,a_{|p|}\) are the products \(X\,Y\,X\,\dots\,X\) separated by the elements of \(p\).

Applying this to the terms in \(M_{X\,Y}(z) = M_{Y\,X}(z)\) and regrouping yields
\begin{equation}
\label{equ_Srelation_1}
M_{X\,Y}(z) 
\;=\; 
M_{Y\,X}(z) 
\;=\; 
C_Y\bigl[z\,M_2(z)\bigr] + 1.
\end{equation}
Applying the same reasoning to \(\phi\bigl(Y\,(X\,Y)^n\bigr)\) gives
\begin{equation}
\label{equ_Srelation_2}
M_1(z) 
\;=\; 
\frac{C_Y\bigl[z\,M_2(z)\bigr]}{z\,M_2(z)}\; M_{X\,Y}(z),
\end{equation}
and by symmetry we obtain
\[
M_2(z)
\;=\;
\frac{C_X\bigl[z\,M_1(z)\bigr]}{z\,M_1(z)}\;M_{Y\,X}(z).
\]

Now let \(\chi(z) := \psi_{X\,Y}^{(-1)}(z)\). Such an inverse exists in the sense of power series because the first term in \(\psi_{X\,Y}(z)\) is \(\phi(X\,Y)\,z\), and \(\phi(X\,Y)\neq 0\). Substituting \(\chi(z)\) in \eqref{equ_Srelation_1} gives
\[
z 
\;=\; 
C_Y\Bigl[\chi(z)\,M_2\bigl(\chi(z)\bigr)\Bigr].
\]
By \eqref{defi_Stransform_alt}, 
\begin{equation}
\label{equ_Srelation_3}
z\,S_Y(z) 
\;=\; 
\chi(z)\,M_2\bigl(\chi(z)\bigr).
\end{equation}
Symmetry yields
\begin{equation}
\label{equ_Srelation_4}
z\,S_X(z) 
\;=\; 
\chi(z)\,M_1\bigl(\chi(z)\bigr).
\end{equation}
From \eqref{equ_Srelation_2} we further get
\begin{align}
\chi(z)\,M_1\bigl(\chi(z)\bigr)\,M_2\bigl(\chi(z)\bigr)
\;=\;
C_Y\Bigl[\chi(z)\,M_2\bigl(\chi(z)\bigr)\Bigr]\,
M_{X\,Y}\bigl(\chi(z)\bigr)
\;=\;
z\,(1 + z),
\label{equ_Srelation_5}
\end{align}
and multiplying \eqref{equ_Srelation_3} and \eqref{equ_Srelation_4}, then using \eqref{equ_Srelation_5}, yields
\[
S_X(z)\,S_Y(z) 
\;=\; 
\chi(z)\,\frac{1 + z}{z}
\;=\; 
S_{X\,Y}(z).
\]
\end{proof}

\index{S-transform!scaling}
\begin{propo}[Scaling property]
\label{propo_S_scaling}
\[
S_{aX}(z) 
\;=\; 
\frac{1}{a}\,S_X(z).
\]
\end{propo}

\begin{proof}
First, we compute \(C_{aX}^{(-1)}(z)\). Suppose \(C_{aX}(u) = z\). By the scaling property for \(C_X\), 
\[
z 
\;=\; 
C_X(a\,u) 
\quad\Longleftrightarrow\quad 
u 
\;=\; 
\frac{1}{a}\,C_X^{(-1)}(z).
\]
Hence,
\[
C_{aX}^{(-1)}(z) 
\;=\; 
\frac{z}{a}\,S_X(z),
\]
and therefore
\[
S_{aX}(z)
\;=\;
\frac{1}{z}\,C_{aX}^{(-1)}(z)
\;=\;
\frac{1}{a}\,S_X(z).
\]
\end{proof}

\section*{Notes}

The \(S\)-transform was introduced in \cite{voiculescu87}, where its multiplicativity was first proved using analytic methods and certain non-commutative random variables from Cuntz algebras. Another analytic proof was given in \cite{haagerup97}. A combinatorial proof appears in \cite{nica_speicher97b} (see also the textbook \cite{nica_speicher06}).
%
%

%


\chapter{Subordination}

In the previous chapter, we used \(R\)-transforms to calculate free additive convolutions. The analytical and combinatorial definitions of the \(R\)-transform are based on a functional inverse and on free cumulants, respectively. Both can be difficult to calculate, and the free cumulant method is restricted to measures for which we can compute all moments. One alternative approach is to use the \emph{subordination functions} instead of the \(R\)-transform. This approach has several advantages, both computational and theoretical, which we will explain later. 

\section{... for free additive convolutions}

The subordination method is based on an interesting identity. Let \(\mu\) and \(\nu\) be two measures on \(\mathbb{R}\) with Cauchy transforms \(G_\mu(z)\) and \(G_\nu(z)\), and let the corresponding \(R\)-transforms be denoted \(R_\mu(z)\) and \(R_\nu(z)\). (For now, we assume that \(\mu\) and \(\nu\) are compactly supported.) Define
\begin{equation}
\label{defi_omega_1}
\omega_1(z) = z - R_\nu\bigl(G_{\mu \boxplus \nu}(z)\bigr),
\end{equation}
and let us calculate \(G_\mu\bigl(\omega_1(z)\bigr)\). Looking ahead, it might be surprising that we will obtain \(G_{\mu \boxplus \nu}(z)\).

Indeed, observe that by Theorems~\ref{theo_R_transform_identity} and~\ref{theo_additivity_R} we have
\[
\begin{aligned}
R_\nu\bigl(G_{\mu \boxplus \nu}(z)\bigr) 
&= R_{\mu \boxplus \nu}\bigl(G_{\mu \boxplus \nu}(z)\bigr) \;-\; R_{\mu}\bigl(G_{\mu \boxplus \nu}(z)\bigr)
\\[6pt]
&= z \;-\; \frac{1}{G_{\mu \boxplus \nu}(z)} \;-\; R_{\mu}\bigl(G_{\mu \boxplus \nu}(z)\bigr).
\end{aligned}
\]
Therefore,
\[
\omega_1(z) 
\;=\; \frac{1}{G_{\mu \boxplus \nu}(z)} \;+\; R_{\mu}\bigl(G_{\mu \boxplus \nu}(z)\bigr),
\]
and, applying again Theorem~\ref{theo_R_transform_identity}, we get
\begin{align}
G_\mu\bigl(\omega_1(z)\bigr) 
&= G_\mu\Bigl(\frac{1}{G_{\mu \boxplus \nu}(z)} + R_{\mu}\bigl(G_{\mu \boxplus \nu}(z)\bigr)\Bigr) \notag
\\[3pt]
&= G_{\mu \boxplus \nu}(z). \label{equ_subordination_1}
\end{align}

Similarly, if we define 
\begin{equation}
\label{defi_omega_2}
\omega_2(z) \;=\; z \;-\; R_\mu\bigl(G_{\mu \boxplus \nu}(z)\bigr),
\end{equation}
then an analogous calculation yields
\begin{equation}
G_\nu \bigl(\omega_2(z)\bigr)  
\;=\; G_{\mu \boxplus \nu}(z).
\label{equ_subordination_2}
\end{equation}

\index{subordination|textbf}
\index{subordination!for additive convolution}
In the theory of functions of a complex variable, the phenomenon exhibited in \eqref{equ_subordination_1} and \eqref{equ_subordination_2} is called \emph{subordination}. The function \(G_{\mu \boxplus \nu}(z)\) is subordinated to \(G_{\mu}(z)\) and \(G_{\nu}(z)\), and \(\omega_1(z)\), \(\omega_2(z)\) are the corresponding \emph{subordination functions}.

Next, we can compute:
\begin{align}
\omega_1(z) + \omega_2(z) 
&= 2z \;-\; R_\nu\bigl(G_{\mu \boxplus \nu}(z)\bigr) \;-\; R_\mu\bigl(G_{\mu \boxplus \nu}(z)\bigr) 
\notag
\\[3pt]
&= 2z \;-\; R_{\mu \boxplus \nu}\bigl(G_{\mu \boxplus \nu}(z)\bigr)
\notag
\\[3pt]
&= z \;+\; \frac{1}{G_{\mu \boxplus \nu}(z)}.
\label{equ_subordination_3}
\end{align}

Hence, using \eqref{equ_subordination_1} and \eqref{equ_subordination_2}, we find the system of equations:
\begin{equation}
\label{system_subordination}
\begin{cases}
\omega_1(z) \;=\; z - \omega_2(z) \;+\;\dfrac{1}{G_\nu\bigl(\omega_2(z)\bigr)}, \\[6pt]
\omega_2(z) \;=\; z - \omega_1(z) \;+\;\dfrac{1}{G_\mu\bigl(\omega_1(z)\bigr)}.
\end{cases}
\end{equation}

This is a system of two equations in the two unknown functions \(\omega_1(z)\) and \(\omega_2(z)\). Note that for large \(|z|\), \(G_{\mu \boxplus \nu}(z) \sim z^{-1}\), and for small \(|z|\), \(R_\nu(z)\) is bounded. From definition \eqref{defi_omega_1}, we see that \(\omega_1(z) \sim z\) for large \(|z|\). Similarly, \(\omega_2(z) \sim z\) for large \(|z|\). These boundary conditions usually determine the solutions of the system \eqref{system_subordination} uniquely.

We now state this result as a theorem. For convenience, define the reciprocal Cauchy transform as 
\[
F_\mu(z) \;=\; \frac{1}{G_\mu(z)}.
\]
A useful property of \(F_\mu(z)\) is that \(\Im\bigl(F_\mu(z) - z\bigr) \geq 0\) for all \(z \in \bC^+\).

\index{subordination!iterative computation}
\begin{theo}
\label{theo_additive_subordination}
Given Borel probability measures \(\mu, \nu\) on \(\R\), there exist unique functions \(\omega_1, \omega_2: \bC^{+} \to \bC^{+}\) such that for all \(z \in \bC^{+}\),
\begin{enumerate}
\item \(\Im\omega_j(z) \ge \Im z\) and \(\displaystyle \lim_{y \uparrow \infty} \frac{\omega_j(i y)}{i y} = 1\) for \(j \in \{1, 2\}\);
\item \(F_{\mu \boxplus \nu}(z) = F_\mu\bigl(\omega_1(z)\bigr) = F_\nu\bigl(\omega_2(z)\bigr)\);
\item \(\omega_1(z) + \omega_2(z) = z + F_{\mu \boxplus \nu}(z)\).
\end{enumerate}
\end{theo}

\noindent
For the details of the proof, we refer to \cite{belinschi_bercovici07}.

\medskip

What is even more remarkable is that we can solve this system iteratively, starting with \(\omega^{(0)}_1(z) = \omega^{(0)}_2(z) = z\) and then iterating until convergence. This iterative scheme yields a major computational advantage for subordination: the subordination functions themselves are straightforward to compute numerically for each \(z\).

Once \(\omega_1(z)\) and \(\omega_2(z)\) are known, we can calculate \(G_{\mu \boxplus \nu}(z)\) using any one of \eqref{equ_subordination_1}, \eqref{equ_subordination_2} or \eqref{equ_subordination_3}, and then determine the density of \(\mu \boxplus \nu\) by using the Stieltjes inversion formula \eqref{formula_Stieltjes_inversion}.

If we define 
\begin{equation}
\label{defi_H}
H_\mu(z) \;=\; F_\mu(z) \;-\; z,
\end{equation}
and define \(H_\nu(z)\) similarly, then we can rewrite system \eqref{system_subordination} in a simpler form:
\begin{equation}
\label{system_subordination2}
\begin{cases}
\omega_1(z) \;=\; z \;+\; H_\nu\bigl(\omega_2(z)\bigr),  \\[6pt]
\omega_2(z) \;=\; z \;+\; H_\mu\bigl(\omega_1(z)\bigr),
\end{cases}
\end{equation}
and this can be written as a single equation for one of the subordination functions, for example,
\begin{equation}
\label{equation_subordination}
\omega_1(z) \;=\; z \;+\; H_\nu\Bigl(z \;+\;  H_\mu\bigl(\omega_1(z)\bigr)\Bigr).
\end{equation}

An important and very useful fact is that the subordination functions are analytic in all of \(\bC^{+}\) and provide an injective mapping \(\bC^+ \to \bC^+\). Moreover, the method of subordination functions applies to free additive convolutions \(\mu \boxplus \nu\) of unbounded probability measures \(\mu\) and \(\nu\). (The existence of these convolutions was originally proved in \cite{bercovici_voiculescu93} using the \(R\)-transform.) The significance of subordination functions was especially emphasized in the work of Belinschi and Bercovici (see, for example, \cite{belinschi_bercovici07}).

\begin{exa}[Sum of two free symmetries]
\label{exa_Bernoulli_sum}
Let \(\mu = \tfrac{1}{2}(\delta_{-1} + \delta_{1})\). Our goal is to compute \(\mu \boxplus \mu\). 
\end{exa}

We have
\[
\begin{aligned}
G_\mu(z) \;&=\; \tfrac{1}{2}\bigl(\tfrac{1}{z + 1} \;+\; \tfrac{1}{z - 1}\bigr)
\;=\; \frac{z}{z^2 - 1}, 
\\[6pt]
H_\mu(z) \;&=\; \frac{1}{G_\mu(z)} - z \;=\; - \frac{1}{z}. 
\end{aligned}
\]
Therefore, by \eqref{system_subordination2}, the equation for \(\omega(z) = \omega_1(z) = \omega_2(z)\) is 
\[
\omega(z) \;=\; z \;-\; \frac{1}{\omega(z)},
\]
which solves to 
\begin{equation}
\label{omega_Bernoulli}
\omega(z) 
\;=\; \frac{z + \sqrt{z^2 - 4}}{2}.
\end{equation}
Then, by \eqref{equ_subordination_3}, 
\[
G_{\mu \boxplus \mu}(z) 
\;=\; \frac{1}{\,2\,\omega(z) - z}
\;=\; \frac{1}{\sqrt{z^2 - 4}}, 
\]
and the corresponding density is that of the arcsine law:
\[
f_{\mu \boxplus \mu}(x) 
\;=\; \frac{1}{\pi \,\sqrt{\,4 - x^2\,}}, 
\quad x \in (-2, 2),
\]
in agreement with the results of Examples~\ref{exe_arcsine} and~\ref{exe_bernoulli}.

By Example~\ref{exe_bernoulli}, the \(R\)-transform of the Bernoulli measure \(\tfrac{1}{2}(\delta_{-1} + \delta_{1})\) is 
\begin{equation}
\label{R_Bernoulli}
R(z) \;=\; \frac{-1 \;+\; \sqrt{\,1 \,+\, 4z^2}}{2z}.
\end{equation}
A significant difference between \eqref{omega_Bernoulli} and \eqref{R_Bernoulli} is that \(\omega(z)\) can be defined in the upper half-plane \(\bC^+\), since its singularities are \(\pm 2\in \mathbb{R}\). By contrast, \(R(z)\) has singularities at \(\pm i/2\) and thus cannot be extended to all of \(\bC^+\). This situation holds in general and illustrates one of the key advantages of the subordination-based approach.

\begin{exa}[Free sum of two Cauchy distributions]
For the Cauchy distribution on \(\R\), we have 
\[
G_\mu(z) \;=\; \frac{1}{z + i},
\quad
H_\mu(z) \;=\; i.
\]
Then,
\[
\omega(z) \;=\; z \;+\; i,
\]
and we calculate 
\[
G_{\mu \boxplus \mu}(z) 
\;=\; G_\mu(z + i) 
\;=\; \frac{1}{z + 2i},
\]
which leads to the density 
\[
f_{\mu \boxplus \mu}(x) 
\;=\; \frac{1}{2\pi\bigl[1 + (x/2)^2\bigr]},
\]
i.e., the Cauchy density scaled by 2. 
\end{exa}


\section{... and conditional expectations}

\index{conditional expectation|textbf}
Let $(\AA, \phi)$ be a non-commutative probability space, where $\AA$ is a von Neumann algebra with a normal, faithful, tracial state (expectation) $\phi$. Let $\BB \subset \AA$ be a unital subalgebra of $\AA$. In classical probability theory, there is a concept of \emph{conditional expectation} mapping random variables in $\AA$ to random variables in $\BB$. It turns out that conditional expectations are also well-defined for non-commutative probability spaces, via the usual axiomatic requirements.

We will use the notation $\E(X \mid \BB)$ for the conditional expectation of $X$ onto $\BB$. It is characterized by:
\begin{enumerate}
\item
$\phi\bigl[\E(X \mid \BB)\bigr] = \phi[X]$,
\item
$\E\bigl(B_1 X B_2 \mid \BB\bigr) \;=\; B_1\,\E(X \mid \BB)\,B_2 \quad \text{for all } B_1, B_2 \in \BB$.
\end{enumerate}
These two properties fully determine $\E(\cdot \mid \BB)$. In particular, if there exists an operator $B \in \BB$ such that $\phi(A\,H) = \phi(B\,H)$ for every $H \in \BB$, then we must have
\[
\E(A \mid \BB) \;=\; B.
\]

 \cite{biane98b} showed that the subordination phenomenon holds at the operator level, not only at the level of probability measures. Biane's approach is combinatorial and links subordination to the theory of free cumulants. Later, \cite{voiculescu2002} generalized Biane's results to more complex situations.

We will use the notation $\E_{\bC[X]}(Y) := \E(Y \mid \bC[X])$, where $\bC[X]$ denotes the unital algebra generated by $X$ over the complex numbers.


\subsection*{Conditional expectation for the resolvent of a sum}
In classical probability, one often wants to compute conditional expectations such as
\[
\E\bigl[h(X + Y) \,\big\vert\, X\bigr],
\]
where $X$ and $Y$ are independent and $h$ is a given function. These identities are relevant in, for example, studying Markov or martingale properties of processes with independent increments. One standard example is
\[
\E\bigl[e^{\,i\,t(X + Y)} \,\big\vert\, Y\bigr]
\;=\;
e^{\,i\,t\,Y}\,\E\bigl[e^{\,i\,t\,X}\bigr]
\;=\;
e^{\,i\,t\,Y}\,\phi_X(t),
\]
where $\phi_X(t)$ is the characteristic function of $X$. By expanding both sides in powers of $t$, one obtains additional identities.

For non-commutative random variables, the following analog holds (due to Biane).

\index{conditional expectation!resolvent of a sum}
\begin{theo}\label{thm:Biane-subordination}
Suppose $X$ and $Y$ are free self-adjoint operators in $\AA$. Then for all $z \in \bC^{+} = \{x + i\,y : y > 0\}$,
\begin{equation}
\label{equ_Biane_subordination}
\E_{\bC[X]}\Bigl[\bigl(z\,I - (X + Y)\bigr)^{-1}\Bigr]
\;=\;
\bigl(\omega(z)\,I - X\bigr)^{-1},
\end{equation}
where $\omega(z)$ is an injective analytic function $\bC^{+} \to \bC^{+}$ that satisfies $\omega(\overline{z}) = \overline{\omega(z)}$.
\end{theo}

Suppose $X$ and $Y$ have distributions $\mu$ and $\nu$, respectively. Then applying $\phi$ to both sides of \eqref{equ_Biane_subordination} gives
\[
G_{\mu \boxplus \nu}(z) \;=\; G_{\mu}\bigl(\omega(z)\bigr).
\]
Comparing with \eqref{equ_subordination_1}, we identify $\omega(z)$ with the subordination function $\omega_1(z)$. Hence, the claims about injectivity and $\omega(\overline{z}) = \overline{\omega(z)}$ follow from known properties of $\omega_1$.

We will prove the result for bounded $X$ and $Y$, although it can be extended to (possibly unbounded) self-adjoint operators affiliated with $\AA$.

First, observe that
\[
z\,I - X - Y
\;=\;
\bigl(I - Y\,(z\,I - X)^{-1}\bigr)\,\bigl(z\,I - X\bigr),
\]
so
\begin{equation}
\label{equ_Biane_identity}
\bigl(z\,I - (X + Y)\bigr)^{-1}
\;=\;
\bigl(z\,I - X\bigr)^{-1}\,\Bigl[I - Y\,\bigl(z\,I - X\bigr)^{-1}\Bigr]^{-1}.
\end{equation}
Denote $R \equiv R_X(z) := (z\,I - X)^{-1}$. Then
\[
\bigl(z\,I - (X + Y)\bigr)^{-1}
\;=\;
R\,\bigl[I - Y\,R\bigr]^{-1}
\;=\;
R \;+\; R\,Y\,R \;+\; R\,Y\,R\,Y\,R \;+\;\cdots.
\]
Applying $\E_{\bC[X]}$ yields
\begin{equation}\label{equ_sum_RYR}
\E_{\bC[X]}\Bigl[\bigl(z\,I - (X + Y)\bigr)^{-1}\Bigr]
\;=\;
R \;+\; R\,\E_{\bC[X]}(Y)\,R
\;+\;
R\,\E_{\bC[X]}(Y\,R\,Y)\,R
\;+\;\cdots,
\end{equation}
which motivates the study of
\[
\E_{\bC[X]}\!\bigl[(Y\,R)^{n - 1}\,Y\bigr]
\;=\;
\E_{\bC[X]}\!\bigl[Y\,R\,Y\,R\,\cdots\,Y\,R\,Y\bigr].
\]
This suggest analyzing expressions of the form
$\phi\bigl((Y\,R)^{n - 1}\,Y\,H\bigr)$ for an arbitrary $H \in \bC[X]$. This can be done using non-crossing partitions and the Kreweras complement. We will repeat here an instructive argument  from \cite{biane98b}.

By applying formula \eqref{eq:phi-alternating-product},
we find that 
\bal{
\phi( (Y R)^{n - 1} Y H) &=\sum_{\pi \in
NC( n) }k_{\pi }(Y, \ldots, Y ) \phi_{K( \pi) }( R,\ldots ,R, H)
\\
&= \sum_{\pi \in
NC( n) }k_{K(\pi) }( Y,\ldots , Y) \phi_{\pi}( R,\ldots , R, H),
}
where in the second line we used that the Kreweras complement is a bijection. In the last formula we think about $\pi$ as a non-crossing partitions of $\{2, 4, \ldots, 2n\}$ (the positions of $R$'s and $H$ in the sequence $ Y, R, Y, \ldots, R, Y, H$ and $K(\pi)$ is a non-crossing partition of $\{1, 3, \ldots, 2 n - 1\}$ (the positions of $Y$'s). 

We split this sum in several parts according to the block of $\pi$ to which $H$ belongs. Suppose this block $b$  has $s$ elements.  Let $0< l_1 <  \ldots < l_s = 2n$ be the even numbers that form this block. For example if we have $n = 4$ 
\bal{
 Y R Y R Y R Y H, 
}
and $\pi$ connects $H$ with the first and third $R$'s, then $s = 3$,  $l_1 = 2$, $l_2 = 6$, and $l_3 = 8$.  For convenience let us also define $l_0 = 0$
Then the partition $\pi$ induces $s$ non-crossing partitions $\pi_1, \ldots, \pi_s$ on the even numbers between $l_{i - 1}$ and $l_i$, some of which might be empty. And $K(\pi)$ is also split into $s$ partitions, which are $\tilde K(\pi_1), \ldots, \tilde K(\pi_s)$. Here $\tilde K(\pi_i)$ denotes the largest non-crossing partition $\pi'$ of the odd numbers between $l_{i - 1}$ and $l_{i}$ that has the property that $\pi_i \cup  \pi'$ is non-crossing. 
Then we write
\bal{
\phi( (Y R)^{n - 1} Y H)  = \sum_{s = 1}^{n} \Big[ \sum_{b = (l_1, \ldots, l_s)} \prod_{j = 1}^s 
\Big(\sum_{\pi \in NC(l_{j - 1}, l_j)}  k_{\tilde K(\pi)}(Y) \phi_\pi(R)\Big)\Big] \phi(R^{s - 1} H)
}
where $NC(l_{j - 1}, l_j)$ denote the set of non-crossing partitions of the set of even numbers
between $l_{j - 1}$ and  $l_j$ (where the empty partitions are allowed), and we write 
$k_{\tilde K(\pi)}(Y)$ and $\phi_\pi(R)$ as shortcuts for $k_{\tilde K(\pi)}(Y, \ldots, Y)$ and $\phi_\pi(R, \ldots, R)$, respectively. 

Since this equality holds for every $H \in \bC[Y]$, then we can conclude that 
\bal{
\E_{\bC[X]}( (Y R)^{n - 1} Y ) = \sum_{s = 1}^{n} \Big[ \sum_{b = (l_1, \ldots, l_s)} \prod_{j = 1}^s 
\Big(\sum_{\pi \in NC(l_{j - 1}, l_j)}  k_{\tilde K(\pi)}(Y) \phi_\pi(R)\Big)\Big] R^{s - 1}
}

Crucially, this implies that 
\begin{equation}
\label{equ_big_sum}
\E_{\bC[X]}( \sum_{n=1}^\infty   (Y R)^{n - 1} Y) = \sum_{s = 1}^\infty \Big(\sum_{l = 1}^\infty \sum_{\pi \in NC(0, 2l)}  k_{\tilde K(\pi)}(Y) \phi_\pi(R)\Big)^s R^{s - 1}
\end{equation}
This suggests defining the function
\bal{
\delta(z) = \sum_{l = 1}^\infty \sum_{\pi \in NC(0, 2l)}  k_{\tilde K(\pi)}(Y) \phi_\pi(R)
}
where $NC(0,2l)$ is the set of relevant non-crossing partitions, and the series converges for $\Im(z)$ (or $|z|$) sufficiently large (e.g. for $|z| \geq C(\|X\| + \|Y\|)$).

Then we can rewrite \eqref{equ_sum_RYR} using \eqref{equ_big_sum} as 
\begin{equation}
\label{equ_cond_expect_sum}
\E_{\bC[X]}\Big[ \big(zI - (X + Y)\big)^{-1}\Big] = R + R \sum_{s = 1}^\infty [\delta(z) R]^s = \frac{R}{I - \delta(z) R}.
\end{equation}
%
Finally, note that by \eqref{equ_Biane_identity},
\[
\Bigl[ z\,I - X - \delta(z)\,I\Bigr]^{-1} = R( I - \delta(z) R)^{-1},
\]
hence,
\[
\frac{R}{\,I - \delta(z)\,R}
\;=\;
\bigl[(z - \delta(z))\,I - X\bigr]^{-1},
\]
which shows \eqref{equ_Biane_subordination} holds with $\omega(z) = z - \delta(z)$. This completes the proof.

\subsection*{Conditional expectation for the resolvent of a product}

Similarly, we can study the conditional expectation for the product of two free operators. Here it is convenient to introduce a modified resolvent, defined by
\[
\Psi_X(z) \;=\; zX\,(I - zX)^{-1} \;=\; zX \;+\;(zX)^{2}\;+\;(zX)^{3}\;+\;\cdots.
\]

Its expected value is the function $\psi_X(x)$ that we defined in Chapter \ref{chapter_Stransform}, 

\begin{equation}
\label{defi_psi_0}
\psi_X(z) \;=\; \phi(\Psi_X(z)) \;=\;
\sum_{k = 1}^\infty \phi\bigl(X^k\bigr)\,z^k
\;=\;
\int_0^\infty \!\frac{z\,t}{1 - z\,t}\,\mathrm{d}\mu_X(t).
\end{equation}

\index{conditional expectation!resolvent of a product}
\begin{theo}
Suppose \(X\) and \(Y\) are positive free operators in \(\AA\). Then
\begin{equation}\label{equ_Biane_subordination2}
\E_{\bC[X]} \Bigl[\Psi_{\,X^{1/2}\,Y\,X^{1/2}}(z)\Bigr]
\;=\;
\Psi_X\bigl(\omega(z)\bigr),
\end{equation}
where \(\omega(z)\) is an analytic function from \(\bC \setminus [0,+\infty)\) to \(\bC \setminus [0,+\infty)\), with the property \(\arg z < \arg \omega(z) < \pi\) for \(z \in \bC^+\). It extends continuously to \(\overline{\bC^+}\cup\{\infty\}\) and maps this set to itself.
\end{theo}

Taking the expectation in \eqref{equ_Biane_subordination2} then yields
\begin{equation}\label{equ_psi_multiplicative}
\psi_{\,X^{1/2}\,Y\,X^{1/2}}(z)
\;=\;
\psi_X\bigl(\omega(z)\bigr).
\end{equation}

\begin{proof}[Sketch of the proof.]
The argument parallels the one for sums of free operators. Observe that
\begin{align*}
\E_{\bC[X]}\!\Bigl[\Psi_{X^{1/2}\,Y\,X^{1/2}}(z)\Bigr]
&\;=\;
zX\,\E_{\bC[X]}[Y]
\;+\;
z^2\,X\,\E_{\bC[X]}\!\bigl[Y\,X\,Y\bigr]
\\
&\;+\;
z^3\,X\,\E_{\bC[X]}\!\bigl[Y\,X\,Y\,X\,Y\bigr]
\;+\;
\cdots,
\end{align*}
and one can adapt the same combinatorial argument as before to compute \(\E_{\bC[X]}\bigl[(YX)^{n-1} Y\bigr]\).

In particular, define
\[
\delta_X(z)
\;=\;
\sum_{l=1}^\infty \;\sum_{\pi \,\in\, NC(0,2l)} k_{\tilde K(\pi)}(Y)\,\phi_\pi\bigl(zX\bigr).
\]
Then, analogously to the derivations in \eqref{equ_big_sum} and \eqref{equ_cond_expect_sum}, we obtain
\[
\E_{\bC[X]}\!\bigl[\Psi_{X^{1/2}\,Y\,X^{1/2}}(z)\bigr]
\;=\;
\sum_{s=1}^\infty \bigl(\delta_X(z)\,zX\bigr)^{s}
\;=\;
\frac{z\,\delta_X(z)\,X}{I - z\,\delta_X(z)\,X}
\;=\;
\Psi_X\!\bigl(z\,\delta_X(z)\bigr).
\]
Hence we may set \(\omega(z) = z\,\delta_X(z)\).

To see that \(\omega(z)\) is analytic in \(\bC^+\), set \(A = X^{1/2}\,Y\,X^{1/2}\). Since \(A\) is positive-definite, its spectrum lies in \(\R^+\). For \(z = r\,e^{i\theta}\) with \(0 < \theta < \pi\), the spectrum of \(\tfrac{z\,A}{I - z\,A}\) belongs to the image of \(\R^+\) under the transformation \(x \mapsto \tfrac{r\,e^{i\theta}\,x}{1 - r\,e^{i\theta}\,x},\) which is a circular arc passing through \(0\) and \(-1\). Crucially, for \(\theta \in (0,\pi)\), this arc remains in the upper half-plane.

Since for any normal operator \(T\in\AA\), the spectrum of \(\E(T|\BB)\) is contained in the closed convex hull of the spectrum of \(T\) (see Lemma~4.3 in \cite{lehner_szpojankowski2021}), applying this to \(\Psi_A(z) = \tfrac{z\,A}{1 - z\,A}\) shows that \(I + \E_{\bC[X]}\!\bigl[\Psi_A(z)\bigr] = I - \omega(z)\,X\) remains invertible for \(z \in \bC^+\), and its inverse maps \(\bC^+\) into \(\bC^-\). Hence \(\omega(z)\) is analytic on \(\bC^+\) and satisfies \(\omega(z)\in\bC^+\) whenever \(z\in\bC^+\). For more refined properties of \(\omega(z)\), we refer to \cite{belinschi_bercovici07}.
\end{proof}

%
\section{... for free multiplicative convolutions}

\index{subordination!for multiplicative convolution}
Next, let \(X\) and \(Y\) have probability measures \(\mu\) and \(\nu\), respectively. Then \(X^{1/2}\,Y\,X^{1/2}\) and \(Y^{1/2}\,X\,Y^{1/2}\) both have law \(\mu\boxtimes\nu\). Hence \eqref{equ_psi_multiplicative} becomes
\begin{equation}\label{equ_psi_mult_2}
\psi_{\mu\boxtimes\nu}(z)
\;=\;
\psi_{\mu}\!\bigl(\omega_1(z)\bigr)
\;=\;
\psi_{\nu}\!\bigl(\omega_2(z)\bigr).
\end{equation}

We want an additional equation in this system, mirroring \eqref{system_subordination} and \eqref{system_subordination2} for additive subordination. To do this, we invoke the \(S\)-transform. Recall that by definition \eqref{defi_Stransform}, the \(S\)-transform can be related to \(\psi_\mu\) via
\begin{equation}\label{defi_Stransform_for_sub}
\psi_\mu^{(-1)}(z)
\;=\;
\frac{z}{z+1}\;S_\mu(z).
\end{equation}
Applying this relation to \eqref{equ_psi_mult_2} yields
\[
\omega_1(z)
\;=\;
\psi_\mu^{(-1)}\!\bigl(\psi_{\mu \boxtimes \nu}(z)\bigr)
\;=\;
\frac{\psi_{\mu \boxtimes \nu}(z)}{\psi_{\mu \boxtimes \nu}(z)+1}\;
S_\mu\!\bigl(\psi_{\mu \boxtimes \nu}(z)\bigr),
\]
\[
\omega_2(z)
\;=\;
\psi_\nu^{(-1)}\!\bigl(\psi_{\mu \boxtimes \nu}(z)\bigr)
\;=\;
\frac{\psi_{\mu \boxtimes \nu}(z)}{\psi_{\mu \boxtimes \nu}(z)+1}\;
S_\nu\!\bigl(\psi_{\mu \boxtimes \nu}(z)\bigr).
\]
Taking the product and using the multiplicativity of \(S\)-transforms,
\begin{equation}\label{equ_omega_product}
\omega_1(z)\,\omega_2(z)
\;=\;
\frac{\psi_{\mu \boxtimes \nu}(z)}{\psi_{\mu \boxtimes \nu}(z) + 1}\;z.
\end{equation}

\index{eta-function|textbf}
\index{$\eta$-function|see{eta-function}}
This motivates introducing
\[
\eta_\mu(z)
\;=\;
\frac{\psi_\mu(z)}{\psi_\mu(z)+1},
\]
often called the \emph{Boolean cumulant generating function}, since its Taylor coefficients are the Boolean cumulants of \(\mu\). (Boolean cumulants are the cumulants which arise when the lattice of interval partitions is used instead of the lattice of non-crossing partitions. This lattice arises when a different concept of independence is used -- Boolean independence instead of free independence.)

Being just a nonlinear reparametrization of \(\psi_\mu\), $\eta_\mu(z)$ also satisfies relations akin to \eqref{equ_psi_mult_2}:
\begin{equation}\label{equ_psi_mult_3}
\eta_{\mu \boxtimes \nu}(z)
\;=\;
\eta_\mu\bigl(\omega_1(z)\bigr)
\;=\;
\eta_\nu\bigl(\omega_2(z)\bigr),
\end{equation}
which by \eqref{equ_omega_product} can be rewritten as
\begin{equation}\label{equ_psi_mult_4}
\eta_\mu\bigl(\omega_1(z)\bigr)
\;=\;
\eta_\nu\bigl(\omega_2(z)\bigr)
\;=\;
\frac{\omega_1(z)\,\omega_2(z)}{z}.
\end{equation}
Hence we arrive at the system
\begin{equation}\label{system_subordination_mult}
\begin{cases}
\omega_1(z) \;=\; z\,\dfrac{\eta_\nu\bigl(\omega_2(z)\bigr)}{\omega_2(z)},\\[6pt]
\omega_2(z) \;=\; z\,\dfrac{\eta_\mu\bigl(\omega_1(z)\bigr)}{\omega_1(z)},
\end{cases}
\end{equation}
which can be used to compute \(\omega_1(z)\) and \(\omega_2(z)\) iteratively.

\section*{Notes}
Subordination in free probability originated in \cite{voiculescu93}, then was extended to conditional expectations in \cite{biane98b} and \cite{voiculescu2000}. It was developed into a key tool for understanding free convolution in \cite{belinschi_bercovici05} and \cite{belinschi_bercovici07}.

A very important fact is that subordination also extends naturally to the matrix-valued setting. Together with the linearization trick, this gives access to properties of polynomials in random variables; see \cite{bstv2014} and \cite{bms2017}. For Boolean cumulants and Boolean subordination, see \cite{lehner_szpojankowski2021}.


%

\chapter{Operator-valued non-commutative random variables}
\section{Basic definitions and operator-valued subordination}

The theory of free probability admits a very important extension to situations where the expectation map has a more involved structure.  In a certain sense, this extension provides a non-commutative analogue of conditional independence in classical probability theory.

We will mostly be interested in a particularly simple example, where our operator-valued random variables can be viewed as matrices whose entries are themselves non-commutative random variables. Concretely, consider a non-commutative probability space $\CC$ with state $\phi$, and form the matrix algebra $\MM = M_n(\bC) \otimes \CC$, which we may identify with $M_n(\CC)$. We define the map $\E: \MM \to \bC_{n\times n}$ by
\[
\E\bigl((x_{ij})\bigr) \;=\; \bigl(\,\phi(x_{ij})\bigr),
\]
that is, we apply the state $\phi$ to each entry of the matrix. We then define a scalar-valued state $\tau$ on $\MM$ by
\[
\tau(X) \;=\; \tr \bigl(\E(X)\bigr),
\]
where $\tr$ is the normalized trace on $n\times n$ matrices, i.e.\ $\tr(A) = \frac{1}{n}\,\Tr(A)$, with $\Tr$ denoting the usual (non-normalized) matrix trace.

More generally, in operator-valued free probability, we assume that we have a unital inclusion of algebras $\BB \subset \MM$ and a linear map $\E: \MM \to \BB$ which is a conditional expectation onto the subalgebra $\BB$. We write $\E_\BB(X) = \E(X \mid \BB)$ and often shorten it to $\E$ if the target $\BB$ is clear from context. The space $\MM$ also has a state, denoted by $\tau$. For the above matrix example, we take $\BB = M_n(\bC)$, and $\E_\BB$ is the map
\[
(x_{ij}) \;\mapsto\; \bigl(\,\phi(x_{ij})\bigr).
\]
In this context, the triple $(\MM, \E, \BB)$ is called an \emph{operator-valued non-commutative probability space}.

\index{freeness with amalgamation|textbf}
\begin{defi}
Let $\AA_1, \AA_2 \subset \MM$ be subalgebras containing $\BB$. We say that $\AA_1$ and $\AA_2$ are \emph{free over $\BB$} (or \emph{free with amalgamation over $\BB$ with respect to $\E$}) if for all $n \ge 1$ and all $x_j \in \AA_{i_j}$ with $i_j \neq i_{j+1}$, whenever $\E[x_j] = 0$ for each $j$, we have
\[
\E[x_1 \, x_2 \,\cdots\, x_n] = 0.
\]
Two random variables $x, y \in \MM$ are called free over $\BB$ if the unital subalgebras $\BB\langle x\rangle$ and $\BB\langle y\rangle$ they generate are free over $\BB$.
\end{defi}

\index{moments!of operator-valued r.v.s.}
The \emph{distribution} of a random variable $x \in \MM$ over $\BB$, denoted $\mu_x$, is given by its collection of \emph{moments}. Here, an $n$th moment becomes a $\BB$-valued multilinear form
\[
m_n^x : \BB^{n-1} \;\longrightarrow\; \BB,
\]
defined by
\[
m_n^x(b_1, \ldots, b_{n-1}) \;=\; \E\bigl[x \,b_1\, x \,b_2 \,\cdots\, x \,b_{n-1}\, x\bigr].
\]
It is known that if $x$ and $y$ are free over $\BB$, then knowing their individual distributions is enough to compute the distribution of $x + y$. This new distribution is denoted by $\mu_x \boxplus \mu_y$ and the operation is called the \emph{free additive convolution} of $\mu_x$ and $\mu_y$. By construction, this convolution is commutative and associative.

Unlike the scalar-valued case, it is far from obvious how to encode the $\BB$-valued moments in a functional transform that linearizes the convolution. Surprisingly, not only is this feasible, but one also obtains a subordination property analogous to the scalar case.

\index{Cauchy transform!operator-valued}
\index{Cauchy transform!fully matricial}
First, define the \emph{(operator-valued) Cauchy transform} of $x \in \MM$:
\[
G_x(b) \;=\; \E\bigl[(b - x)^{-1}\bigr],
\]
whenever $b - x$ is invertible in $\MM$. A more refined version, the \emph{fully matricial Cauchy transform}, is defined for each $n \ge 1$ by
\[
G_x^{(n)}(b) \;=\; (\E \otimes \mathrm{Id}_n)\bigl[(\,b \;-\; x \otimes I_n\,)^{-1}\bigr],
\]
where $b \in M_n(\BB) = \BB \otimes M_n(\bC)$. The fully matricial transforms encode all moments of $x$, so in principle one recovers the distribution of $x$ from the family $\{G_x^{(n)}\}$.

A suitable domain for $G_x$ is the \emph{operator upper half-plane} 
\[
\bH^+(\BB) \;=\; \bigl\{\, b \in \BB : \Im(b) > 0 \bigr\},
\]
where $\Im(b) = \frac{1}{2i}(b - b^*)$. On this domain, for self-adjoint $x$ the Cauchy transform $G_x(b)$ is well-defined and lands in $\bH^-(\BB) := -\bH^+(\BB)$. For brevity, we mainly focus on the $n=1$ case, $G_x^{(1)} = G_x$.

We also define the \emph{reciprocal Cauchy transform} and the $h$-transform:
\[
F_x(b) \;=\; G_x(b)^{-1},
\quad\quad
h_x(b) \;=\; F_x(b)\;-\; b.
\]
It can be shown that
\[
h_x\bigl(\bH^+(\BB)\bigr)\;\subseteq\;\overline{\bH^+(\BB)}.
\]

While $G_x$ and $F_x$ are often easier to handle, a full analogue of Voiculescu’s $R$-transform in the operator-valued setting is more complicated, making a direct calculation of $G_{x+y}$ for free $x$ and $y$ challenging. However, the remarkable \emph{subordination} phenomenon still applies and simplifies such calculations considerably.

\index{subordination!operator-valued}
Operator-valued subordination was first developed in \cite{voiculescu2000}, with a convenient version due to \cite{bms2017}:

\begin{theo}[Subordination for matrix free additive convolutions]
\label{theo_matrix_subordination}
Let $(\MM, \E, \BB)$ be a $C^*$-operator-valued non-commutative probability space, and let $x, y \in \MM$ be two self-adjoint random variables that are free over $\BB$. Then there is a unique pair of Fr\'echet-analytic maps
\[
\omega_1, \omega_2 : \bH^+(\BB) \;\longrightarrow\; \bH^+(\BB)
\]
such that for every $b \in \bH^+(\BB)$:
\begin{enumerate}
\item[(1)] $\Im\bigl(\omega_j(b)\bigr) \;\ge\; \Im(b)$ for $j \in \{1,2\}$;
\item[(2)] $F_x\bigl(\omega_1(b)\bigr) + b \;=\; \omega_1(b) + \omega_2(b) \;=\; F_y\bigl(\omega_2(b)\bigr) + b$;
\item[(3)] $G_x\bigl(\omega_1(b)\bigr) \;=\; G_y\bigl(\omega_2(b)\bigr) \;=\; G_{x+y}(b)$.
\end{enumerate}
\end{theo}

This result directly generalizes Theorem~\ref{theo_additive_subordination} from the scalar case. For a proof, see \cite{bms2017}. Its importance lies in showing that once $G_x$ and $G_y$ are known, we can compute $G_{x+y}$, and hence deduce all moment (and spectral) information for $x+y$.

Returning to our main example with $\BB = M_n(\bC)$, we note that the subordination maps $\omega_1$ and $\omega_2$ are now matrix-valued functions of a matrix variable $b$. Fortunately, they can be approximated efficiently in practice. As in the scalar case, $\omega_1$ and $\omega_2$ can be characterized as a unique fixed point of an appropriate self-map of $\bH^+(\BB)\times \bH^+(\BB)$.  Concretely, one can also solve for $\omega_1$ via iteration:
\[
w \;\longmapsto\;
h_y\bigl(\,h_x(w) + b\bigr) + b,
\]
where $h_x(b) = F_x(b) - b$ and $h_y(b) = F_y(b) - b$. In practice, this iteration converges to the subordination function quite rapidly.


\section{Operator-valued semicircle}
We will use the following definition for operator-valued (or “matrix-valued”) semicircle random variables.

\index{semicircle r.v.!operator-valued}
\begin{defi}
A \emph{matrix-valued (``multivariate'') semicircle random variable} $S$ is a linear combination of free standard semicircular r.v.s $s_i$ with self-adjoint coefficients $b_i \in \BB = M_n(\bC)$:
\begin{equation}
\label{equ_defi_semicircle}
S = \sum_{i = 1}^k b_i \otimes s_i,
\end{equation}
\end{defi}

This definition parallels that of a multivariate normal random variable, but in the context of free (noncommutative) probability.

For the multivariate semicircle $S$, we can compute its covariance function $\eta: M_n(\bC) \to M_n(\bC)$:
\[
\eta(b):= \E[S\,b\,S] \;=\; \sum_{i = 1}^{k} b_i \,b \, b_i.
\]
If the covariance function $\eta$ is known, then we can compute the Cauchy transform of $S$.

\begin{theo}
\label{theo_Cauchy_matrix_semicircle}
Let $S \in M_n(\AA)$ be a multivariate semicircle random variable with covariance function $\eta(b)$. Then the matrix Cauchy transform $G(b) := G_S(b)$ satisfies 
\begin{equation}
\label{equ_Cauchy_matrix_semicircle}
b\,G(b) \;=\; I_n + \eta\bigl(G(b)\bigr)\,G(b)
\end{equation}
for all $b \in M_n^+(\bC)$. 
\end{theo}
This remarkable theorem can be proved combinatorially by extending free cumulants to the operator-valued setting, as developed in \cite{speicher98}. For the proof, see Section 6.3 of \cite{speicher2019}. In Section \ref{sec:partial_proof} below, we will give a proof for two simpler cases, when $b = z I_n$, $\Im z > 0$, and when $b = b_0 + z I_n$, $b_0$ is self-adjoint and $\Im z > 0$. 

For $b \in M_n^+(\bC)$, one can solve the equation \eqref{equ_Cauchy_matrix_semicircle} by an appropriate iterative method. In fact, for $b = z I$, with $z \in \bC^+$,  in \cite{hfs2007} it was shown that this equation has a unique solution that satisfy property  $\Im G(z) < 0$ for all $\{z: \Im z > 0\}$, and the following algorithm has good convergence properties. Define $W(z) := i\,G(z)$. Then the equation \eqref{equ_Cauchy_matrix_semicircle} becomes 
\[
-\,i\,z\,W \;+\; W\,\eta(W) \;=\; I,
\]
which motivates iterating the map 
\begin{equation}
\label{hfs_map0}
W \;\mapsto\; \FF_z(W) \;:=\; \bigl[-\,i\,z\,I \;+\; \eta(W)\bigr]^{-1}.
\end{equation}
While this scheme already converges, a further modification with improved stability is:
\begin{equation}
\label{hfs_map}
W \;\mapsto\; \GG_z(W) \;=\; \tfrac12\Bigl(W \;+\; \FF_z(W)\Bigr).
\end{equation}

The convergence of this algorithm is based on the Earle-Hamilton Theorem that allows one to show the remarkable property that the maps \eqref{hfs_map0} and \eqref{hfs_map} are contractions in the Carath\'eodory-Reiffen-Finsler metric on certain subdomains of $M_n^+(\bC)$ and therefore converge to a fixed point.   

This provides a numerical method for computing $G_S(z)$. We can also regard $S$ as an element of the noncommutative probability space $M_n(\bC)\otimes \AA$, equipped with the scalar expectation $\tau(x) = \tr\,\E(x)$. Since $S$ is self-adjoint, it has an associated probability measure $\mu_S$, which can be obtained by defining the scalar Cauchy transform $H(z) := \tr\bigl(G_S(z)\bigr)$ and applying the usual Stieltjes inversion formula.

\begin{exa}
Consider
\[
S \;=\; i \begin{bmatrix}
0 & 2\,s_1 + s_3 & s_2\\[6pt]
-\,\bigl(2\,s_1 + s_3\bigr) & 0 & -\,s_3\\[6pt]
-\,s_2 & s_3 & 0 
\end{bmatrix},
\]
where $s_1, s_2, s_3$ are free (standard) semicircle elements. This $S$ can be written in the form \eqref{equ_defi_semicircle} with
\[
b_1 \;=\;
\begin{bmatrix}
0 & 2i & 0\\
-2i & 0 & 0\\
0 & 0 & 0
\end{bmatrix}, 
\quad
b_2 \;=\; 
\begin{bmatrix}
0 & 0 & i\\
0 & 0 & 0\\
-\,i & 0 & 0
\end{bmatrix},
\quad
b_3 \;=\;
\begin{bmatrix}
0 & i & 0\\
-\,i & 0 & -\,i\\
0 & i & 0
\end{bmatrix}.
\]
One can then apply the iterative scheme above to compute the distribution of $S$. The result of such a computation is shown in Figure~\ref{fig:exa_matrix_semicircle}.

\begin{figure}[ht]
\centering
\includegraphics[width=0.8\textwidth]{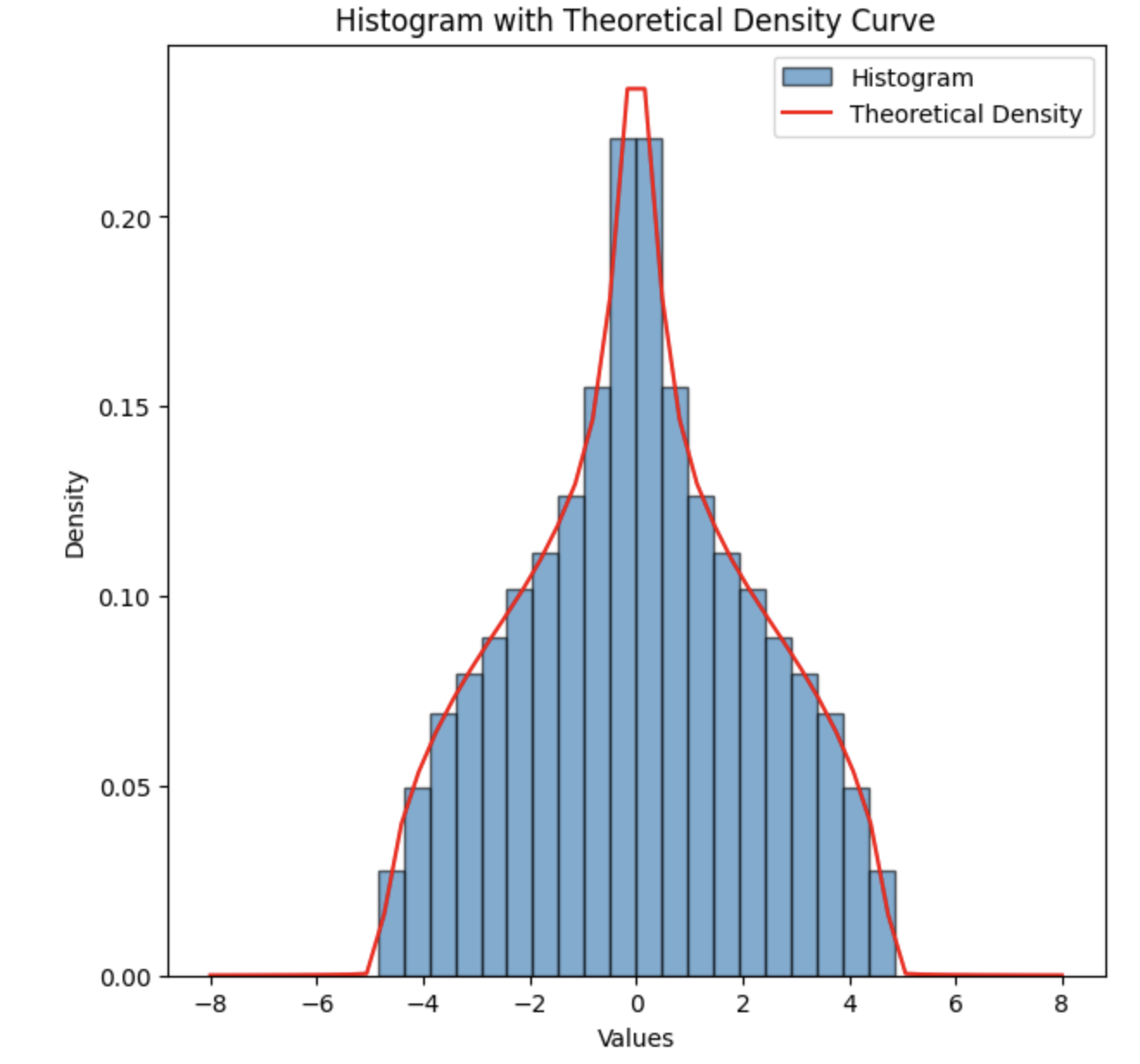}
\caption{Distribution of the matrix semicircle $S$ and its Gaussian matrix counterpart.}
\label{fig:exa_matrix_semicircle}
\end{figure}

Here, the histogram was generated by drawing eigenvalues from 10 realizations of a Gaussian matrix model in which each $s_i$ is approximated by an independent Hermitian $200\times 200$ Gaussian random matrix with entries distributed as $N(0, 1/\sqrt{200})$.
\end{exa}


\section{Linearization trick}

\index{linearization trick|textbf}
The ``linearization trick'' is used in several areas of mathematics. Its basic idea is that in order to study the spectral properties of a (non-linear) polynomial in non-commutative variables \(X_1, \ldots, X_n\), one can instead study a linear polynomial
\[
c_0 + c_1 X_1 + \cdots + c_n X_n,
\]
provided one is willing to allow the coefficients \(c_i\) to be matrices.

There are various flavors of this trick, and linearizations are typically not unique. Because of this non-uniqueness, one often strives for a version in which the dimension of the matrix is minimal, or else one that enjoys a particular structural property. 

As an elementary example, recall that the roots of a polynomial can be found (and typically are found in numerical applications) by taking the eigenvalues of its companion matrix. The linearization trick can be viewed as a more sophisticated extension of that example.

Here, we explain a version of the linearization trick which is suitable for self-adjoint polynomials (i.e., polynomials invariant under the \(\ast\)-map) and which ensures that the matrix coefficients themselves are self-adjoint.

Let \(\mathbb{C}\langle X_1, \ldots, X_n\rangle\) denote the algebra of non-commutative polynomials in variables \(X_1, \ldots, X_n\) over \(\mathbb{C}\). 

\begin{defi}\label{defi_linearization}
A matrix
\[
L \;=\; b_0 \;+\; b_1 \otimes X_1 \;+\;\cdots+\; b_n \otimes X_n,
\quad b_i \in M_N(\mathbb{C}),
\]
is called a \emph{linearization} of a non-commutative polynomial \(p \in \mathbb{C}\langle X_1, \ldots, X_n\rangle\) if \(L\) satisfies
\[
L \;=\; 
\begin{bmatrix}
0 & u \\
v & Q
\end{bmatrix},
\quad
Q \text{ is invertible in } M_{N - 1}(\mathbb{C}) \otimes \mathbb{C}\langle X_1, \ldots, X_n\rangle, 
\quad
p \;=\; -\,u\,Q^{-1}\,v.
\]
\end{defi}

The condition that \(Q\) be invertible and \(p = -\,u\,Q^{-1}\,v\) ensures that \(L\) is invertible if and only if \(p\) is invertible. This follows from the well-known Schur-complement argument, which we state after a few examples.

\begin{exa}\label{exa_X_linear}
Formally, \(X\) by itself is not a linearization of the polynomial \(p(X) = X\). Instead, we can use
\[
L \;=\; 
\begin{bmatrix}
0 & X \\
1 & -1
\end{bmatrix}.
\]
\end{exa}

\begin{exa}\label{exa_monomial_linear}
Let \(p = X_1 X_2 X_3 X_4\). A suitable linearization is
\[
L \;=\; 
\begin{bmatrix}
0 & 0 & 0 & X_1 \\
0 & 0 & X_2 & -1 \\
0 & X_3 & -1 & 0 \\
X_4 & -1 & 0 & 0
\end{bmatrix}.
\]
Indeed,
\[
Q \;=\; 
\begin{bmatrix}
0 & X_2 & -1 \\
X_3 & -1 & 0 \\
-1 & 0 & 0
\end{bmatrix}
\]
is invertible in \(M_{3}(\mathbb{C}) \otimes \mathbb{C}\langle X_2, X_3\rangle\) (it has determinant 1), and a direct calculation confirms that \(-\,u\,Q^{-1}\,v = X_1 X_2 X_3 X_4.\)
\end{exa}

This example easily generalizes to any monomial \(X_{i_1} X_{i_2}\cdots X_{i_n}\).

\begin{exa}\label{exa_polynomial_linear}
If \(p = p_1 + \cdots + p_k\), where each \(p_i\) is a monomial with linearization
\[
L_{p_i} 
\;=\; 
\begin{bmatrix}
0 & u_i \\
v_i & Q_i
\end{bmatrix},
\]
then
\[
L_p 
\;=\;
\begin{bmatrix}
0 & u_1 & \cdots & u_k \\
v_1 & Q_1 &        &     \\
\vdots &     & \ddots &     \\
v_k &     &        & Q_k
\end{bmatrix}
\]
is a linearization of \(p\).
\end{exa}

\medskip

\index{Schur-complement lemma|textbf}
We now formulate the fundamental Schur-complement lemma, used in proving the central statement about linearizations. Let \(\Lambda(z)\) be the matrix in \(M_N(\mathbb{C})\) which is zero everywhere except for its \((1,1)\)-entry, which is set to \(z\):
\[
\Lambda(z) 
:= 
\begin{bmatrix}
z & 0 & \ldots & 0\\
0 & 0 & \ldots & 0\\
\vdots & \vdots &  & \vdots\\
0 & 0 & \ldots & 0
\end{bmatrix}.
\]

\begin{theo}[Belinschi--Mai--Speicher]\label{theo_linearization}
Let 
\[
L_p 
\;=\; 
b_0 \otimes 1 
\;+\; 
b_1 \otimes X_1 
\;+\;\cdots+\; 
b_n \otimes X_n 
\;\;\in\; 
M_N(\mathbb{C}) \otimes \mathbb{C}\langle X_1,\ldots,X_n\rangle
\]
be a linearization of \(p \in \mathbb{C}\langle X_1,\ldots,X_n\rangle\). Define \(P = p(x_1,\ldots,x_n)\) and 
\[
L_P 
\;=\; 
b_0 \otimes 1 
\;+\; 
b_1 \otimes x_1 
\;+\;\cdots+\; 
b_n \otimes x_n 
\;\;\in\; 
M_N(\mathbb{C}) \otimes \mathcal{A},
\]
where \(x_1,\ldots,x_n \in \mathcal{A}\) for some unital complex algebra \(\mathcal{A}\).

Then the following two statements are equivalent:
\begin{enumerate}
\item[(i)] \(z - p\) is invertible in \(\mathcal{A}\).
\item[(ii)] \(\Lambda(z) - L_P\) is invertible in \(M_N(\mathbb{C}) \otimes \mathcal{A}\).
\end{enumerate}
Moreover, if (i) and (ii) hold, then
\begin{equation}
\label{equ_linearization}
\bigl[(\Lambda(z) - L_P)^{-1}\bigr]_{1,1}
\;=\;
(z - p)^{-1}.
\end{equation}
\end{theo}

\index{Schur-complement lemma}
The proof is based on the following well-known result about Schur complements. 
\begin{propo}
\label{propo_Schur_complement}
Let $\AA$ be a unital complex algebra and let 
\[
x = \bmatr{a & b \\ c & d} \in M_{k + l}(\AA),
\]
 where $a, b , c, d$ are $k\times k$, $k \times l$, $l \times k$, and $l \times l$ matrices, respectively, with entries in $\AA$. Assume that $d$ is invertible in $M_l(\AA)$. 
 
  Then, (i) $x$ is invertible in $M_{k + l}(\AA)$ if and only if (ii) the Schur complement $a - b d^{-1} c$ is invertible in $M_k(\AA)$. Moreover, if (i) and (ii) are satisfied, then 
  \[
  \Big[\bmatr{a & b \\ c & d}^{-1}\Big]_{1, 1} = (a - b d^{-1} c)^{-1},
  \]
  where $[X]_{1, 1}$ denotes the $k \times k$ upper-left corner of the matrix $X \in M_{k + l}(\AA)$. 
\end{propo}
Comment: In fact, one can write a more detailed formula:
\begin{equation}
\label{equ_Schur_general}
 \bmatr{a & b \\ c & d}^{-1} = \bmatr{0 & 0 \\ 0 & d^{-1}} + \bmatr{I \\ -d^{-1} c} (a - b d^{-1} c)^{-1}
  \bmatr{I & -bd^{-1}},
\end{equation}
which we will use later. This Proposition is a standard result in Linear Algebra so we omit the proof.

\begin{proof}[Proof of Theorem \ref{theo_linearization}]
We apply Proposition \ref{propo_Schur_complement} to matrix 
\[
\Lambda(z) - L_P = \bmatr{z & -\tilde u \\ - \tilde v& -\tilde Q}, 
\]
where $\tilde u, \tilde v, \tilde Q$ are obtained from $u, v, Q$ in Definition \ref{defi_linearization} by applying evaluation homomorphism $X_1 \to x_1, \ldots, X_n \to x_n$. 
By assumption,  $Q$ is invertible in $M_{N-1}(\bC) \otimes \bC(X_1, \ldots X_n)$, which implies that $\tilde Q$ is invertible in $M_{N - 1}(\AA)$. The Schur complement for $\Lambda(z) - L_P$ is $z + \tilde u \tilde Q^{-1} \tilde v  = z - P$ and therefore the equivalence of (i) and (ii) in Proposition \ref{propo_Schur_complement} implies the equivalence of (i) and (ii) in Theorem \ref{theo_linearization}. Moreover, formula \ref{equ_linearization} follows from the corresponding formula in Proposition \ref{propo_Schur_complement}. 
\end{proof}

%

\subsection*{Self-adjoint linearizations}

\index{linearization, self-adjoint}
The linearizations in Examples~\ref{exa_X_linear}, \ref{exa_monomial_linear}, and \ref{exa_polynomial_linear} are not necessarily self-adjoint, even if the original polynomial \(p\) is self-adjoint. We now explain how to convert such linearizations into self-adjoint ones. While this systematic procedure might not always give a matrix of minimal dimension, it works in general.

Suppose \(p\) is a self-adjoint polynomial. In particular, we can write \(p = q + q^\ast\). If 
\[
L_q 
\;=\; 
\begin{bmatrix}
0 & u\\ 
v & Q
\end{bmatrix}
\]
is a linearization of \(q\), consider the self-adjoint matrix
\[
L 
\;=\; 
\begin{bmatrix}
0 & u & v^\ast\\
u^\ast & 0 & Q^\ast\\
v & Q & 0
\end{bmatrix}.
\]
Since
\[
\begin{bmatrix}
0 & Q^\ast\\
Q & 0
\end{bmatrix}^{-1}
\;=\;
\begin{bmatrix}
0 & Q^{-1}\\
(Q^\ast)^{-1} & 0
\end{bmatrix},
\]
a direct calculation shows 
\[
-\,\begin{bmatrix}
u & v^\ast
\end{bmatrix}
\begin{bmatrix}
0 & Q^{-1}\\
(Q^\ast)^{-1} & 0
\end{bmatrix}
\begin{bmatrix}
u^\ast\\
v
\end{bmatrix}
\;=\;
-\;u\,Q^{-1}\,v \;-\; v^\ast\,(Q^\ast)^{-1}\,u^\ast
\;=\;
q + q^\ast
\;=\;
p,
\]
so \(L\) is a self-adjoint linearization of \(p\).

\subsection*{Applications to numerical algorithms}
We will next use such linearizations in a numerical procedure for finding the distribution of a self-adjoint polynomial in non-commutative variables.

Suppose \(p\) is self-adjoint and has a self-adjoint linearization 
\[
L 
\;=\; 
b_0 \;+\; b_1 \otimes X_1 \;+\;\cdots+\; b_n \otimes X_n.
\]
Introducing a small regularization parameter \(\varepsilon > 0\), define
\[
\Lambda_\varepsilon(z)
\;:=\;
\begin{bmatrix}
z & 0 & \cdots & 0\\
0 & i\varepsilon & \cdots & 0\\
\vdots & \vdots & & \vdots\\
0 & 0 & \cdots & i\varepsilon
\end{bmatrix}.
\]
By~\eqref{equ_linearization}, we can approximate the Cauchy transform of \(p\):
\begin{equation}\label{equ_linearization2}
(z - p)^{-1}
\;=\;
\lim_{\varepsilon \to 0}
\Bigl[\bigl(\Lambda_\varepsilon(z) - L_P\bigr)^{-1}\Bigr]_{1,1}
\;=\;
\lim_{\varepsilon \to 0}
\Bigl[\bigl(B_\varepsilon(z) - (b_1 \otimes x_1 + \cdots + b_n \otimes x_n)\bigr)^{-1}\Bigr]_{1,1},
\end{equation}
where \(B_\varepsilon(z) = \Lambda_\varepsilon(z) - b_0\).

The important observation is that the matrix in the square brackets can be computed via the iterative procedure described in Theorem~\ref{theo_matrix_subordination}, as long as the (matrix-valued) variables \(b_i \otimes x_i\) are free. While Theorem~\ref{theo_matrix_subordination} is stated for the sum of two matrix variables, it can be generalized to more variables simply by adding one free variable at a time. 

In order to apply Theorem~\ref{theo_matrix_subordination}, we need to compute the matrix-valued Cauchy transform of \(b_i \otimes x_i\) given that we know the scalar Cauchy transform of each \(x_i\). We discuss how to carry out that step in the next section.

%

\section{Calculation of the matrix Cauchy transform}

In order to apply the linearization technique, one must be able to compute the matrix Cauchy transform \( G_{B \otimes x}(A) \) for \( A \in M_n^+(\mathbb{C}) \). For a semicircular \( x \),  one can refer to equation in Theorem~\ref{theo_Cauchy_matrix_semicircle} and algorithms that give a solution for this equation. Below we outline the algorithm required to carry out this computation in a general case. 

Assume that \( B \otimes x \in M_n(\mathbb{C}) \otimes \mathcal{A} \) and we wish to calculate
\[
  G_{B \otimes x}(A) \;=\; \mathbb{E} \bigl[ (A \otimes 1 \;-\; B \otimes x)^{-1} \bigr],
\]
where \( A \in M_n^+(\mathbb{C}) \). In applications, \( A \) may depend on a complex parameter \( z \).

The two main ideas are:
\begin{enumerate}
\item The operators \( C \otimes 1 \) and \( I \otimes x \) commute for any \( C \), so we can choose suitable bases (in \( M_n(\mathbb{C}) \) and in \( \mathcal{A} \)) to diagonalize both operators.
\item We can use the properties of the conditional expectation to take certain matrix multiplications outside the expectation. In particular, if \( B \) is invertible, then
\[
  \mathbb{E}\bigl[ (A\otimes 1 \;-\; B \otimes x)^{-1}\bigr] 
  \;=\;
  \mathbb{E}\Bigl[ \bigl(B^{-1} A \otimes 1 \;-\; I \otimes x \bigr)^{-1}\Bigr] \; B^{-1}.
\]
By choosing a basis in which matrix \( B^{-1} A \) is diagonal, one can compute 
\(\mathbb{E}\bigl[ (B^{-1} A \otimes 1 \;-\; I\otimes x)^{-1}\bigr]\), and the problem is then essentially solved.
\end{enumerate}

A complication arises if \( B \) is not invertible. In practice one can  regularize the problem by adding small noise to $B$. However, below we explain how to handle this issue by using the Schur inversion formula.

\subsection*{Case of non-invertible \(\boldsymbol{B}\)}

Choose unitary matrices \( U_1 \) and \( U_2 \) such that 
\[
   U_1^\ast\, B \,U_2
   \;=\;
   \begin{bmatrix}
       D & 0 \\[6pt]
       0 & 0 
   \end{bmatrix},
\]
where \( D = \operatorname{diag}(\lambda_1, \dots, \lambda_r) \) is an invertible diagonal \( r \times r \) matrix, and the \(0\) blocks have appropriate dimensions. If \( B \) is self-adjoint, one can take \( U_1 = U_2 \). Importantly, \( B \) does not depend on \( z \), so \( U_1 \), \( U_2 \), and \( D \) need only be computed once.

Next, let
\[
   U_1^\ast\, A \,U_2
   \;=\;
   \begin{bmatrix}
       a_{11} & a_{12} \\[3pt]
       a_{21} & a_{22}
   \end{bmatrix},
\]
where \( a_{11} \) is an \( r \times r \) matrix. We assume \( a_{22} \) is invertible. (Typically, this generically holds when \( A \) depends on a small smoothing parameter \( \varepsilon \),  and it also holds when \( B \) is self-adjoint and \( A \) lies in the upper half-plane \( M_n^+(\mathbb{C}) \) (i.e. \(\frac{A - A^*}{2i} > 0\).) In the latter case $U^\ast A U \in  M_n^+(\mathbb{C})$ and therefore $a_{22} \in M_{n - r}^+(\bC)$, which ensures that $a_{22}$ is invertible.

We identify \( U \) with \( U \otimes 1 \) and use Schur’s formula \eqref{equ_Schur_general} to invert \( U^\ast \bigl(A\otimes 1 \;-\; B\otimes x\bigr) U \). Let
\[
   S \;=\; a_{11} \;-\; a_{12}\,a_{22}^{-1}\,a_{21}.
\]
Then
\begin{equation}\label{equ_matrix_G1}
  G_{B\otimes x}(A) 
  \;=\;
  \mathbb{E}\bigl[(A - B \otimes x)^{-1}\bigr]
  \;=\;
  U_2
  \begin{bmatrix}
     I & 0 \\[3pt]
     -\,a_{22}^{-1}a_{21} & I
  \end{bmatrix}
  \,M\,
  \begin{bmatrix}
     I & -\,a_{12}\,a_{22}^{-1} \\[3pt]
     0 & I
  \end{bmatrix}
  U_1^\ast,
\end{equation}
where
\[
   M 
   \;=\;
   \begin{bmatrix}
       \mathbb{E}\bigl[(S \;-\; D\otimes x)^{-1}\bigr] & 0 \\[4pt]
       0 & a_{22}^{-1}
   \end{bmatrix}.
\]

It remains to compute 
\(\mathbb{E}\bigl[(S - D \otimes x)^{-1}\bigr]\). Since \( D \) is invertible and diagonal, calculate $D^{-1} S$.  Assume \( D^{-1} S \) is diagonalizable. (This is true for a generic choice of \( A \); if \(D^{-1}S\) is not diagonalizable, one can still work with its Jordan form, though the resolvent formula may become more complicated.) Then we can write
\[
  D^{-1} S
  \;=\;
  V \,\operatorname{diag}(\mu_1,\dots,\mu_r)\,V^{-1}.
\]
Then
\begin{align}\label{equ_matrix_G2}
  \mathbb{E}\bigl[(S \;-\; D \otimes x)^{-1}\bigr]
  &\;=\;
  \mathbb{E}\bigl[\bigl(D^{-1} S \;-\; I \otimes x\bigr)^{-1}\bigr]\;D^{-1}
  \\
  &\;=\;
  V\,\operatorname{diag}\bigl(G_x(\mu_1), \dots, G_x(\mu_r)\bigr)\,V^{-1}\,D^{-1}. \notag
\end{align}
Combining \eqref{equ_matrix_G1} and \eqref{equ_matrix_G2} yields the desired formula for \( G_{B\otimes x}(A) \). Thus, one obtains a complete procedure for evaluating the matrix Cauchy transform in the non-invertible case as well.

\section{Proof of Theorem \ref{theo_Cauchy_matrix_semicircle}}
\label{sec:partial_proof}

We begin by proving Theorem~\ref{theo_Cauchy_matrix_semicircle} for the special case \(b = zI\) with \(\Im(z) > 0\).

Let \(X = (s_{ij})_{i,j=1}^d\) be a semicircle family, and define
\[
\sigma(ij, kl)\;:=\;\phi\bigl(s_{ij}\,s_{kl}\bigr).
\]
We also introduce the \emph{covariance mapping}
\(\eta\colon M_d(\mathbb{C}) \to M_d(\mathbb{C})\) via 
\[
\eta(B) \;=\; \mathbb{E}\,[\,X\,B\,X\,]
\;=\; \mathrm{id} \otimes \phi\!\bigl(X\,B\,X\bigr).
\]
Concretely, if \(B = (b_{kl})\), then
\[
\bigl[\eta(B)\bigr]_{ij} 
\;=\; \sum_{k,l=1}^d 
\sigma\bigl(i\,k,\;l\,j\bigr)\;b_{kl}.
\]

\medskip

\index{Wick's formula}
Using the standard non-crossing Wick-type formula for a semicircle family (cf.~\eqref{equ_free_Wick}), one obtains
\begin{align}
\bigl[\mathbb{E}(X^m)\bigr]_{ij}
&\;=\; 
\sum_{i_2,\dots,i_m=1}^d
\phi\bigl(s_{i\,i_2}\,\cdots\,s_{i_{m}\,j}\bigr)
\notag
\\
&\;=\;
\sum_{\pi \in NC_2(m)}\,\sum_{i_1,\dots,\,i_{m+1}=1}^d
\delta_{i_1,i}\,\delta_{i_{m+1},j}
\prod_{(p,q)\in\pi}
\sigma\bigl(i_p\,i_{p+1},\,i_q\,i_{q+1}\bigr),
\notag
\end{align}
where \(NC_2(m)\) is the set of all non-crossing pairings of \(\{1,\dots,m\}\).  In the language of matrix-valued free cumulants, this can be rewritten as
\begin{equation}
\label{equ_matrix_moment}
\mathbb{E}(X^m)
\;=\;
\sum_{\pi \in NC_2(m)}
\kappa_\pi,
\end{equation}
where
\[
[\kappa_\pi]_{ij}
\;=\;
\sum_{\,i_1,\dots,i_{m+1}=1}^d
\delta_{i_1,i}\,\delta_{i_{m+1},j}
\prod_{(p,q)\in\pi}
\sigma\bigl(i_p\,i_{p+1},\;i_q\,i_{q+1}\bigr).
\]

\begin{wrapfigure}{r}{0.40\textwidth}
\vspace{-3ex}
\centering
\includegraphics[width=0.40\textwidth]{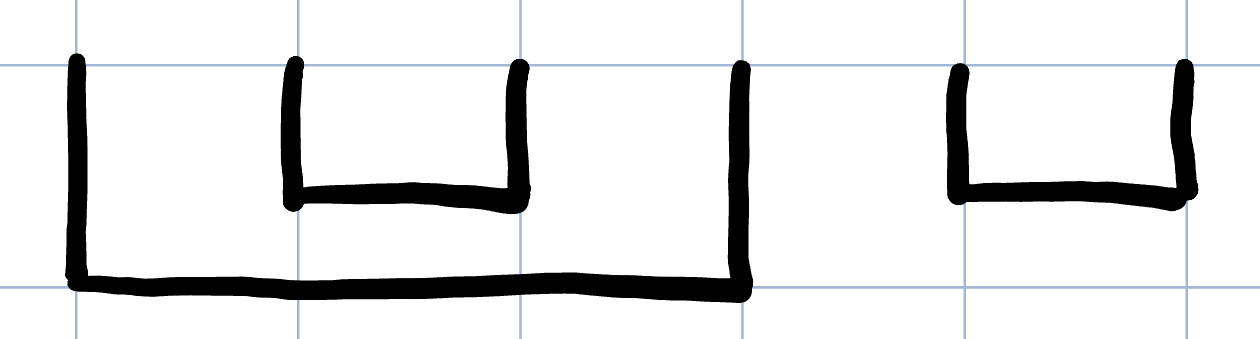}
\caption{A pairing \(\pi = \{(1,4),(2,3),(5,6)\}\).}
\label{fig:pairing142356}
\end{wrapfigure}

The matrix-valued free cumulants \(\kappa_\pi\) can be computed by iterating~\(\eta\).  As an illustration, consider 
\(\pi = \{\,(1,4),(2,3),(5,6)\}\); see Figure~\ref{fig:pairing142356}. Then
\[
[\kappa_\pi]_{ij}
\;=\;
\sum_{\,i_2,\dots,i_6}
\sigma\bigl(i\,i_2,\;i_4\,i_5\bigr)\,
\sigma\bigl(i_2\,i_3,\;i_3\,i_4\bigr)\,
\sigma\bigl(i_5\,i_6,\;i_6\,j\bigr).
\]
The pairings \((2,3)\) and \((5,6)\) factor out via sums over \(i_3\) and \(i_6\), giving
\[
[\kappa_\pi]_{ij}
\;=\;
\sum_{\,i_2,\,i_4,\,i_5}
\sigma\bigl(i\,i_2,\;i_4\,i_5\bigr)\,
\bigl[\eta(I)\bigr]_{i_2,i_4}
\;\bigl[\eta(I)\bigr]_{i_5,j}.
\]
Next, the outer pair \((1,4)\) effectively applies \(\eta\) once more, leading to
\[
[\kappa_\pi]_{ij}
\;=\;
\Bigl[\,
\eta\!\bigl(\eta(I)\bigr)\,\eta(I)\Bigr]_{ij}.
\]
Hence,
\[
\kappa_\pi
\;=\;
\eta\!\bigl(\eta(I)\bigr)\,\eta(I).
\]


More generally, a pair in $\pi$ corresponds to an application of the covariance mapping $\eta$ to the quantity recursively calculated inside the pair. Using this observation, consider the formula \eqref{equ_matrix_moment} and split the sum according to the element $r$ to which $1$ is paired. Then, using linearity of $\eta$, we obtain the following recursive formula. 
\[
\mathbb{E}(X^m)
\;=\;
\sum_{k=0}^{m-2}
\eta\bigl(\mathbb{E}(X^k)\bigr)
\;\mathbb{E}(X^{m-k-2}).
\]
Defining the moment-generating series
\[
M(z)
\;=\;
\sum_{m=0}^{\infty}
\mathbb{E}[X^m]\;z^m
\]
yields the functional equation
\[
M(z)
\;=\;
I 
\;+\;
z^2\;\eta\!\bigl(M(z)\bigr)\;M(z).
\]
By letting \(G(z) := \tfrac{1}{z}M\!\bigl(\tfrac{1}{z}\bigr)\), one obtains
\[
z\,G(z)
\;=\;
I 
\;+\;
\eta\!\bigl(G(z)\bigr)\;G(z),
\]
which is precisely \eqref{equ_Cauchy_matrix_semicircle} in the case \(b=zI\).

\paragraph{Generalization to \(b\in M_n^+(\mathbb{C})\).}
A generalization of the above argument extends to any matrix \(b\in M_n^+(\mathbb{C})\).  If we set
\[
G(b) 
\;=\; 
\mathbb{E}\bigl[(\,b - X)^{-1}\bigr],
\]
then the same non-crossing partition reasoning shows
\[
b\,G(b)
\;=\;
I
\;+\;
\eta\bigl(G(b)\bigr)\,G(b).
\]

Below is a variant of this statement that proves convenient for calculating the distribution of the \emph{biased} semicircle variables.

\begin{theo}[Biased semicircle]
\label{theo_distr_semicircle}
Let \(X\in M_d(\mathbb{C}) \otimes \mathcal{A}\) be a matrix-valued semicircle random variable with covariance map \(\eta\).  Let \(b_0 \in M_d(\mathbb{C})\) be self-adjoint, and let \(z\in\mathbb{C}^+\).  Define
\[
b(z)
\;=\;
z\,\bigl[z\,I \;-\; b_0\bigr]^{-1}.
\]
Then the Cauchy transform \(G(z) \equiv G_{\,b_0 + X}(z\,I)\) satisfies the functional equation
\begin{equation}
\label{equ_biased_semicircle}
z\,G
\;=\;
b(z)\,\Bigl[
\,I + \eta\bigl(G\bigr)\,G
\Bigr].
\end{equation}
\end{theo}

\begin{proof}
Define
\[
m_k(b)
\;=\;
\mathbb{E}\bigl[b \,(X\,b)^k\bigr],
\]
in particular, \(m_0(b)=b\).  Since $b$ and $X$ are free over $M_d(\bC)$, we can apply a modification of the argument above and obtain the recursion
\[
m_n(b)
\;=\;
b 
\sum_{\,k=0}^{\,n-2}
\eta\!\bigl(m_k(b)\bigr)\;m_{\,n-k-2}(b).
\]
The generating series
\[
M(z,b)
\;=\;
\sum_{\,n=0}^{\infty}
m_n(b)\;z^n
\]
then satisfies
\[
M(z,b)
\;=\;
b\;
\Bigl[
\,I + z^2\,\eta\bigl(M(z,b)\bigr)\,M(z,b)
\Bigr].
\]
Letting
\[
G(z,b)
\; := \;
\tfrac1z\,
M\!\bigl(\tfrac1z,\,b\bigr)
\]
translates the above into
\begin{equation}
\label{equ_Gzb}
z\,G(z,b)
\;=\;
b\,
\Bigl[
\,I + \eta\!\bigl(G(z,b)\bigr)\,G(z,b)
\Bigr].
\end{equation}
Note that
\begin{align*}
G(z, b) &= \E[b z^{-1} + b X b z^{-2} + b (X b)^2 z^{-3} + \ldots] = \E bz^{-1} [I - X b z^{-1}]^{-1} \notag
\\ 
&= \E\Big[ z b^{-1} - X\Big]^{-1}. \label{Gzb}
\end{align*}
and observe that
\begin{align*}
G_{\,b_0 + X}(z\,I)
&\;=\;
\mathbb{E}\bigl[z\,I - b_0 - X\bigr]^{-1}
\\
&\;=\;
\mathbb{E}\bigl[z b^{-1} - X\bigr]^{-1}
\;\equiv\;
G\!\bigl(z,b(z)\bigr).
\end{align*}
Hence it satisfies \eqref{equ_biased_semicircle}.
\end{proof}

The equation \eqref{equ_biased_semicircle} can be solved numerically by iterating the map
\[
G
\;\mapsto\;
\frac{1}{2}\;
\Bigl[
G
\;+\;
\bigl(z\,I - b\,\eta(G)\bigr)^{-1}\,b
\Bigr]
\]
until convergence.

\begin{theo}[Eigenvalue distribution of a polynomial in semicircle variables]
\label{theo_semicircle_polynomial}
Suppose \(p\) is a self-adjoint non-commutative polynomial in the semicircle random variables \(s_{ij}\), and let 
\[
L_p \;=\; b_0 \;+\; X
\]
be its self-adjoint \emph{linearization}, where \(b_0 \in M_d(\mathbb{C})\) and \(X\) is matrix-valued semicircle with covariance~\(\eta\).  Then the Cauchy transform
\[
G_p(z)
\;=\;
\phi\bigl[(\,z - p)^{-1}\bigr]
\;=\;
\lim_{\varepsilon \to 0}
\Bigl[
\,G_X\bigl(z,\;b_\varepsilon(z)\bigr)
\Bigr]_{1,\,1},
\]
where 
\(
b_\varepsilon(z)
\;=\;
z\;\bigl[
\Lambda_\varepsilon(z)\,-\,b_0
\bigr]^{-1}.
\)
\end{theo}

\begin{proof}
From \eqref{equ_linearization2} and the definition of the conditional expectation \(\mathbb{E}\), we have
\[
\phi\bigl[(\,z-p)^{-1}\bigr]
\;=\;
\lim_{\varepsilon\to 0}
\Bigl[
\mathbb{E}\bigl(\Lambda_\varepsilon(z)\;-\;b_0\;-\;X\bigr)^{-1}
\Bigr]_{1,\,1}.
\]
By \eqref{equ_Gzb}, this is
\[
\phi\bigl[(\,z-p)^{-1}\bigr]
\;=\;
\lim_{\varepsilon\to 0}
\Bigl[
G_X\bigl(z,\;b_\varepsilon(z)\bigr)
\Bigr]_{1,1},
\]
as claimed.
\end{proof}

\medskip

\noindent
\emph{Remark.}~In this way, one can numerically compute the eigenvalue distribution of the polynomial \(p\) by evaluating 
\(\bigl[G_X\bigl(z,b_\varepsilon(z)\bigr)\bigr]_{1,1}\) for $x + i \eps$ (i.e. closely to the real axis) and then applying the Stieltjes inversion formula to extract the distribution.

%
\section{Examples}

\begin{figure}[ht]
    \centering
    \includegraphics[width=0.6\textwidth]{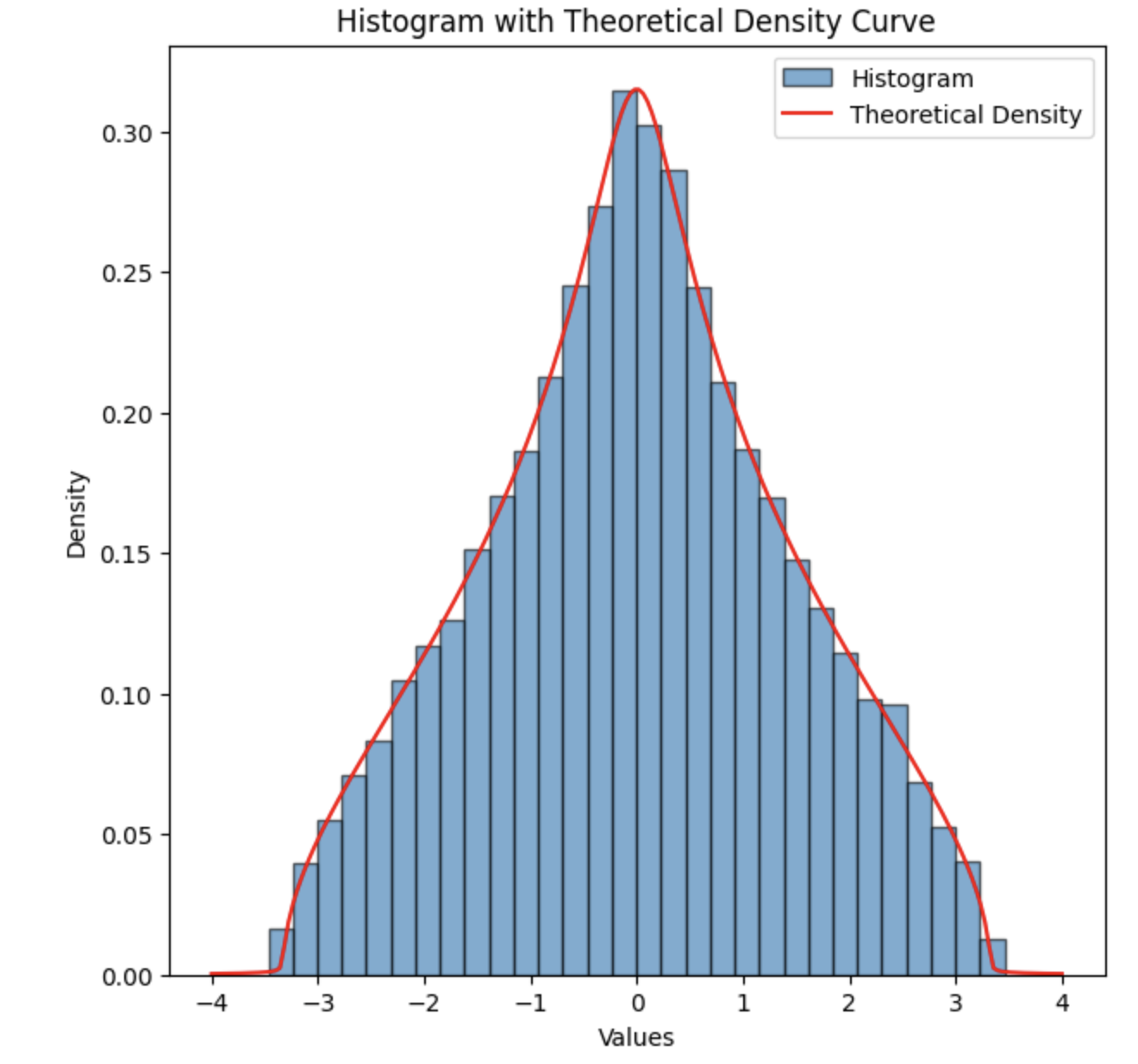}
    \caption{The distribution of the anticommutator of two semicircles and the histogram of eigenvalues of its Gaussian matrix counterpart.}
    \label{fig:exa_anticommutator}
\end{figure}

\index{anti-commutator|textbf}
\begin{exa}[Anti-commutator] 
The non-commutative polynomial $p = X_1 X_2 + X_2 X_1$ has a self-adjoint linearization 
\bal{
L &= \bmatr{ 0 & X_1 & X_2 \\
X_1 & 0 & -1 \\
X_2 & -1 & 0. 
} 
\\
&= \bmatr{ 0 & 0 & 0 \\
0 & 0 & -1 \\
0 & -1 & 0. 
} 
+ \bmatr{ 0 & 1 & 0 \\
1 & 0 & 0 \\
0 & 0 & 0. 
} \otimes X_1 
+ \bmatr{ 0 & 0 & 1 \\
0 & 0 & 0 \\
1 & 0 & 0. 
}  \otimes X_2
}
Let $s_1$ and $s_2$ are two free semicircle random variables. We can compute the distribution of $s_1 s_2 + s_2 s_1$ using equation  \eqref{equ_linearization2} and Theorem \ref{theo_matrix_subordination}. 
However, it is easier to apply Theorem \ref{theo_semicircle_polynomial}.
\end{exa}

 The result is shown in Figure 
 \ref{fig:exa_anticommutator}. 
In this figure, the density of the measure of $s_1 s_2 + s_2 s_1$ is compared with the histogram of the eigenvalue distribution for $X_1 X_2 + X_2 X_1$, where $X_1$ and $X_2$ are independent hermitian Gaussian $200\times 200$ matrices with entries distributed as $N(0, 1/\sqrt{200})$.

\begin{figure}[ht]
    \centering
    \includegraphics[width=0.6\textwidth]{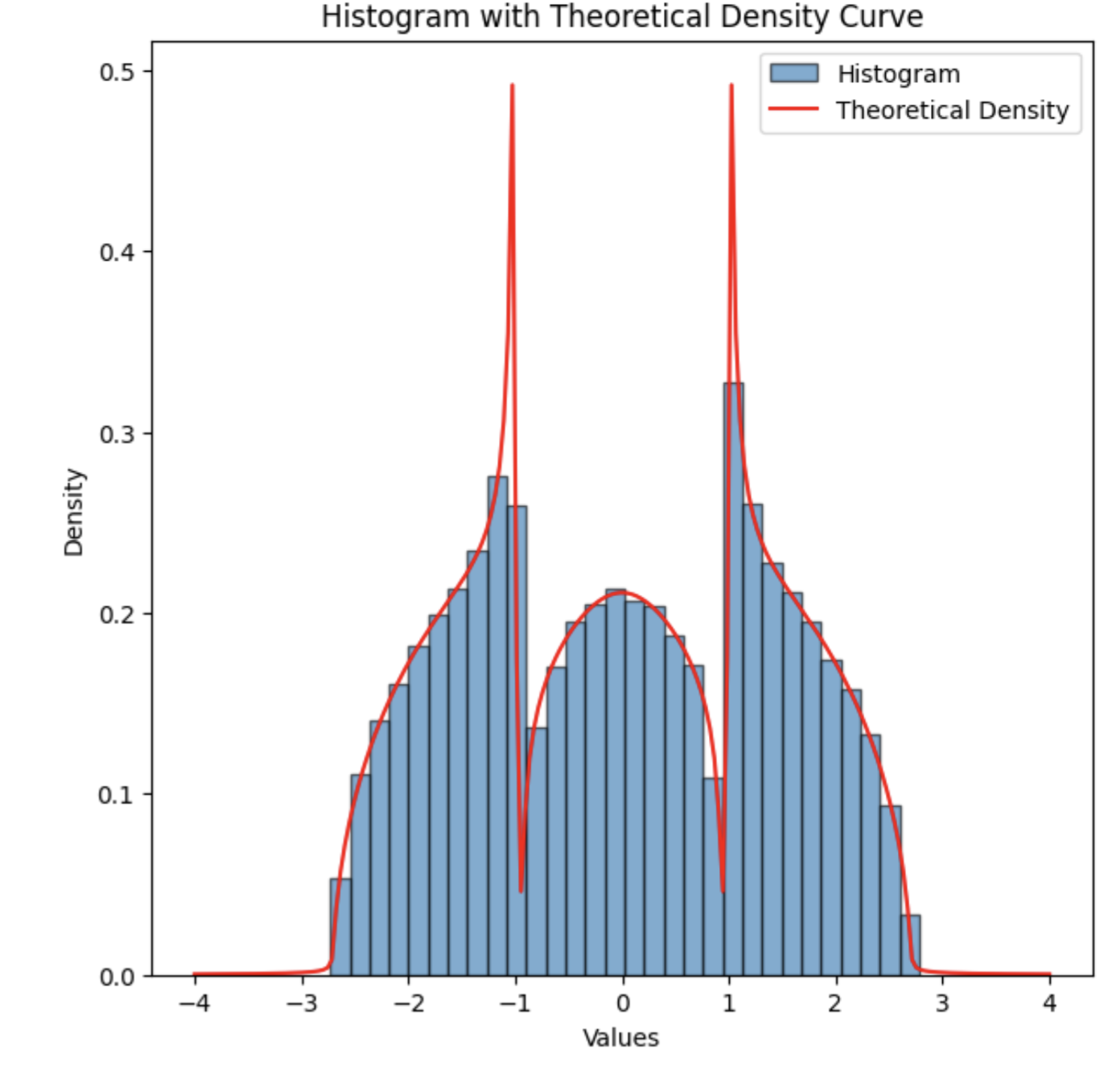}
    \caption{The distribution of a biased matrix semicircle variable and the histogram of the eigenvalues of its Gaussian matrix counterpart.}
    \label{fig:exa_biased_semicircle}
\end{figure}

\begin{exa}
Here we consider a ``biased matrix semicircle'' from the previous example  
\bal{
L &= \bmatr{ 0 & X_1 & X_2 \\
X_1 & 0 & -1 \\
X_2 & -1 & 0
}, 
}
where $X_1$ and $X_2$ are free semicircle and ask what is its distribution as a random variable in $(\AA, \tr \otimes \phi)$. The density of this distribution can be computed using Theorem \ref{theo_distr_semicircle} and the result is shown in Figure \ref{fig:exa_biased_semicircle}.
\end{exa}

\begin{figure}[ht]
    \centering
    \includegraphics[width=0.6\textwidth]{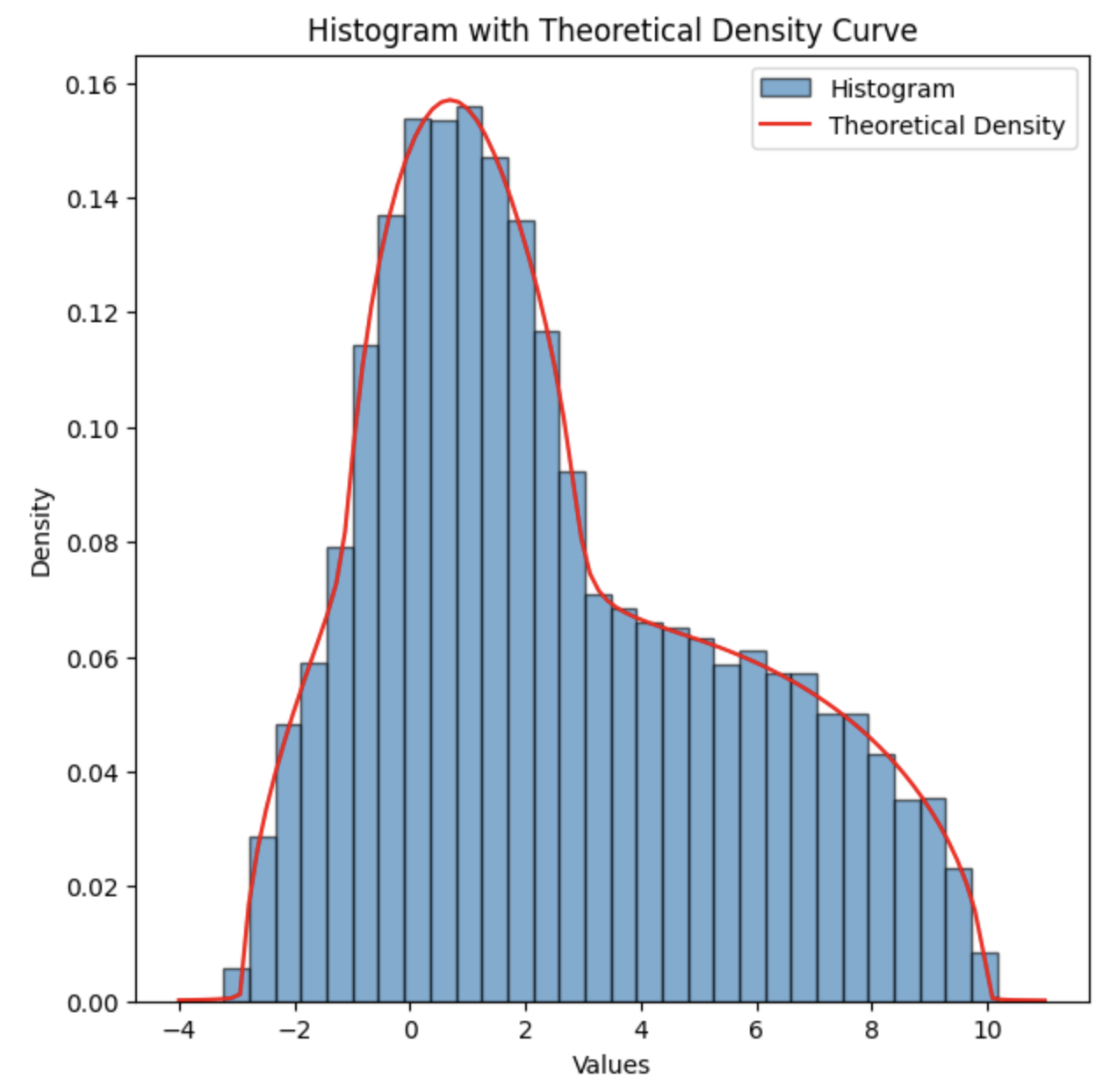}
    \caption{The distribution of a deformed anti-commutator and its Gaussian matrix counterpart.}
    \label{fig:exa_anticommutator_deformed}
\end{figure}

\index{anti-commutator!deformed}
\begin{exa}
Here we consider the deformed anticommutator
$$
p(X, Y) = X Y + Y X + X^2
$$
(This is an Example 5.2 from \cite{bms2017}).

This polynomial has a nice linearization:
$$
L = \begin{bmatrix}
0 & X & \frac{1}{2}X + Y \\
X & 0 & -1 \\
\frac{1}{2}X + Y & - 1 & 0
\end{bmatrix}
=  \begin{bmatrix}
0 & 0 & 0 \\
0 & 0 & -1 \\
0 & - 1 & 0
\end{bmatrix} +
\begin{bmatrix}
0 & 1 & \frac{1}{2} \\
1 & 0 & 0 \\
\frac{1}{2} & 0 & 0
\end{bmatrix} X +
\begin{bmatrix}
0 & 0 &  1 \\
0 & 0 & 0 \\
1 & 0 & 0
\end{bmatrix} Y.
$$

Let $X$ have the distribution $\mu_X = \frac{1}{4}(2 \delta_{-2} + \delta_{-1} + \delta_{+1})$ and $Y$ be the standard semicircle. In this situation we have to use equation  \eqref{equ_linearization2} and Theorem \ref{theo_matrix_subordination}. The result is shown in Figure \ref{fig:exa_anticommutator_deformed}.

\end{exa}

%

\section*{Notes}

This chapter draws extensively on material from \cite{bms2017} and on the Ph.D. thesis of Tobias Mai. In particular, Theorems~\ref{theo_matrix_subordination} and \ref{theo_linearization} are taken from \cite{bms2017}. A proof of Theorem~\ref{theo_Cauchy_matrix_semicircle} can be found in Chapter~9 of \cite{mingo_speicher2017} and in \cite{fobs2006}.

All of the examples presented here are also from \cite{bms2017}.

\chapter{$R$-diagonal variables}

\section{Definition and basic properties}
\index{semicircle r.v.!free cumulants}

Recall from Example~\ref{exa_semicircle_cumulants} that for the standard semicircle variable \(s\), one has
\[
\kappa_{2}(s) = 1
\quad\text{and}\quad
\kappa_{n}(s) = 0 \quad\text{for all } n \neq 2.
\]

We call a variable \(c\) a \emph{standard circle} variable if \(c\) has the same $\ast$-distribution as \(\tfrac{s_{1} + i\,s_{2}}{\sqrt{2}}\), where \(s_{1}\) and \(s_{2}\) are two standard semicircle variables.

\index{circle distribution!free cumulants}

\begin{exe}
Derive the following formulas for the free cumulants of the set \(\{c, c^{*}\}\):
\begin{itemize}
\item[(i)] \(\kappa_{n}(\dots) = 0\) for \(n \neq 2\).
\item[(ii)] \(\kappa_{2}(c,c) \;=\; \kappa_{2}(c^{*},c^{*}) = 0\).
\item[(iii)] \(\kappa_{2}(c,c^{*}) \;=\; \kappa_{2}(c^{*},c) = 1.\)
\end{itemize}
\end{exe}

A natural next question is: what are the free cumulants of a Haar unitary, i.e.\ a unitary operator whose probability distribution is uniform on the unit circle? The answer is given by the following theorem.

\index{Haar unitary!free cumulants}

\begin{theo}
Let \(u\) be a Haar unitary. Then
\[
\kappa_{2n}\bigl(u,u^{*},\ldots,u,u^{*}\bigr)
\;=\;
\kappa_{2n}\bigl(u^{*},u,\ldots,u^{*},u\bigr)
\;=\;
(-1)^{\,n-1}\,C_{n-1},
\]
where \(C_{n}\) denotes the \(n\)-th Catalan number. All other free cumulants of a Haar unitary are zero.
\end{theo}

\begin{proof}
First, note that if the number of \(u\) in the argument of a free cumulant differs from the number of \(u^{*}\), that cumulant must be zero. Indeed, we can write
\[
\kappa_{n}
\;=\;
\sum_{\pi \leq 1_n} \mu(\pi,1_{n})\,E_{\pi},
\]
and one checks that if the total number of \(u\) differs from the number of \(u^{*}\) then \(E_{\pi}=0\) for all $\pi$. 

Next, suppose there are equally many \(u\) and \(u^{*}\) but the argument of $\kappa_{2n}$ contains two consecutive \(u\) or two consecutive \(u^{*}\). We show by induction that such cumulants also vanish. For example, consider
\[
\kappa_{2n}(\dots,\,u^{*},\,u,\,u,\,\dots).
\]
(Other cases are similar.) Since \(u^{*}u=1\), we apply Theorem~\ref{theorem_cumulants_of_products} and Lemma~\ref{lemma_freeness_and_constant}, yielding
\[
0
\;=\;
\kappa_{2n-1}(\dots,\,1,\,u,\dots)
\;=\;
\sum_{\substack{\pi \in NC(2n) \\ \pi \vee \sigma = 1_{2n}}}
\kappa_{\pi}(\dots,\,u^{*},\,u,\,u,\dots),
\]
where \(\sigma\) is the partition pairing \(u^{*}\) with the first \(u\), leaving all other elements single. Let $u^{\ast }$ be in
position $m,$ and the two $u$ be in positions $m+1$ and $m+2,$ respectively.
One partition $\pi$ that connects all blocks in $\sigma $ is $1_{2n}.$ All other
non-crossing partitions with this property consist of exactly two blocks one of which
contains $m,$ and another contains $m+1.$ If the block that contains $m+1$ also contains  $m+2,$ then the cumulant $\kappa_\pi$ is zero
by inductive assumption. Therefore,  $m+2$ must be contained in the same block as $m.$ However, in this case $m+1$ is a singleton block
because the partition is non-crossing, and since $\kappa_1(u) = 0$ this implies that the cumulant $\kappa_\pi$ is also
zero. It follows that
\[
k_{2n}(\dots,u^{*},\,u,\,u,\dots)=0.
\]

We now compute
\[
\kappa_{2n}\bigl(u,u^{*},\dots,u,u^{*}\bigr).
\]
As before,
\[
0
\;=\;
\kappa_{2n-1}(1,\,u,u^{*},\dots)
\;=\;
\sum_{\substack{\pi \in NC(2n)\\ \pi \vee \sigma = 1_{2n}}}
\kappa_{\pi}\bigl(u,u^{*},u,u^{*},\dots\bigr),
\]
where \(\sigma\) pairs the first \(u\) with the first \(u^{*}\). Among the partitions contributing to the sum, one is \(1_{2n}\), and the others each having exactly two blocks, one connecting the first \(u\) with some later \(u^{*}\), and the other connecting the first \(u^{*}\) with some later \(u\). All elements in  these blocks must alternate between \(u\) and \(u^{*}\).  Let the first $u^{\ast }$ in the block to which the first $u$ connects  be in position $2p.$ (See Figure \ref{figure_haar_cumulants} for
illustration.) 

\begin{figure}[tbph]
\begin{center}
\includegraphics[width=10cm]{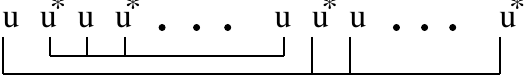}
\end{center}
\caption{Partition $\protect\pi $ in the recursion for free cumulants of
Haar unitaries}
\label{figure_haar_cumulants}
\end{figure}

This yields the recursion
\[
\kappa_{2n}
\;=\;
-\sum_{p=1}^{n}
\kappa_{2(n-p+1)}\,\kappa_{2(p-1)},
\]
which one recognizes as that of \(\,(-1)^{n-1}C_{n-1}\).
\end{proof}

A generalization of Haar unitaries based on this property is called an \emph{$R$-diagonal variable}.

\index{R-diagonal variables|textbf}
\begin{defi}
\label{definition_R_diagonal_1}
An element \(X\) in a non-commutative probability space is called \emph{$R$-diagonal} if all free cumulants of \(\{X,X^{*}\}\) vanish unless the arguments strictly alternate \(X, X^{*}, X, X^{*},\dots\) or \(X^{*}, X, X^{*}, X,\dots\).  Concretely, the only possibly non-zero cumulants are of the form
\[
k_{2n}(X,\,X^{*},\,\ldots,\,X,\,X^{*}) =: \alpha_n(X)
\quad\text{and}\quad
k_{2n}(X^{*},\,X,\ldots,X^{*},\,X) =: \beta_n(X).
\]
\end{defi}
The sequences $\alpha_n(X)$ and $\beta_n(X)$ are called the \emph{determining sequences} of $R$-diagonal element $X$. If $\alpha_n(X) = \beta_n(X)$ then $X$ is called \emph{tracial} $R$-diagonal element. This obviously always holds in the tracial probability space but may be violated in non-tracial spaces.

It turns out that $R$-diagonal variables admit two useful representations. We present the first below; the second appears in the exercises at the end of this section.

\begin{theo}\label{Theorem_R_diagonal_2}
Suppose \(U\) is a Haar unitary, \(H\) is any bounded operator, and \(U\) is free from \(H\). Then \(X = U H\) is $R$-diagonal.
\end{theo}

\begin{figure}[tbph]
\begin{center}
\includegraphics[width=7cm]{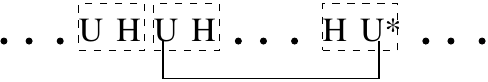}
\end{center}
\caption{}
\label{figure_R_diagonal_1}
\end{figure}

\noindent\textbf{Proof.}
We must show the free cumulants of \(X=UH\) match the $R$-diagonal requirement.  That is, we want:
\begin{itemize}
\item[(1)] all cumulants with an odd total number of arguments to vanish,
\item[(2)] any cumulant in which \(X\) or \(X^{*}\) repeats consecutively also vanishes.
\end{itemize}

\begin{figure}[tbph]
\begin{center}
\includegraphics[width=7cm]{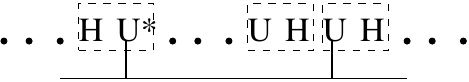}
\end{center}
\caption{{}}
\label{figure_R_diagonal_2}
\end{figure}

For example, consider the cumulant that has two consecutive $X$ in the argument $k_n( \ldots, X, X, \ldots)  = k_{n}(\dots,\,UH,\,UH,\dots)$. Using Theorem~\ref{theorem_cumulants_of_products},
\[
k_{n}(\dots, UH,\,UH, \dots)
\;=\;
\sum_{\substack{\pi \vee \sigma = 1_{2n}}}
k_{\pi}(\dots,U,\,H,\,U,\,H,\dots),
\]
where \(\sigma\) pairs each \(U\) with the \(H\) immediately after it.

Consider the block of $\pi ,$ that contains the second $U.$ If this blocks
starts with this second $U,$ then we have the situation as in Figure \ref%
{figure_R_diagonal_1}. Since $U$ is Haar unitary, the last element in this
block must be $U^{\ast }.$ Otherwise, the cumulant $k_{\pi }$ is zero.
However if the last elements in this block of $\pi $ is $U^{\ast },$ then it
is easy to see that $\pi $ can connect only those blocks of $\sigma $ which
are located between the first $U$ and the last $U^{\ast }$ of $\pi .$ In
particular, it cannot connect these blocks to the first $UH$ block depicted
in the Figure. Hence $\pi \vee \sigma \neq 1_{2n}$ and such a $\pi $ does
not enter the sum.

Next, consider the possibility that the second $U$ is not the first in the
block of $\pi $ that contains this $U.$ This situation is illustrated in
Figure \ref{figure_R_diagonal_2}. Consider then the first $U^{\ast }$ on the
left of this $U,$ which belongs to this block$.$ Such a $U^{\ast }$ must
exist by the properties of the Haar unitary elements, or we would have $%
k_{\pi }=0.$ Then again, it is clear that $\pi $ cannot connect the 
block $UH$ which is between those $U^\ast$ and $U$ to any other block of $\sigma $. Hence, $%
\pi \vee \sigma \neq 1_{2n}$ and such a $\pi $ does not enter the sum. It follows that for every $\pi$ in the sum $\kappa_\pi = 0$. 
Hence \(k_{n}(\dots,UH,\,UH,\dots)=0.\) Similar arguments cover all other cases. \(\quad\square\)

In the other direction:


\begin{theo}\label{Theorem_R_diagonal_3}
Suppose \(X\) is an $R$-diagonal element in a \textbf{tracial} $C^{\ast}$-probability space. Then \(X\) can be represented (in distribution) by \(UH\), where \(U\) is a Haar unitary, \(H\) is a positive operator having the same distribution as $ |X| := \sqrt{\,X^{*}\,X}$, and $U$ and $H$ are free. 
\end{theo}

Borrowing the terminology from linear algebra we call the probability distribution of $|X|$, $\mu_{|X|}$, the \emph{singular value distribution} of element $X$. The theorem shows that the distribution of $X$ (i.e., its determining sequence and all its *-moments) can in principle be computed from its singular value distribution. 

\begin{proof}
By Theorem~\ref{Theorem_R_diagonal_2}, \(UH\) is $R$-diagonal whenever \(U\) is a Haar unitary free from \(H\). So, to show \(X \overset{d}{=} UH\), we just need to ensure their non-vanishing cumulants coincide. Let
\[
a_{s} 
\;=\;
k_{2s}\bigl(X,X^{*},\ldots,X,X^{*}\bigr)
\;=\;
k_{2s}\bigl(X^{*},X,\ldots,X^{*},X\bigr),
\]
where the second equality uses the tracial property.

We claim that 
\begin{equation}
k_{n}\left( X^{\ast }X,\ldots ,X^{\ast }X\right) =\sum_{\substack{ \pi \in
NC\left( n\right)  \\ \pi =\left\{ V_{1},\ldots ,V_{r}\right\} }}a_{\left|
V_{1}\right| }\ldots a_{\left| V_{r}\right| }.
\label{formula_cumulants_XXstar}
\end{equation}

\begin{figure}[tbph]
\begin{center}
\includegraphics[width=7cm]{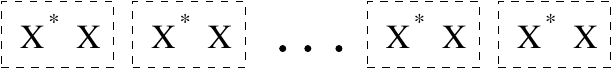}
\end{center}
\caption{{}Partition $\protect\sigma .$}
\label{figure_partition_sigma}
\end{figure}

Indeed, 
\begin{equation}
\label{sum_cumulants}
k_{n}\left( X^{\ast }X,\ldots ,X^{\ast }X\right) =\sum_{\pi \vee \sigma
=1_{2n}}k_{\pi }\left( X^{\ast },X,\ldots ,X^{\ast },X\right) ,
\end{equation}%
where as usual the sum is over non-crossing $\pi $ and $\sigma $ is the
partition which pairs $X^{\ast }$ with the following $X$ (see Figure \ref%
{figure_partition_sigma}).

Let us consider only those $\pi$ in this sum for which $\kappa_\pi \ne 0$. The blocks of this partitions must connect $X$ and $X^\ast$ in alternating order. We claim that all $\pi$ with $\kappa_\pi \ne 0$ in this sum can be put in a bijection with the set of non-crossing partitions of $\{1, 2, \ldots, n\}$ and that $\kappa_\pi = a_{|
V_{1}| }\ldots a_{| V_{r}|}$, where $V_{1}, \ldots, V_{r}$ are blocks of the partition in $NC(n)$ which corresponds to $\pi$.

\begin{figure}[th]
\begin{center}
\includegraphics[width=7cm]{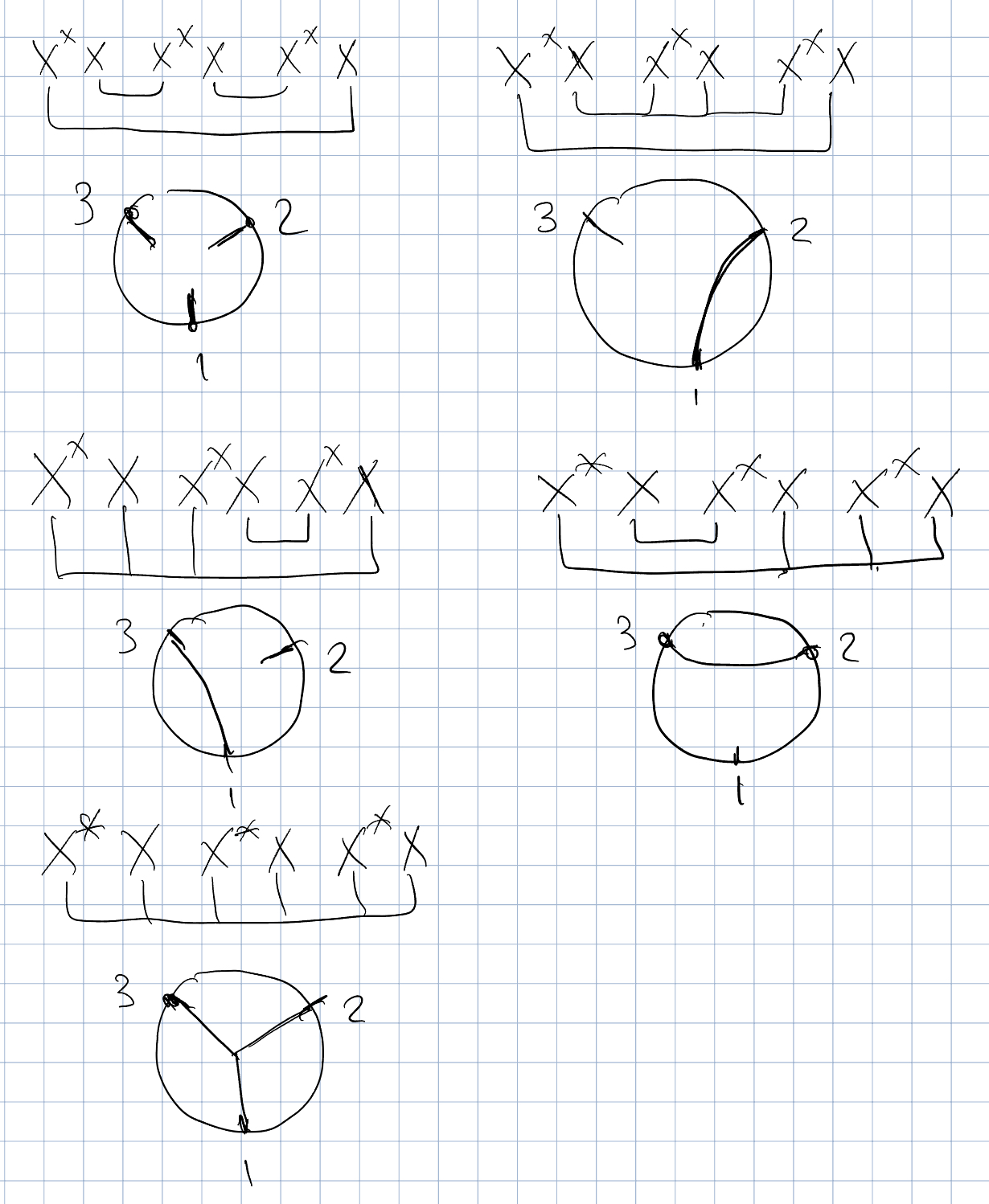}
\end{center}
\caption{Bijection between $\pi$ in the sum \eqref{sum_cumulants}, for which $\kappa_\pi \ne 0$, and elements of $NC(n)$, for $n = 3$.}
\label{fig:bijection_R_diagonal}
\end{figure}

The bijection is illustrated in Figure \ref{fig:bijection_R_diagonal}.

Here is a description of the bijection. We think about the set $\{1, 2, \ldots, n \}$ as arranged on the circle. $k$-th $X$ and $k +1$-st $X^{\ast }$ will correspond to element $k$ for $1 \leq k \leq n - 1$ and $n$-th $X$ and $1$-st $X^\ast$ will correspond to element $n$. We claim that a $\pi$ with the property described above correspond to a non-crossing partition of $\{1, 2, \ldots, n\}$. (and the corresponding cumulants is $a_{\left| V_{1}\right|
}\ldots a_{\left| V_{r}\right| },$ where $V_{i}$ are the blocks of the
resulting partition from $NC\left( n\right) .$ 

Indeed, consider a block $b$ of the partition $\pi.$ Suppose for every $X$ that belongs
to $b$, $X$ is connected to an $X^{\ast }$ on the right. Then it
must be that this $X^{\ast }$ is the immediate neighbor of $X$ on the right.
Otherwise, the blocks of $\sigma $ between these $X$ and $X^{\ast }$ would
be disconnected from $b$. It follows that this block corresponds to a block of a partition of $\{1, \ldots, n\}$.

Next suppose that some $X$ of $b$ is not connected to any $X^{\ast }$ on the right.
Then, if the cumulant is not-zero, then it must be that this is the last $%
X $ in $b$ and that $b$ starts with $X^{\ast }.$ Then, it is clear that in
this case all blocks of $\sigma $ on the right of this last $X$ and all
blocks of $\sigma $ on the left of $X^{\ast }$ are disconnected from block $%
b.$ Hence, the only possible case is when this $X$ is
the $X$ in the position $2n$ and the corresponding $X^{\ast }$ is in the
position $1$. It follows that this block of $\pi$ also corresponds to a valid block of a partition of $\{1, \ldots, n\}$. We can also observe that the resulting partition is non-crossing.  

One can also see that conversely every partition of $\{1, \ldots,  n\}$ is in correspondence with $\pi$ that connects all blocks of $\sigma$. This proves the bijection and the formula \ref{formula_cumulants_XXstar}.


We can re-write formula \ref{formula_cumulants_XXstar} as follows: 
\begin{equation*}
k_{n}\left( X^{\ast }X,\ldots ,X^{\ast }X\right) =a_{n}+\sum_{\substack{ \pi
\neq 1_{n}  \\ \pi =\left\{ V_{1},\ldots ,V_{r}\right\} }}a_{\left|
V_{1}\right| }\ldots a_{\left| V_{r}\right| }.
\end{equation*}%
This form of the formula makes it clear that we can calculate the cumulants $%
a_{n}$ from the sequence \ of $k_{n}\left( X^{\ast }X,\ldots ,X^{\ast
}X\right) .$


%
%
Hence the distribution of \(X^{*}X\) completely determines the sequence \(\{a_{s}\}\).  But \((UH)^{*}(UH)=H^{2}\) also has that same distribution, and $UH$ is $R$-diagonal by Theorem \ref{Theorem_R_diagonal_2}, so \(\{UH,(UH)^{*}\}\) and \(\{X,X^{*}\}\) must share identical non-vanishing cumulants.  Therefore, \(UH\) and \(X\) have the same \(*\)-distribution. 

\end{proof}

In some sense, $R$-diagonal can be thought of as non-commutative non-selfadjoint analogues of symmetric random variables. Here are some results about $R$-diagonal variables formulated as exercises. See Notes section for the sources of these results. 

Theorem \ref{Theorem_R_diagonal_2} can be generalized as follows. 
\begin{exe}[\textbf{R-diagonalization}]
\label{exe_SA_Rdiagonal_0}
 Let $a$ and $b$ be elements in a *-probability space $(\AA,\phi)$ such that $a$ is $R$-diagonal and $a$ and $b$ are *-free (i.e., $\{a,a^\ast\}$ and $\{b, b^\ast\}$ are freely independent). Then $a b$ is $R$-diagonal.
 \end{exe}
 
 \begin{exe}[\textbf{Characterization of $R$-diagonal elements through unitary invariance}]
 \label{exe_SA_Rdiagonal_0.5}
 Let $a$ be an element in *-probability $(\AA, \phi)$ and $u$ be a Haar unitary in $(\AA, \phi)$ such that  $a$ and $\{u, u^\ast\}$ are free. Then $a$ is R-diagonal if and only if $a$ has the same *-distribution as $ua$.   
 \end{exe}
 
Exercises  \ref{exe_SA_Rdiagonal_0.5} and \ref{exe_conjugation_by_unitary} can be used to prove the following result. 
 \begin{exe}[\textbf{Surprising freeness}]
  \label{exe_SA_Rdiagonal_0.75}
Let $a$ be $R$-diagonal. Then $a a^\ast$ and $a^\ast a$ are freely independent. 
\end{exe}

\index{even random variable|textbf}
Somewhat surprisingly, there is another representation for $R$-diagonal elements. A self-adjoint random variable is called \emph{even} if it has a symmetric distribution


\begin{exe}[\textbf{Product of free even elements is $R$-diagonal.}]
\label{exe_product_even}
Let $x$ and $y$ be two free even self-adjoint random variable in a *-probability space. Then $x y$ is a tracial $R$-diagonal element. Moreover, for the determining sequence we have
\[
\alpha_n(xy) = \sum_{\substack{\pi, \sigma \in NC(n) \\ \sigma \leq K(\pi)}} \alpha_\pi(x) \alpha_\sigma(y), 
\]
where $\alpha_n(x) = \kappa_{2n}(x)$ and  $\alpha_\pi(x)$ are defined by multiplicativity.
\end{exe}

\begin{exe}[\textbf{Representation using free symmetry}]
\label{exe_SA_Rdiagonal_2}
Show that every $R$-diagonal operator $b$ can be represented (in distribution) by \(a x\), where \(x\) if a self-adjoint even element with $\alpha_n(x) = \alpha_n(b)$ (that is, $b$ and $x$ have the same determining sequence), $a$ is a self-adjoint even element with distribution \(\frac12(\delta_{-1}+\delta_{1})\) and $a$ and $x$ are free. 

Show that the probability distribution of $x$ is the symmetrization of the singular value distribution of $b$. 
\end{exe}

It is not hard to see from these results that the product or the sum of two free $R$-diagonal elements remains $R$-diagonal. 

Indeed, Exercise \ref{exe_SA_Rdiagonal_0} implies that a product of two *-free $R$-diagonal random variables $x$ and $y$ is $R$-diagonal. Moreover, one can calculate the determining sequence of the product. 

\begin{exe}[\textbf{Product of $R$-diagonal random variables}]
\label{exe_R_diagonal_product}
Let $x$ and $y$ be two *-free $R$-diagonal elements in a tracial *-probability space $(\AA, \phi)$. 
Then
\[
\alpha_n(xy) = \sum_{\substack{\pi, \sigma \in NC(n) \\ \sigma \leq K(\pi)}} \alpha_\pi(x) \alpha_\sigma(y).
\]
\end{exe}

Now, for the sum of *-free $R$-diagonal random variables, we can use the definition and multilinearity of cumulants. For example. 
\[
\kappa_2(a + b, a + b) = \kappa_2(a, a) + \kappa_2(b,b) + \kappa_2(a,b)  + \kappa_2(b, a) = 0,
\]
where the first two terms vanish because $a$ and $b$ are R-diagonal and the last two terms vanish because $a$ and $b$ are free. Generalizing this argument we can show the result in the following exercise.

\begin{exe}[\textbf{Sum of $R$-diagonal random variables}]
\label{exe_R_diagonal_sum_0}
Let $a$ and $b$ be two $R$-diagonal elements which are *-free. Show that $a + b$ is also $R$-diagonal. Calculate its determining sequence. 
\end{exe}

Using Exercises  \ref{exe_SA_Rdiagonal_2} and \ref{exe_R_diagonal_sum_0}, we can also calculate the singular value distribution of the sum. 
\begin{propo}[Singular value distribution of a sum of $R$-diagonals.]
\label{propo_R_diagonal_sum_1}
Let $a$ and $b$ be two *-free $R$-diagonal elements in a tracial $C^\ast$ probability space. Then, 
\[
\tilde \mu_{|a + b|} = \tilde \mu_{|a|} \boxplus \tilde \mu_{|b|}, 
\]
where $\tilde \mu$ denotes the symmetrization of  the measure $\mu$: $\tilde \mu$ is a symmetric measure with the property $\mu(A) + \mu(-A) = \tilde \mu(A) + \tilde \mu(-A)$.  
\end{propo}

\begin{proof}
Let $x$ and $y$ be two even elements that represent $a$ and $b$ according to Exercise  \ref{exe_SA_Rdiagonal_2}. We can assume them to be free of each other. Then the statement about the determining sequence of $a + b$ in Exercise \ref{exe_R_diagonal_sum_0} implies that $a + b$ is represented by the even element $x + y$ and the conclusion of the theorem follows. 

\end{proof}

In addition, powers of an $R$-diagonal element also remain $R$-diagonal. To prove that fact, we first prove a general lemma about $R$-diagonal products.

\begin{lemma}\label{lemma_products}
Let $a_{1}, a_{2}, \dots, a_{n}$ be free $R$-diagonal variables in a $C^{*}$-probability space $\mathcal{A}_{1}$.  Let $A_{1}, A_{2}, \dots, A_{n}$ be positive self-adjoint random variables in another probability space $\mathcal{A}_{2}$, each $A_{i}$ having the same distribution as $\sqrt{\,a_{i}^{*}\,a_{i}}$.  Suppose $U$ is a Haar unitary in $\mathcal{A}_{2}$, free from $\{A_{1},\dots,A_{n}\}$.  Set
\[
\Pi \;=\; a_{n}\,\cdots\,a_{1}
\quad\text{and}\quad
X \;=\; U\,A_{n}\,\cdots\,U\,A_{1}.
\]
Then \(\Pi\) and \(X\) have the same \(*\)-distribution (\(\Pi \cong X\)).
\end{lemma}
In general we will write $Y \cong X$ iff $Y$ and $X$ have the same *-distribution. The key point is that $A_i$ are not necessarily free from each other and that at each alternating place it is the same $U$ in the product, not $U_1, \ldots, U_n$.

\begin{proof}
By Theorem \ref{Theorem_R_diagonal_3},
\[
\Pi
\;\cong\;
u_{n}\,\lvert a_{n}\rvert
\;\cdots\;
u_{1}\,\lvert a_{1}\rvert,
\]
where \(u_{i}\) are Haar unitaries, all *-free from each other and from \(\{\lvert a_{1}\rvert,\dots,\lvert a_{n}\rvert\}\).  

The sequence $(|a_1|, \ldots , |a_n|)$ has the same joint distribution as the sequence $(v_i A_i v_i^\ast)$, where $v_i$ are *-free Haar unitaries, free of $\{A_1, \ldots, A_n, u_1, \ldots, u_n\}$.   

We obtain 
\begin{eqnarray*}
\Pi &\cong &u_{n}v_{n}A_{n}v_{n}^{\ast }u_{n-1}v_{n-1}A_{n-1}\ldots
u_{1}v_{1}A_{1}v_{1}^{\ast } \\
&\cong &v_{1}^{\ast }u_{n}v_{n}A_{n}v_{n}^{\ast }u_{n-1}v_{n-1}A_{n-1}\ldots
u_{1}v_{1}A_{1}.
\end{eqnarray*}%
The variables $U_{i}:=v_{i+1}^{\ast }u_{i}v_{i}$ are
Haar-distributed. (Here $v_{n+1} = v_{1}.$) They are also
free. 

Indeed, we can replace each of $u_{i}$ with a product of two free
Haar-unitaries, $u_{i}^{\prime }$ and $u_{i}^{\prime \prime }.$ This will
not change the distribution. Then variables $u_{i}^{\prime \prime }v_{i}$
are $\ast $-free from all $v_{i}$. (This holds by Theorem 1 in \cite{ryan98}%
) The same is true for $v_{i+1}^{\ast }u_{i}^{\prime }.$ This implies that $%
v_{i+1}^{\ast }u_{i}^{\prime }u_{i}^{\prime \prime }v_{i}$ is $\ast $-free
from all $v_{j+1}^{\ast }u_{j}^{\prime }u_{j}^{\prime \prime }v_{j}$ for $%
j\neq i.$

Therefore,
\[
\Pi
\;\cong\;
U_{n}\,A_{n}\,U_{n-1}\,A_{n-1}\,\cdots\,U_{1}\,A_{1},
\]
where \(U_{1},\dots,U_{n}\) are Haar unitaries, mutually free and free from $\{A_1, \ldots A_n\}$.

A crucial step is to verify that this has the same \(*\)-distribution as
\[
X
\;=\;
U\,A_{n}\,U\,A_{n-1}\,\cdots\,U\,A_{1},
\]
where $U$ is a Haar unitary, which is free from all $A_{i}.$

Let $U^{i}$ denote the identical copies of $U$ with superscript $i$ showing
the position of this copy of $U$ in the product $UA_{n}UA_{n-1}\ldots
UA_{1}. $ Hence $U^{i}=U^{j}=U$ for all $i$ and $j,$ and we write: 
\begin{eqnarray*}
X &=&U^{n}A_{n}U^{n-1}A_{n-1}\ldots U^{1}A_{1}, \\
X^{\ast } &=&A_{1}U^{1\ast }A_{2}U^{2\ast }\ldots A_{n}U^{n\ast }
\end{eqnarray*}%
We are interested in a formula for the $\ast $-moments of $X$. More
particularly, we want to show that this formula is exactly the same as the
formula for the corresponding moments of $\Pi .$

Let us explain the reasoning by an example. Consider the following $\ast $%
-moment: 
\begin{eqnarray*}
&&E\left( X^{\ast }XXX^{\ast }\right) \\
&=&E\left( \underset{X^{\ast }}{\underbrace{A_{1}U^{1\ast }A_{2}U^{2\ast
}\ldots A_{n}}}\underset{X}{\underbrace{A_{n}\ldots U^{1}A_{1}}}\underset{X}{%
\underbrace{U^{n}A_{n}\ldots U^{1}A_{1}}}\underset{X^{\ast }}{\underbrace{%
A_{1}U^{1\ast }\ldots A_{n}U^{n\ast }}}\right)
\end{eqnarray*}%
By applying formula from Theorem \ref{theorem_expectation_products}, we
obtain the following expression: 
\begin{eqnarray}
&&E\left( X^{\ast }XXX^{\ast }\right)  \label{formula_expansion} \\
&=&\sum_{\pi \in NC\left( 4n-2\right) }[\kappa _{\pi }\left( U^{1\ast
},\ldots ,U^{(n-1)\ast },U^{n-1},\ldots ,U^{1},U^{n},\ldots ,U^{1},U^{1\ast
},\ldots ,U^{n\ast }\right) \times  \notag \\
&&\times E_{K\left( \pi \right) }\left( A_{2},\ldots ,\left( A_{n}\right)
^{2},\ldots ,A_{1},A_{n},\ldots ,\left( A_{1}\right) ^{2},\ldots
,A_{n},A_{1}\right) ].  \notag
\end{eqnarray}

All these formulas are the same for $E\left( \Pi ^{\ast }\Pi \Pi \Pi ^{\ast
}\right) $ except that we have to use variables $U_{j}$ instead of $%
U^{j}$ throughout these expressions.

We claim that if a block of a partition $\pi $ connects a variable from the
subset $\left\{ U^{i},U^{i\ast }\right\} $ with a variable from the subset $%
\left\{ U^{j},U^{j\ast }\right\} $ (where $i\neq j)$, then the cumulant $%
\kappa _{\pi }$ is equal to zero. Note that this is true for variables $%
U_{i} $ and $U_{j},$ because $U_{i}$ is $\ast $-free from $U_{j}$ by
assumption. However, the validity of this claim for $U^{i}$ and $U^{j}$
needs a proof because $U^{i}=U^{j}.$

In order to explain this fact, we need an additional layer of notation.
Namely, we will use a subscript to describe a position of an $X$ in the
product. For example, we will write $E\left( X_{1}^{\ast
}X_{2}X_{3}X_{4}^{\ast }\right) $ instead of $E\left( X^{\ast }XXX^{\ast
}\right) .$ It will be understood, however, that all $X_{i}$ are identical.
Correspondingly, we use notation $U_{\alpha }^{i}$ and $U_{\alpha }^{i\ast }$
in order to specify the position of variables $U$ and $U^{\ast }.$ For
example, $U_{3}^{2\ast }$ means the second $U^{\ast }$ in the third $X.$

In addition, we will use variables $\varepsilon \in \left\{ \emptyset ,\ast
\right\} $ in order to have a uniform notation for $U$ and $U^{\ast }.$ That
is, $U_{\alpha }^{i,\varepsilon }$ means $U_{\alpha }^{i}$ if $\varepsilon
=\emptyset ,$ and it means $U_{\alpha }^{i\ast }$ if $\varepsilon =\ast .$

\begin{figure}[tbph]
\begin{center}
\includegraphics[width=7cm]{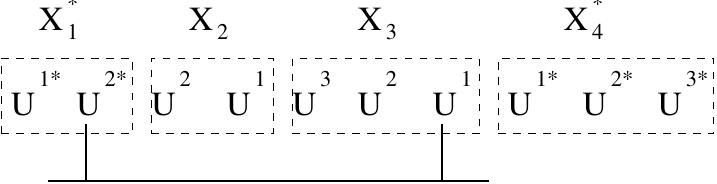}
\end{center}
\caption{{}Example of $k_{\protect\pi }.$}
\label{figure_R_diagonal_3}
\end{figure}

Suppose that a partition $\pi $ has a block $\mathfrak{b,}$ such that
variables $U_{\alpha }^{i,\varepsilon _{1}}$ and $U_{\beta }^{j,\varepsilon
_{2}}$ belong to $\mathfrak{b}$ and that $i\neq j.$ In addition, suppose
that among all the pairs of variables $U_{\alpha }^{i,\varepsilon _{1}}$ and 
$U_{\beta }^{j,\varepsilon _{2}}$ that satisfy this property we choose a
pair of variables with the smallest possible distance between them. (Here,
the distance between variables is understood as the difference in their
positions in the sequence $U_{1}^{1,\varepsilon _{1}},\ldots
,U_{m}^{n,\varepsilon _{m}},$ which is the argument of the cumulant $\kappa
_{\pi }.$) This choice implies, in particular that $U_{\alpha
}^{i,\varepsilon _{1}}$ and $U_{\beta }^{j,\varepsilon _{2}}$ are neighbors
in the block $\mathfrak{b}$, that is, that there are no variables between
them that belong to the same block $\mathfrak{b}$.

An example is shown in Figure \ref{figure_R_diagonal_3}, where $\alpha =1,$ $%
\beta =3,$ $i=2,$ $j=1.$

If $\varepsilon _{1}=\varepsilon _{2},$ then $\kappa _{\pi }=0$ because $U$
is a Haar unitary and a cumulant of a Haar unitary is non-zero if and only
if the sequence of $U$ and $U^{\ast }$ in its argument is alternating and
has the same number of $U$ and $U^{\ast }$.

Hence, we can assume that $\varepsilon _{1}\neq \varepsilon _{2}$ and $%
\alpha <\beta .$ Consider the case when $\varepsilon _{1}=\ast $ and $%
\varepsilon _{2}=\emptyset .$ (The other case is similar.)

Let us count those $U^{\ast }$ in the argument of the cumulant $k_{\pi
}$ which are between $U_{\alpha }^{i,\ast }$ and $U_{\beta }^{j}$ and which
come from $X_{\alpha }^{\ast }.$ Their number is 
\begin{equation*}
c_{2}\left( \ast \right) =\left\{ 
\begin{array}{cc}
n-i, & \text{if }\varepsilon _{\alpha +1}=\ast , \\ 
n-1-i, & \text{if }\varepsilon _{\alpha +1}=\emptyset .%
\end{array}%
\right.
\end{equation*}
In our example this is $3 - 2 - 1 = 0$. Note that $U^{3 \ast}$ in $X_1^\ast$ was cancelled out by $U^3$ in $X_2$. 

Similarly, we can count those $U$ that are between $U_{\alpha }^{i,\ast }$
and $U_{\beta }^{j}$ and that come from $X_{\beta }.$ This count is 
\begin{equation*}
c_{2}\left( \emptyset \right) =\left\{ 
\begin{array}{cc}
n-1-j, & \text{if }\varepsilon _{\beta -1}=\ast , \\ 
n-j, & \text{if }\varepsilon _{\beta -1}=\emptyset .%
\end{array}%
\right.
\end{equation*}
In our example, this is $3 - 1 = 2$. 

Let $S$ be the sequence of $X_{\alpha +1}^{\varepsilon \left( \alpha
+1\right) },\ldots ,X_{\beta -1}^{\varepsilon \left( \beta -1\right) }.$ (If 
$\alpha +1>\beta -1,$ then $S$ is empty.) That is, we look on the sequence
of $X$ and $X^{\ast }$ between $X_{\alpha }^{\ast }$ and $X_{\beta }$ in the
defining sequence of the moment. In our example $S = (X_2).$

 Let $p$ be the number of substrings $X^\ast X$ in $S$,  
$s^{\ast }$ and $s$ be the number of $X^\ast$ and $X$, respectively, that do not belong to these substrings. 

In our example, $p = s^\ast = 0$, $s = 1$. For the sequence 
\[
S = X^\ast X X^\ast X^\ast X X X X^\ast X^\ast X X^\ast X^\ast,  
\]
we have $p = 3$, $s^\ast = 4$, and $s = 2$ and the length of the sequence is $l = 2 p +  s^\ast + s = 2 \times 3 + 4 +  2 = 12$.


Let us now count variables $U^{\ast }$ in the argument of the cumulant $%
k_{\pi }$ that arise from this sequence. (Note that some variables $U^{n\ast
}$ are cancelled out when the product $X^{\ast }X$ occurs.). By taking into account that we have $X^\ast$ befor and $X$ after the sequence, it is easy to
count that the number of $U^{\ast }$ in $S$ is 
\begin{equation*}
c_{1}( \ast) =\left\{ 
\begin{array}{cc}
( n-1) p + n s^\ast, & \text{if }\varepsilon _{\beta
-1}=\emptyset , \\ 
( n-1) p + n s^\ast - 1 , & \text{if 
}\varepsilon _{\beta -1}=\ast .%
\end{array}%
\right.
\end{equation*}%
A similar count for variables $U$ gives 
\begin{equation*}
c_{1}( \emptyset ) =\left\{ 
\begin{array}{cc}
( n-1) p + ns, & \text{if }\varepsilon _{\alpha +1}=\ast , \\ 
( n-1) p +n s-1  , & \text{if }%
\varepsilon _{\alpha +1}=\emptyset .%
\end{array}%
\right.
\end{equation*}

We intend to compute the difference between the number of $U^{\ast }$ and $U$
which are between $U_{\alpha }^{i,\ast }$ and $U_{\beta }^{j}.$ This number
is $\Delta :=c_{1}\left( \ast \right) +c_{2}\left( \ast \right) -c_{1}\left(
\emptyset \right) -c_{2}\left( \emptyset \right) .$ We need to consider the
four possible combinations of $\varepsilon _{\alpha +1}$ and $\varepsilon
_{\beta -1}.$

For example, if $\varepsilon _{\alpha +1}=\varepsilon _{\beta -1}=\emptyset
, $ then 
\begin{align*}
\Delta &= ( n-1) p +ns^\ast +n-1-i \\
&-[ ( n-1) p +n s - 1 +n-j]
\\
&=n\left( s^{\ast }-s\right) +j-i.
\end{align*}%
It turns out that the result is the same for all other combinations.

Since $1\leq i,j\leq n$ and $i\neq j,$ therefore $\Delta \neq 0.$

To summarize, the number of $U^{\ast }$ is always different from the number
of $U$ between $U_{\alpha }^{i\ast }$ and $U_{\beta }^{j}$. These $U^{\ast }$
and $U$ cannot be connected to a $U^{\ast }$ or $U$ outside of $U_{\alpha
}^{i\ast }$ and $U_{\beta }^{j}$ because this would result in a crossing
partition. Hence, all these $U^{\ast }$ and $U$ are split by blocks of $\pi
, $ which contain only these $U^{\ast }$ and $U.$ It follows that one of
these blocks must have an unequal number of $U^{\ast }$ and $U.$ This
implies that $\kappa _{\pi }=0$ because $u$ is a Haar unitary.

If we write a formula for a $\ast $-moment of $\Pi ,$ which is similar to
formula (\ref{formula_expansion}), then we get the same expression as in ( %
\ref{formula_expansion}), except that we have to write $U_{i}$, $U_{i}^{\ast
}$ instead of $U^{i}$ and $U^{i\ast },$ respectively. In this case, if a
block of a partition $\pi $ connects to variables with different indices,
then the cumulant $\kappa _{\pi }$ is zero by assumption of $\ast $-freeness
of $U^{i}$ and $U^{j}$ for different $i$ and $j.$ The remaining cumulants of
variables $U_{i}$ coincide with the corresponding cumulants of variables $%
U^{i}$ because all $U_{i}$ and $U^{i}$ have the same Haar distribution.

Hence, all $\ast $-moments of $\Pi $ are the same as $\ast $-moments of $X$,  and $\Pi$ and $X$ coincide in distribution. 
\end{proof}

A key special case occurs if all \(a_{i}\) are identically distributed. Then \(\Pi = a_{n}\cdots a_{1}\) has the same distribution as \(\bigl(U\,\lvert a\rvert \bigr)^{n}\), where \(U\) is free from \(\lvert a\rvert\). From this, one easily shows:

\begin{theo}
\label{theo_Rdiagonal_powers}
If $X, X_1, \ldots, X_n$ are identically distributed $R$-diagonal elements, and $X_1, \ldots, X_n$ are free, then $X^{n}$ is $R$-diagonal and has the same distribution as the product $X_1, \ldots, X_n$.
\end{theo}

\begin{proof}
Let $X = U |X|$. In Lemma \ref{lemma_products} we can take all $A_i = |X|$ and $a_i = X_i$. Then the conclusion of the lemma shows that 
\[
X^{n} \;\cong\; X_{n}\,\cdots\,X_{1}.
\]
A product of free $R$-diagonal variables is again $R$-diagonal, so \(X^{n}\) is $R$-diagonal.
\end{proof}

\section*{Notes}
$R$-diagonal random variables were introduced in \cite{nica_speicher97} who also proved the results in Theorems \ref{Theorem_R_diagonal_2} and \ref{Theorem_R_diagonal_3}. 

For results in Exercises \ref{exe_SA_Rdiagonal_0}, \ref{exe_SA_Rdiagonal_0.5},   \ref{exe_SA_Rdiagonal_0.75}, \ref{exe_product_even}, \ref{exe_SA_Rdiagonal_2} see Proposition 15.8, Theorem 15.10, Corollary 15.11, Theorem 15.17, and Corollary 15.18  in \cite{nica_speicher06}. 

Exercise \ref{exe_R_diagonal_product} correspond to Exercise 15.24 in \cite{nica_speicher06}, which also gives a formula for the non-tracial case. 

A different proof of Proposition \ref{propo_R_diagonal_sum_1} is given in 
\cite{haagerup_larsen00}. It is based on the following interesting result. 

\begin{lemma}[Lemma 3.4 in \cite{haagerup_larsen00}]
Let $a$, $x$, $y$ be free even elements in a non-commutative probability space $(\AA, \phi)$.  Assume that $a^2 = 1$. Then the sets $\{a x, x a\}$ and $\{a y, y a\}$ are free. 
\end{lemma}

Theorem \ref{theo_Rdiagonal_powers} was proved in \cite{haagerup_larsen00} (Proposition 3.10) using the following surprising result.

\begin{lemma}[Lemma 3.7 in \cite{haagerup_larsen00}]
Let $(\AA, \phi)$ be a non-commutative *-probability space, let $u$ be a Haar unitary in $\AA$. Assume that $S$ is a set in $\AA$ such that $S$ and $\{u\}$ are *-free. 

Then for any natural number $n$ we have that 
\begin{enumerate}
\item the sets $S$, $u S u^\ast, u^2 S (u^\ast)^2, \ldots$ are *-free,  
\item the sets $S$, $u S u^\ast, \ldots, u^{n - 1} S (u^\ast)^{n -1}$, $\{u^n\}$ are *-free, 
\item 
the sets $u S u^\ast, \ldots, u^{n} S (u^\ast)^{n}$, $\{u^n\}$ are *-free.
\end{enumerate} 
\end{lemma}

 \cite{haagerup_larsen00} also proved that if an $R$-diagonal element $a$ is invertible, then $a^{-1}$ is also $R$-diagonal and showed how one can calculate its distribution. 

The result of Theorem \ref{theo_Rdiagonal_powers} also appears as Exercise~15.25 in \cite{nica_speicher06}. Our proof of Theorem  \ref{theo_Rdiagonal_powers} is inspired by the proof of Proposition 15.22 in \cite{nica_speicher06} in which they showed that any power $a^r$ of an R-diagonal element $a$ is $R$-diagonal. 


\section{Brown measure of $R$-diagonal variables}

\index{Brown measure}

 It is not straightforward how one should define a generalization
of the eigenvalue distribution for infinite-dimensional, non-normal
operators. One interesting definition is that of the Brown measure.
It is defined only for operators in von Neumann algebras and uses the fact
that in these algebras one can define an analogue of the determinant, which
is called the Fuglede--Kadison determinant.

\index{determinant!Fuglede-Kadison} 
\index{Fuglede-Kadison determinant|textbf}

\begin{defi}
Let $X$ be an element of a tracial $W^{\ast}$-probability
space $\left( \mathcal{A},\tau \right)$. Then the \emph{Fuglede--Kadison determinant}
of $X$ is defined as
\[
  \det X := \exp \biggl[
    \frac{1}{2}\,\tau\bigl(\log \bigl(X^{\ast }X\bigr)\bigr) 
  \biggr] = \exp \biggl(\int_0^\infty \log t \, d\mu_{|X|} (t)\biggr) \in [0, \infty),
\]
where $|X| = (X^\ast X)^{1/2}$. 
\end{defi}
For a random variable $X$ that has non-trivial kernel, $\mu_{|X|}$ puts a nonzero mass on $x = 0$. In this case the integral is $-\infty$ and $\det(X) = 0$. However, $\det(X)$ can be zero even if the kernel is non-zero and $X$ is invertible as an unbounded operator. 

Very often, the FK determinant of $X$ is denoted as $\Delta(X)$.

Here is a list of some of the basic properties of the FK determinant. They are analogous to the properties of the \emph{absolute value of the usual matrix determinant}.  

\begin{theo}
\label{propo_FKdet_properties}
\begin{enumerate}
\item 
$\det(XY) = \det(X) \det(Y)$.
\item 
$\det(X^\ast) = \det(X)$.
\item $\det \bigl( e^{X}\bigr) =\bigl| e^{\tau (X) }\bigr| =\exp
\bigl( \Re \tau( X) \bigr) .$
\item 
If $X > 0$, that is,  $X$ is positive, then $\det(X) \leq \tau(X)$. 
\end{enumerate}

\end{theo}

(A technical assumption for the operators in this proposition is that the operators are closed, densely defined  operators affiliated to $\AA$.)


\begin{exa} Let us consider the algebra of $N\times N$ matrices
with
\[
  \tau(X)=\frac{1}{N}\,\sum_{i=1}^{N} X_{ii}.
\]
Then we can write the Fuglede--Kadison determinant as
$\det X = \bigl(\prod_{i=1}^{N} s_i\bigr)^{1/N}$, where $s_i$ are the
singular values of the matrix $X$ (i.e.\ the square roots of the
eigenvalues of $X^\ast X$).  From linear algebra, we get
\[
  \det X
  =\Bigl[\mathrm{Det}\bigl(X^\ast X\bigr)\Bigr]^{\tfrac{1}{2N}}
  =\bigl\lvert\mathrm{Det}(X)\bigr\rvert^{\tfrac{1}{N}},
\]
where $\mathrm{Det}(X)$ is the usual (classical) determinant. In particular
that
\[
  \log \det\bigl(z - X\bigr)
  =\frac{1}{N}\,\sum_{i=1}^N \log\bigl\lvert z - \lambda_i\bigr\rvert,
\]
where $\lambda_i$ are the eigenvalues of $X$ taken with multiplicities
equal to the number of times that $\lambda_i$ is repeated on the diagonal
of the Jordan form of $X$.
\end{exa}

Note that in this example,
$\log \det\bigl(z - X\bigr)$ is precisely the logarithm of the
absolute value of the characteristic polynomial, divided by the degree of
that polynomial. In a more general situation, we can think of the function
$\log \det\bigl(z - X\bigr)$ as a suitable generalization of the
logarithm of the modulus of a characteristic polynomial.

\begin{exa}
\label{exercise_determinant_1_minus_tS}
Suppose $a$ is a random variable with distribution
$\{\delta_{-1} + \delta_{1}\} / 2$, and let $z$ be a complex number such
that $\lvert z \rvert < 1$. 
\[
  \log \det(1 - z a) = \tfrac12 \log | 1 - z^2 | =  
\tfrac14 \log \bigl(1 -  2\,\mathrm{Re}(z^2) + \lvert z \rvert^4 \bigr).
\]
Indeed, $1 - z a$ takes values $1 - z$ or $1 + z$ with
probability $\tfrac12$.  Consequently, $\bigl(1 - z a\bigr)^\ast (1 -  z a)$ has the measure 
that put equal weights $\tfrac12$ on  $\lvert 1 - z \rvert^2$ and $\lvert 1 + z \rvert^2$.  Then 
\begin{align*}
\tfrac12\, \tau \bigl[\log\bigl((1 - z a)^\ast (1 - z a)\bigr)\bigr] &= \tfrac14 \log\big(|1 + z|^2 |1 - z|^2\big) 
= \tfrac12 \log| 1 - z^2|.
\end{align*}
 And the second equality is an identity.
\end{exa}

\index{L-function|textbf}
\begin{defi}
The \emph{L-function} of a variable $X$ is defined as
\[
  L_X(z)
  :=\log \det\bigl(z - X\bigr)
  = \tau\bigl[\log \lvert z - X \rvert\bigr],
\]
where
$\lvert z - X \rvert =\bigl[\,(z - X)^\ast( z - X)\bigr]^{1/2}$.
\end{defi}

Since the infinite-dimensional determinant is well-defined for all bounded
operators in a tracial $W^{\ast}$-algebra (and can, with some effort,
be extended to unbounded operators), we can \emph{exploit} this fact to
define the spectral measure of these operators. Indeed, in finite dimensions
the determinant is the product of eigenvalues. Hence, if we know
$\det(z - X)$ for all $z$, then we know the characteristic
polynomial and can recover the eigenvalues as its zeros. In the
infinite-dimensional case, the situation is more subtle: the function
$L_X(z) = \log \det\bigl(z - X \bigr)$ can fail to be
harmonic on an open set.  Interpreting $L_X(z)$ as a potential, this
failure means there is a continuous “distribution of charges” on that open
set. This suggests the following definition.

\index{Brown measure|textbf}
\begin{defi}
Let $X$ be a bounded random variable in a tracial $W^{\ast}$-probability
space $( \mathcal{A},\tau )$. Then its \emph{Brown measure} is a
measure $\mu_X$ on the complex plane $\mathbb{C}$ defined by
\[
  \mu_X
  =\frac{1}{2\pi}\,\Delta\,L_X(z)\,\mathrm{d}x\,\mathrm{d}y,
\]
where $z = x + i\,y$, $\Delta = \partial_x^2 + \partial_y^2$ is
the Laplace operator, and the equality holds in the sense of
(Schwartz) distributions.
\end{defi}

We will not prove here that this construction gives a well-defined
probability measure, and instead refer to the original paper by Brown. We
only mention that the intuitive reason for this is that the function
$L_X(z)$ is subharmonic, and that the Laplacian of $L_X(z)$
measures the extent to which $L_X(z)$ fails to be harmonic.

To gain intuition for how this definition works, consider once again the
finite-dimensional case of $N \times N$ matrices, where
\[
  L_X(z)
  =\frac{1}{N}\,\sum_{i=1}^{N}\,\log\bigl\lvert z - \lambda_i \bigr\rvert.
\]
Then $L_X(z)$ is harmonic everywhere outside the points
$\lambda_i$. Let $f$ be a smooth function which is constant in some
neighborhood of $\lambda_1$ and zero at the other $\lambda_i$.  Let $C$
be a circle of radius $\varepsilon$ around $\lambda_1$:
\[
  C =\bigl\{\, z : \lvert z - \lambda_1\rvert \le \varepsilon \bigr\}.
\]
By one of Green’s formulas,
\[
  \frac{1}{2\pi}\,\iint_C f\,\Delta L_X\,\mathrm{d}x\,\mathrm{d}y
  =\frac{1}{2\pi}\,\oint_{\partial C} f\,
    \bigl(\nabla L_X \cdot n\bigr)\,\mathrm{d}s,
\]
where $n$ is the outward normal vector to the contour $\partial C$. We
can compute $\nabla L_X \cdot n = (N\,\varepsilon)^{-1}$. Hence, the
integral on the right equals $f(\lambda_1)/N$, and we conclude that
\[
  \mu_X =\frac{1}{N}\,\sum_{k=1}^N \delta_{\lambda_k},
\]
where $\delta_x$ \emph{denotes} the Dirac distribution concentrated at the
point $x$.

We now list (without proof) some of the properties of the Brown measure:

\smallskip
\noindent
\textbf{(i)}~It is the unique measure such that
\[
  L_X(z)
  =\int_{\mathbb{C}} \log\bigl\lvert z - \lambda\bigr\rvert \,
  \mathrm{d}\mu_X(\lambda).
\]

\smallskip
\noindent
\textbf{(ii)}~For every integer $n \ge 0$, we have
\[
  \tau \bigl(X^n\bigr)
  =\int z^n\,\mathrm{d}\mu_X(z).
\]

\smallskip
\noindent
\textbf{(iii)}~The Brown measure of a normal operator $X$ coincides with
its usual spectral probability distribution.

\bigskip

For an $R$-diagonal operator, the Brown measure is invariant under
rotations around the origin in the complex plane.  We can thus write it as a
product of its radial and polar parts.  Let $\rho_X$ denote the
\emph{radial part} of the Brown measure.  Concretely, if $A$ is a thin
annulus between circles with \emph{radii} $r$ and $r + \mathrm{d}r$,
then the measure of this annulus is
$\rho_X\bigl([r,\,r+\mathrm{d}r]\bigr)$.  By a slight abuse of notation we
use the same symbol $\rho_X$ for this distribution function and its
density.

It is useful to have a formula for $\rho_X$ in terms of the $L$-function
$L_X(z)$.  Note that if $f(r,\varphi)$ depends only on $r$, then
$\Delta f = f_{rr} + \tfrac{1}{r}f_r$.

\begin{lemma}
\label{lemma_radial_part_computation}
Let $z = r\,e^{i\varphi}$ and suppose $L_X(z)$ depends only
on $r$, i.e.\ $L_X(z) = L(r)$. Then the radial part of the Brown
measure density can be computed as
\[
  \rho_X(r)
  =\bigl(r\,L'(r)\bigr)'.
\]
\end{lemma}

\begin{proof}
Consider the circular annulus $A$ with inner radius $R$ and outer radius
$R + \Delta R$.  Then
\[
  \mu_X(A)
  =\int_0^{2\pi}\int_{R}^{R+\Delta R}
      \frac{1}{2\pi}\,\bigl(L_{rr} + \tfrac{1}{r}L_r\bigr)\,r\,
      \mathrm{d}r\,\mathrm{d}\varphi
  =\int_{R}^{R+\Delta R}
      \bigl(r\,L_{rr} + L_{r}\bigr)\,\mathrm{d}r.
\]
Since $\mu_X(A)$ should approximate $\rho_X(R)\,\Delta R$ for small
$\Delta R$, it follows that
$\rho_X(r) = r\,L_{rr}(r) + L_{r}(r)$, which can be rewritten as
$\rho_X(r) = \bigl(r\,L'(r)\bigr)'$.
\end{proof}

\begin{coro}
\label{corollary_radial_part_computation}
Let $F(r)$ denote the distribution function corresponding to the polar
measure $\rho_X$. Then
\[
  F(R) - F(0)
  = R\,L'_X(R).
\]
\end{coro}

Another natural distribution associated to an $R$-diagonal operator $X$ is that of the squared singular values, i.e. the probability distribution of $X^*X$.  We denote this measure by $\sigma_X$, with density (when it exists) also written $\sigma_X$.

Typically one writes $X=U\,H$, where $U$ is Haar-unitary and $H\ge0$.  Then
\[
  X^*X \;=\; H^2,
\]
so $\sigma_X$ is simply the spectral distribution of $H^2$.

\begin{figure}[tbph]
  \centering
  \includegraphics[width= \textwidth]{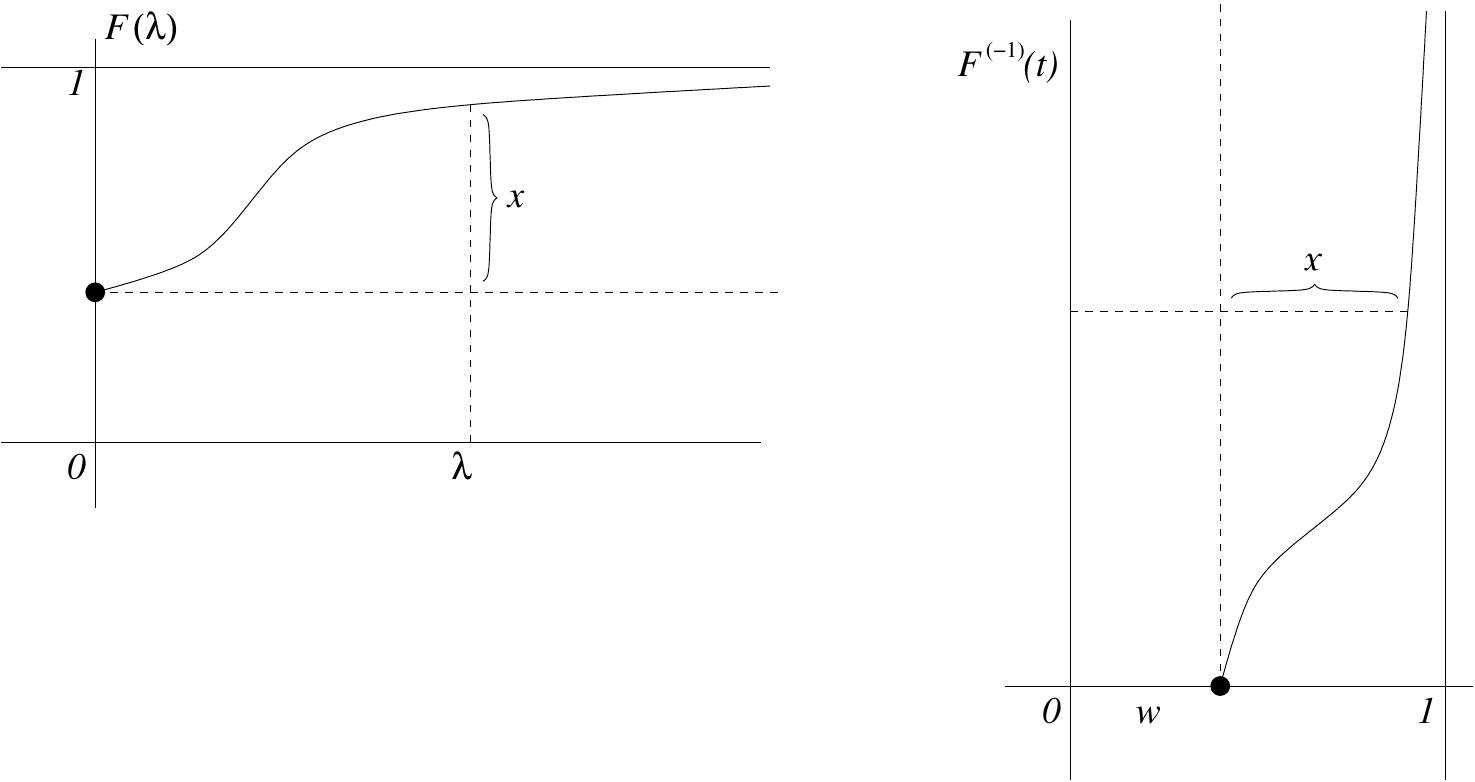}
  \caption{Distribution function of the radial measure and its inverse}
  \label{figure_Haagerup_Larsen1}
\end{figure}

We can then ask how the distributions $\sigma_X$ and the radial law $\rho_X$ relate.  The answer is given by a  theorem in \cite{haagerup_larsen00}. Before formulating it, let us introduce some notation. Write 
\begin{itemize}
 \item $F(r) \equiv F_X(r)$, $r \geq 0$, for the distribution function of the radial measure $\rho_X$, 
\item $F^{-1}(t)$ for its functional inverse on the range of $F$, and
\item $S(z)$ for the $S$-transform of $\sigma_X$.  
\end{itemize}

\index{Haagerup--Larsen theorem|textbf}
\begin{theo}[\cite{haagerup_larsen00}]
\label{theo_Haagerup_Larsen}
Let $X$ be $R$-diagonal.  

Set 
\[
  w \;=\;\sigma_X(\{0\}) \;=\;\rho_X(\{0\}).
\]
Then for all $x\in (0,1-w]$,
\[
  F^{-1}(w + x) \;=\;\sqrt{\frac{1}{\,S(x-1)\,}}.
\]
Moreover, the upper edge of $\supp(\rho_X)$ is $\sqrt{\mathbb{E}[X^*X]}$.
\end{theo}

\medskip
\noindent\textbf{Example.}
Let $X=U\,S$ with $S$ semicircular.  Then $X^*X=S^2$, so $\sigma_X$ is free Poisson with parameter $\lambda=1$  Its $S$-transform is
\[
  S(z)=\frac1{1+z},
\]
and $w = 0$. Then,  
\[
  F^{-1}(t)=\sqrt{t}, 
  \quad
  F(r)=r^2,\quad r\in[0,1].
\]
Hence
\[
  \rho_X(r)=F'(r)=2r,
\]
and the Brown measure is uniform on the unit disk:
\[
  \mu_X(z)
  =\frac1\pi\,r\,dr\,d\theta
  =\frac1\pi\,dx\,dy,
\quad |z|\le1,
\]
and zero outside.

\medskip
Finally, for comparison:  if $X$ is an $N\times N$ Ginibre matrix with i.i.d.\ complex Gaussian entries of variance $(2N)^{-1}$, then its empirical eigenvalue law converges to the uniform law on the unit disk (Ginibre 1965), and this universality extends well beyond the Gaussian case (Tao–Vu).

Before giving a proof of Theorem \ref{theo_Haagerup_Larsen}, let us collect some facts that will be useful in the proof. 

 Let us define $g_X(t): = \tau (I - i t X)^{-1}$. 
\begin{exe}
\label{exercise_g_S_and_A}%
For a r.v. $S$ that have distribution $\mu_S = \tfrac12 (\delta_1 + \delta_{-1})$ we can calculate:
\[
g_S(t) = \frac{1}{1 + t^2}.
\]
More generally, if $A$ is a self-adjoint variable with symmetric distribution $\mu_A$, then 
\[
g_A(t) =
\int_{\mathbb{R}}\frac{1}{1 - i\,t\,x}\,d\mu_A(x)
=
\int_{\mathbb{R}}\frac{1}{1 + t^2 x^2}\,d\mu_A(x).
\]
\end{exe}

It will be convenient to use the notation 
\[
L(X) :=L_{X}( 0) =\log \det (X)
\]

\begin{lemma}
\label{lemma_det_1minusA}
Let $a$ be a random variable in a $W^\ast$-probability space $(\mathcal A, \tau)$.  If 
\[
\tau(a^k)=0\quad\text{for all integers }k>0,
\]
and the spectral radius of $a$ satisfies $\rho(a)<1$, then
\[
L(1-a)=0.
\]
\end{lemma}

This conclusion remains valid even when $\rho(a)=1$; the more general proof may be found in \cite{haagerup_larsen00}.

\begin{proof}
Since \(\rho(a)<1\), the logarithm
\[
x \;=\;\log(1-a)
\;=\;\sum_{k=1}^\infty \frac{a^k}{k}
\]
converges in norm.  Taking expectation term by term gives
\[
\tau(x)
\;=\;\sum_{k=1}^\infty \frac{\tau(a^k)}{k}
\;=\;0.
\]
On the other hand, by definition of the Fuglede–Kadison determinant and property (3) of Theorem \ref{propo_FKdet_properties}, 
\[
L(1-a) = L\bigl(e^x\bigr)
= \log \exp \Re\tau(x) = \log 1
\;=\;0.
\]
\qedhere
\end{proof}

\begin{propo}
\label{proposition_HL_key identity}
Let $A$ be a self-adjoint r.v. with symmetric distribution $\mu_A$ and $S$ is the r.v. with distribution $\tfrac12(\delta_1 + \delta_{-1})$. Suppose that $A$ and $S$ are free and 
$z $ and $t$ are related by the following equation $1+z^2 t^2 =g_A(t) ^{-1}.$ Then, 
\begin{align}
L( A+z S) =& \log (z) -\frac{1}{2}\log \bigl[ 1+( z t) ^2\bigr] +L( 1-itA) . \notag \\
=&\log (z) -\frac{1}{2}\log \bigl[ 1+( z t) ^2\bigr]  \notag \\
&+\frac{1}{2}\int \log \bigl[ 1+x^{2}t^{2}\bigr] d\mu _A(x) .
\end{align}
\end{propo}

\begin{proof}

We intend to write $A+z S$ as a product. We start
with the following expression: 
\begin{equation}
\left( 1-it_{1}A\right) \left( 1-\frac{\left( 1-it_{1}A\right) ^{-1}}{%
g_{A}\left( t_{1}\right) }\right) \left( 1-\frac{\left( 1-it_{2}S\right)
^{-1}}{g_{S}\left( t_{2}\right) }\right) \left( 1-it_{2}S\right) ,
\label{expression_in_HL_theorem}
\end{equation}%
where $g_{A}\left( t_{1}\right) =\tau\bigl[ \left( 1-it_{1}A\right) ^{-1}\bigr]
$ and $g_{S}\left( t_{2}\right) =\tau\bigl[ \left( 1-it_{2}S\right) ^{-1}\bigr]
.$

The main idea of the proof of Proposition \ref{proposition_HL_key identity}
is to make expression (\ref{expression_in_HL_theorem}) close to $\left(
1-it_{1}A\right) \left( 1-it_{2}S\right) +A+z S$ by a suitable choice
of $t_{1}$ and $t_{2}.$ Since 
\begin{equation*}
\tau\left( 1-\frac{\left( 1-it_{1}A\right) ^{-1}}{g_{A}\left( t_{1}\right) }%
\right) =\tau\left( 1-\frac{\left( 1-it_{2}S\right) ^{-1}}{g_{S}\left(
t_{2}\right) }\right) =0,
\end{equation*}%
and $A$ is free from $S,$ hence 
\begin{equation*}
\tau \left( \left( 1-\frac{\left( 1-it_{1}A\right) ^{-1}}{g_{A}\left(
t_{1}\right) }\right) \left( 1-\frac{\left( 1-it_{2}S\right) ^{-1}}{%
g_{S}\left( t_{2}\right) }\right) \right) ^{k}=0,
\end{equation*}%
and we will be able to use Lemma \ref{lemma_det_1minusA} in our calculations
provided that the spectral radius of the operator is smaller or equal than $%
1.$

We can re-write expression (\ref{expression_in_HL_theorem}) as

\begin{eqnarray*}
&&\left( 1-it_{1}A\right) \left( 1-it_{2}S\right) -\frac{1-it_{2}S}{%
g_{A}\left( t_{1}\right) }-\frac{1-it_{1}A}{g_{S}\left( t_{2}\right) }+\frac{%
1}{g_{A}\left( t_{1}\right) g_{S}\left( t_{2}\right) } \\
&=&\left( 1-it_{1}A\right) \left( 1-it_{2}S\right) +\frac{it_{1}}{%
g_{S}\left( t_{2}\right) }\left( A+\frac{t_{2}g_{S}\left( t_{2}\right) }{%
t_{1}g_{A}\left( t_{1}\right) }S\right) +f\left( t_{1},t_{2}\right) ,
\end{eqnarray*}%
where 
\begin{equation*}
f\left( t_{1},t_{2}\right) =\frac{1}{g_{A}\left( t_{1}\right) g_{S}\left(
t_{2}\right) }-\frac{1}{g_{A}\left( t_{1}\right) }-\frac{1}{g_{S}\left(
t_{2}\right) }.
\end{equation*}

We impose on $t_{1}$ and $t_{2}$ the requirement that 
\begin{equation}
\frac{t_{2}g_{S}\left( t_{2}\right) }{t_{1}g_{A}\left( t_{1}\right) }%
=z ,  \label{requirement1}
\end{equation}%
and 
\begin{equation}
f\left( t_{1},t_{2}\right) =0.  \label{requirement2}
\end{equation}
It turns out that it is possible to satisfy the two previous equations by a
suitable choice of $t_{1}$ and $t_{2}$ if $z $ is not too large. We will address the choice of $t_1$ and $t_2$ later. Then 
\begin{eqnarray*}
&&( 1- i t_1 A) \Biggl[ 1-\Bigl( 1-\frac{( 1- i t_1 A) ^{-1}}{g_A( t_1) }\Bigr) \Bigl( 1-\frac{( 1- i t_2S)^{-1}}{g_S( t_2) }\Bigr) \Biggr] ( 1-i t_{2}S) \\
&=&-\frac{t_{1}}{g_{S}\left( t_{2}\right) }\left( A+z S\right) = h\left( z \right) \left( A+z S\right) ,
\end{eqnarray*}%
where 
\begin{equation*}
h( z ) :=-\frac{it_{1}}{g_{S}( t_2) }= - i t_1 \bigl(1 + (t_2)^2\bigr)
\end{equation*}
by Exercise \ref{exercise_g_S_and_A}.

An application of Lemma \ref{lemma_det_1minusA} shows that 
\begin{equation*}
\det \left[ 1-\left( 1-\frac{\left( 1- i t_{1}A\right) ^{-1}}{g_{A}\left(
t_{1}\right) }\right) \left( 1-\frac{\left( 1- i t_{2}S\right) ^{-1}}{%
g_{S}\left( t_{2}\right) }\right) \right] =1,
\end{equation*}
provided that the spectral radius of 
\begin{equation}
\label{equ_operator_ha_ha}
\left( 1-\frac{\left( 1-it_{1}A\right) ^{-1}}{g_{A}\left( t_{1}\right) }%
\right) \left( 1-\frac{\left( 1-it_{2}S\right) ^{-1}}{g_{S}\left(
t_{2}\right) }\right)
\end{equation}%
is less than or equal to 1.

Assume for the moment that this fact is established. Then, we see that%
\begin{align}
\det \left( A+z S\right) &=\frac{1}{\left| h\left( z \right)
\right| }\det \left( 1-it_{1}A\right) \det \left( 1-it_{2}S\right)
 \notag \\
&=\frac{1}{t_{1}\sqrt{1+\left( t_{2}\right) ^{2}}}\det \left(
1-it_{1}A\right) \label{formula_det_A_plus_lambda_S},
\end{align}%
where we used definition of $h(z )$ and Exercise \ref{exercise_determinant_1_minus_tS} to write $\det \bigl( 1-it_2 S\bigr) = \sqrt{1 + (t_2)^2}$.

The next step is to figure out the dependence of $t_{1}$ and $t_{2}$ on $%
z .$ 
\begin{exe}
Equations (\ref{requirement1}) and (\ref{requirement2}) can be solved for $g_A(t_1)$ and $t_2$ as follows:
\begin{equation}
\label{requirement3}
g_{A}( t_1) =\frac{1}{1+( z t_1) ^2} \text{ and } t_2 = \frac{1}{z t_1} 
\end{equation}
\end{exe}
That is, if we choose $t_1$ and $t_2$ in such a way that \eqref{requirement3} is satisfied, then we ensure that (\ref{requirement1}) and (\ref{requirement2}) hold true. 


Then, \eqref{formula_det_A_plus_lambda_S} implies that 
\begin{equation*}
\det \left( A+z S\right) =\frac{z }{\sqrt{1+\left( z
t\right) ^{2}}}\det \left( 1-itA\right) ,
\end{equation*}%
where $z $ and $t$ are related by the expression: 
\begin{equation*}
g_{A}\left( t\right) :=\int \frac{1}{1+t^{2}x^{2}}d\mu _{A}\left( x\right) =%
\frac{1}{1+\left( z t\right) ^{2}}.
\end{equation*}%
We can also write this as 
\begin{equation*}
\det \left( A+z S\right) =z \sqrt{g_{A}\left( t\right) }\det
\left( 1-itA\right) ,
\end{equation*}

If we take logarithm we get 
\begin{equation*}
L\left( A+z S\right) =\log \left( z \right) -\frac{1}{2}\log %
\left[ 1+\left( z t\right) ^{2}\right] +L\left( 1-itA\right) .
\end{equation*}%
This completes proof of Proposition \ref{proposition_HL_key identity} subject to the proof that the operator in \eqref{equ_operator_ha_ha} has the spectral radius $\leq 1$. 

\end{proof}

The proof that the operator  in \eqref{equ_operator_ha_ha} has the spectral radius $\leq 1$ is based on 
\begin{propo}
\label{spectral_radius_product}
Let $a$ and $b$ be two bounded random variables in $W^{\ast }$-probability $(\AA, \tau)$
space with the faithful trace $\tau$. Assume that $a$ and $b$ are $\ast$-free and centered:  $\tau(a) =\tau( b) =0.$ Then the
following formula holds for spectral radius of the product: 
\begin{equation*}
\rho ( ab) =\| a \|_2 \| b\|_2,
\end{equation*}
where $\|x \|_2 := \sqrt{\tau (x^\ast x)}$
\end{propo}

\begin{proof} Assume without loss of generality that $\|a\| = \|b\| = 1$. We need to show that $\rho(a b) = 1$. 
Let $\MM_a$ and $\MM_b$ be subalgebras of von Neumann probability space $(\MM, \tau)$ generated by $a$ and $b$, respectively. Let $(\HH_a, \xi_a)$ and $(\HH_b, \xi_b)$ be Hilbert spaces corresponding by the GNS construction to $\MM_a$ and $\MM_b$ and let $(\HH, \xi) = (\HH_a, \xi_a) \ast (\HH_b, \xi_b)$. Then, 
\[
\HH = \bC \xi \oplus \bigoplus_{\substack{n \in \N \\ j_1 \ne \ldots \ne j_n}} \HH_{j_1}^o \otimes \ldots \otimes \HH_{j_n}^o, 
\]
where $\HH_{j}^o = \{\xi_j\}^\perp \subset \HH_j$. 

Note that 
\[
\HH_{j_1}^o \otimes \ldots \otimes \HH_{j_n}^o = [\MM_{j_1}^o \otimes \ldots \otimes \MM_{j_n}^o \xi], 
\quad j_1 \ne \ldots \ne j_n,
\]
where $\MM_j^o$ denote the set of centered elements of $\MM_j$, $j = a, b$, and $[S]$ denotes the closed linear span of a set $S$. Put
\begin{align*}
\KK_0 &= \bC \xi, \\
\KK_n &= [\underbrace{\MM_a^o \MM_b^0 \ldots}_{n} \xi], \quad \LL_n = [\underbrace{\MM_b^o \MM_a^o \ldots}_{n} \xi], \quad n \in \N, \\
\KK &= \bigoplus_{n = 0}^\infty \KK_n, \quad \LL = \bigoplus_{n = 1}^\infty \LL_n, 
\end{align*}
Then $\HH = \KK \oplus \LL$, and $a b \KK_n \subset \KK_{n + 2}$, $n = 0, 1, \ldots$. Hence the matrix of $ab$ in the decomposition $\KK \oplus \LL$ is 
\[
a b = \bmatr{ R & S \\ 0 & T}
\]
Note that $R(a_1 b_1  \ldots a_n b_n) = a b a_1 b_1 \ldots a_n b_n$ so $R$ corresponds to tensoring from the left by $a \xi_a \otimes b \xi_b$. Since $\|a \|_2 = \| b\|_2 = 1$, this is an isometric map of $\KK_n$ into $\KK_{n + 2}$. Hence, $R$ maps $\KK$ isometrically into $\KK$. It follows that $\| R^p \| = 1$ for all $p \in \N$. 

Next, we have 
\[
(a b)^\ast = \bmatr{ R^\ast & 0 \\ S^\ast & T^\ast},
\]
so $\LL$ is invariant under $T^\ast = b^\ast a^\ast$. Then, using $\|a^\ast\|_2 = \|a\|_2 = 1$ and $\|b^\ast\|_2 = \|b\|_2 = 1$, we find that $T^\ast$ is an isometry of $\LL$ into $\LL$ and $\|T^p\| = \|(T^\ast)^p\| = 1$ for all $p$. Then, 
\[
(a b)^p = \bmatr{R^p & 0 \\ 0 & T^p} + \sum_{r = 0}^{p - 1} \bmatr{0 & R^{p - r - 1} S T^r \\ 0 & 0}.
\]
Hence, 
\[
1 \leq \|(a b)^p\| \leq 1 + p \|S\| \leq 1 + p \| a b\|, 
\]
and therefore $\rho(ab) = \lim_{p \to \infty} \|(a b)^p\|^{1/p} = 1$ as desired.

\end{proof}

\section*{Exercises}
\begin{exe}
If $s$ is a standard semicircular, then $\det (s) = e^{-1/2}$. 
\end{exe}

\section*{Notes}
The Fuglede--Kadison determinant was defined in \cite{fuglede_kadison52} which also developed its main properties, including properties in Prop. \ref{propo_FKdet_properties}. They worked mostly with bounded operators with a bounded inverse (which they call regular operators). For extension of these properties to a more general class of operators, including some unbounded operators, see \cite{haagerup_schultz2007}. Property 3. in Prop. \ref{propo_FKdet_properties} is Lemma 3.4 in \cite{nayak2018}.
This is Proposition 4.1 in \cite{haagerup_larsen00}.

In the context of random matrices, an analogue of the Haagerup-Larsen theorem was proved in \cite{guionnet_krishnapur_zeitouni11}.

\appendix{}

\chapter{Distributions of Self-Adjoint Variables and Their Transforms}

Note: The Cauchy transform is defined as follows:%
\begin{equation*}
G\left( z\right) :=E\left( \frac{1}{z-X}\right) ;
\end{equation*}%
the K-function is the inverse of the Cauchy transform,%
\begin{equation*}
K\left( z\right) :=G^{\left( -1\right) }\left( z\right) ;
\end{equation*}%
the S-function is 
\begin{equation*}
S\left( z\right) :=\frac{1+z}{z}\psi ^{\left( -1\right) }\left( z\right) ,
\end{equation*}%
where $\psi ^{\left( -1\right) }\left( z\right) $ is the inverse of the
moment-generating function 
\begin{equation*}
\psi \left( z\right) :=E\left( \frac{zX}{1-zX}\right) .
\end{equation*}%
If $R\left( z\right) :=zK\left( z\right) -1$ (this is a variant of the
definition of the $R$-transform)$,$ then $S\left( z\right) $ can also be
defined by 
\begin{equation*}
S\left( z\right) =\frac{1}{z}R^{\left( -1\right) }\left( z\right) .
\end{equation*}

\begin{tabular}{ll}
& \textbf{Semicircle} \\ 
Density \& atoms & $\frac{1}{2\pi }\sqrt{4-x^{2}}\chi _{\lbrack -2,2]}\left(
x\right) $ \\ 
Cauchy transform & $\frac{1}{2}\left( z-\sqrt{z^{2}-4}\right) $ \\ 
K-function & $\frac{1}{z}+z$ \\ 
S-function & $\pm \frac{1}{\sqrt{z}}$ \\ 
Moments & $m_{2k}=\frac{1}{k+1}\binom{2k}{k},$ $m_{2k+1}=0.$ \\ 
Free Cumulants & $c_{1}=1,$ $c_{i}=0$ for $i\neq 1$%
\end{tabular}

\bigskip
\footnote{Free Poisson distribution differs from the Marchenko--Pastur distribution by a scaling. So you can find different formulas for the density and the Cauchy transform of the Marchenko--Pastur distribution in other sources.}
\begin{tabular}{ll}
& \textbf{Free Poisson} \\ 
Density \& atoms & $\frac{\sqrt{4x-\left( 1-\lambda +x\right) ^{2}}}{2\pi x}%
\chi _{\left[ \left( 1-\sqrt{\lambda }\right) ^{2},\left( 1+\sqrt{\lambda }%
\right) ^{2}\right] }\left( x\right) $ \\ 
& and an atom at 0 with mass $\left( 1-\lambda \right) $ if $\lambda <1$ \\ 
Cauchy transform & $\frac{1-\lambda +z-\sqrt{\left( 1-\lambda +z\right)
^{2}-4z}}{2z}$ \\ 
K-function & $\frac{1}{z}+\frac{\lambda }{1-z}$ \\ 
S-function & $\frac{1}{\lambda +z}$ \\ 
Moments &  \\ 
Free Cumulants & $c_{i}=\lambda $ for all $i.$%
\end{tabular}

\bigskip

\begin{tabular}{ll}
& \textbf{Bernoulli }I \\ 
Density \& atoms & $p\delta _{1}+q\delta _{0}$ \\ 
Cauchy transform & $\frac{z-q}{\left( z-1\right) z}$ \\ 
K-function & $\frac{1+z-\sqrt{\left( 1+z\right) ^{2}-4qz}}{2z}$ \\ 
S-function & $\frac{1+z}{p+z}$ \\ 
Moments & $m_{i}=p$ for all $i$ \\ 
Free Cumulants & 
\end{tabular}

\bigskip

\begin{tabular}{ll}
& \textbf{Bernoulli }II \\ 
Density \& atoms & $\frac{1}{2}\left( \delta _{1}+\delta _{-1}\right) $ \\ 
Cauchy transform & $\frac{z}{z^{2}-1}$ \\ 
K-function & $\frac{1+\sqrt{1+4z^{2}}}{2z}$ \\ 
S-function & $\pm \sqrt{1+\frac{1}{u}}$ \\ 
Moments & $m_{2k}=1$ and $m_{2k+1}=0$ \\ 
Free Cumulants & $c_{2k}=\left( -1\right) ^{k-1}2^{k}\frac{\left(
2k-3\right) !!}{k!}$%
\end{tabular}

\bigskip

\begin{tabular}{ll}
& \textbf{Arcsine} \\ 
Density \& atoms & $\frac{\chi _{\left[ -2,2\right] }\left( t\right) }{\pi 
\sqrt{4-t^{2}}}$ \\ 
Cauchy transform & $\frac{1}{\sqrt{z^{2}-4}}$ \\ 
K-function & $\frac{1}{z}\sqrt{1+4z^{2}}$ \\ 
S-function & $\frac{\pm \sqrt{u(2+u)}}{2\left( 1+u\right) }$ \\ 
Moments & $m_{n}=\binom{2k}{k}$ if $n=2k;$ $=0$ if $n=2k+1$ \\ 
Free Cumulants & 
\end{tabular}

\bigskip

\begin{tabular}{ll}
& \textbf{Cauchy} \\ 
Density \& atoms & $\frac{1}{\pi }\frac{1}{1+t^{2}}$ \\ 
Cauchy transform & $\frac{1}{z+i}$ \\ 
K-function & $\frac{1}{z}-i$ \\ 
S-function & $i$ \\ 
Moments & $E\left| X\right| =\infty $ \\ 
Free Cumulants & 
\end{tabular}


\chapter{Scaling Properties}

It is useful to know how various functions related to the Cauchy transform
of random variable $X$ behave under scaling and translation. For convenience
of reference we collect these results in the following propositions:

\index{scaling properties|textbf}

\begin{propo}
\label{proposition_scaling_properties} \ \newline
(i) $G_{aX}\left( z\right) =a^{-1}G_{X}\left( z/a\right) ;$\newline
(ii) $R_{a X}\left( z\right) =a R_{X}\left( a z\right) ;$%
\newline
(iii) $S_{aX}\left( z\right) =a^{-1}S_{X}\left( z\right) $
\end{propo}
\begin{proof}
(i) is obvious from the definition of the Cauchy transform, (ii) was proved in (\ref{scaling_R}), and (iii) was proved in Proposition \ref{propo_S_scaling}.
\end{proof}

\begin{propo}
\label{proposition_translation_properties} \ \newline
(i) $G_{X+b}\left( z\right) =G_{X}\left( z-b\right) ;$\newline
(ii) $R_{X+b}\left( z\right) =R_{X}\left( z\right) +b;$\newline
\end{propo}

These properties directly follow from definitions.

\bibliographystyle{plainnat}
\bibliography{comtest}

\begin{thebibliography}{46}
\providecommand{\natexlab}[1]{#1}
\providecommand{\url}[1]{\texttt{#1}}
\expandafter\ifx\csname urlstyle\endcsname\relax
  \providecommand{\doi}[1]{doi: #1}\else
  \providecommand{\doi}{doi: \begingroup \urlstyle{rm}\Url}\fi

\bibitem[Aigner(1979)]{aigner1979}
Martin Aigner.
\newblock \emph{Combinatorial Theory}.
\newblock Springer Verlag, 1979.

\bibitem[Belinschi and Bercovici(2005)]{belinschi_bercovici05}
S.~T. Belinschi and H.~Bercovici.
\newblock Partially defined semigroups relative to multiplicative free
  convolution.
\newblock \emph{International Mathematics Research Notices}, 2005:\penalty0
  65--101, 2005.

\bibitem[Belinschi and Bercovici(2007)]{belinschi_bercovici07}
S.~T. Belinschi and H.~Bercovici.
\newblock A new approach to subordination results in free probability.
\newblock \emph{Journal d'Analyse Math{\'e}matique}, 101:\penalty0 357--365,
  2007.

\bibitem[Belinschi et~al.(2017)Belinschi, Mai, and Speicher]{bms2017}
Serban Belinschi, Tobias Mai, and Roland Speicher.
\newblock Analytic subordination theory of operator-valued free additive
  convolution and the solution of a general random matrix problem.
\newblock \emph{Journal f{\"u}r die reine und angewandte Mathematik (Crelles
  Journal)}, \penalty0 (732):\penalty0 21--53, 2017.
\newblock \doi{10.1515/crelle-2014-0138}.

\bibitem[Belinschi et~al.(2014)Belinschi, Speicher, Treilhard, and
  Vargas]{bstv2014}
Serban~T. Belinschi, Roland Speicher, John Treilhard, and Carlos Vargas.
\newblock Operator-valued free multiplicative convolution: Analytic
  subordination theory and applications to random matrix theory.
\newblock \emph{International Mathematics Research Notices}, pages 5933--5958,
  2014.

\bibitem[Bercovici and Voiculescu(1993)]{bercovici_voiculescu93}
H.~Bercovici and D.~Voiculescu.
\newblock Free convolution of measures with unbounded support.
\newblock \emph{Indiana University Mathematics Journal}, 42:\penalty0 733--773,
  1993.

\bibitem[Bercovici and Voiculescu(1995)]{bercovici_voiculescu95}
H.~Bercovici and D.~Voiculescu.
\newblock Superconvergence to the central limit and failure of the
  \mbox{C}ramer theorem for free random variables.
\newblock \emph{Probability Theory and Related Fields}, 102:\penalty0 215--222,
  1995.

\bibitem[Bercovici et~al.(1999)Bercovici, Pata, and Biane]{bercovici_pata99}
Hari Bercovici, Vittorino Pata, and Philippe Biane.
\newblock Stable laws and domains of attraction in free probability theory.
\newblock \emph{Annals of Mathematics}, 149:\penalty0 1023--1060, 1999.
\newblock with an appendix by \mbox{P}hilippe \mbox{B}iane.

\bibitem[Biane(1998)]{biane98b}
Philippe Biane.
\newblock Processes with free increments.
\newblock \emph{Mathematische Zeitschrift}, 227:\penalty0 143--174, 1998.

\bibitem[Brillinger(1980)]{brillinger80}
David~R. Brillinger.
\newblock \emph{Time Series. Data Analysis and Theory}.
\newblock Holden-Day Series in Statistics. Holden-Day, Inc., expanded edition,
  1980.

\bibitem[Collins(2002)]{collins02}
Beno\^it Collins.
\newblock Moments and cumulants of polynomial random variables on unitary
  groups, the \mbox{I}tzykson-\mbox{Z}uber integral and free probability.
\newblock arxiv:math-ph/0205010v2, 2002.

\bibitem[Collins(2003)]{collins2003}
Beno\^it Collins.
\newblock Moments and cumulants of polynomial random variables on unitary
  groups, the \mbox{I}tzykson-\mbox{Z}uber integral and free probability.
\newblock \emph{International Mathematics Research Notices}, pages 953--982,
  2003.

\bibitem[Collins and \'Sniady(2006)]{collins_sniady2006}
Beno\^it Collins and Piotr \'Sniady.
\newblock Integration with respect to the \mbox{H}aar measur on unitary,
  orthogonal and symplectic group.
\newblock \emph{Comm. Math. Phys.}, 264:\penalty0 773--795, 2006.

\bibitem[Collins et~al.(2022)Collins, Matsumoto, and Novak]{cmn2022}
Beno\^it Collins, Sho Matsumoto, and Jonathan Novak.
\newblock The weingarten caclulus.
\newblock \emph{Notices of AMS}, 69:\penalty0 734--745, 2022.

\bibitem[Far et~al.(2006)Far, Oraby, Bryc, and Speicher]{fobs2006}
R.~Rashidi Far, T.~Oraby, W.~Bryc, and R.~Speicher.
\newblock Spectra of large block matrices.
\newblock cs.IT/0610045, 2006.

\bibitem[Fuglede and Kadison(1952)]{fuglede_kadison52}
Bent Fuglede and Richard~V. Kadison.
\newblock Determinant theory in finite factors.
\newblock \emph{Annals of Mathematics}, 55:\penalty0 520--530, 1952.

\bibitem[Guionnet et~al.(2011)Guionnet, Krishnapur, and
  Zeitouni]{guionnet_krishnapur_zeitouni11}
Alice Guionnet, Manjunath Krishnapur, and Ofer Zeitouni.
\newblock The single ring theorem.
\newblock \emph{Annals of Mathematics}, 174:\penalty0 1189--1217, 2011.
\newblock preprint available at arxiv:0909.2214.

\bibitem[Haagerup(1997)]{haagerup97}
Uffe Haagerup.
\newblock On \mbox{V}oiculescus \mbox{$R$}- and \mbox{$S$}-transforms for free
  non-commuting random variables.
\newblock In Dan-Virgil Voiculescu, editor, \emph{Free Probability Theory},
  volume~12 of \emph{Fields Institute Communications}, pages 127--148. American
  Mathematical Society, 1997.

\bibitem[Haagerup and Larsen(2000)]{haagerup_larsen00}
Uffe Haagerup and Flemming Larsen.
\newblock Brown\mbox{$'$}s spectral distribution measure for
  \mbox{$R$}-diagonal elements in finite von \mbox{N}eumann algebras.
\newblock \emph{Journal of Functional Analysis}, 176:\penalty0 331--367, 2000.

\bibitem[Haagerup and Schultz(2007)]{haagerup_schultz2007}
Uffe Haagerup and Hanne Schultz.
\newblock Brown measures of unbounded operators affiliated with a finite von
  \mbox{N}eumann algebra.
\newblock \emph{Math. Scand.}, 100:\penalty0 209--263, 2007.

\bibitem[Helton et~al.(2007)Helton, Far, and Speicher]{hfs2007}
J.~William Helton, Reza~Rashidi Far, and Roland Speicher.
\newblock Operator-valued semicircular elements: Solving a quadratic matrix
  equation with positivity constraints.
\newblock \emph{International Mathematics Research Notices}, 2007, 2007.

\bibitem[Hiai and Petz(2000)]{hiai_petz00}
Fumio Hiai and Denes Petz.
\newblock \emph{The Semicircle Law, Free Random Variables And Entropy},
  volume~77 of \emph{Mathematical Surveys and Monographs}.
\newblock American Mathematical Society, 1 edition, 2000.

\bibitem[Lehner and Szpojankowski(2021)]{lehner_szpojankowski2021}
Franz Lehner and Kamil Szpojankowski.
\newblock Boolean cumulants and subordination in free probability.
\newblock \emph{Random Matrices: Theory and Applications}, 10, 2021.
\newblock https://arxiv.org/abs/1907.11442.

\bibitem[Marcenko and Pastur(1967)]{marchenko_pastur67}
V.A Marcenko and L.A Pastur.
\newblock Distribution of eigenvalues of some sets of random matrices.
\newblock \emph{Mathematics in U.S.S.R}, 1:\penalty0 507--536, 1967.

\bibitem[Mingo and Speicher(2017)]{mingo_speicher2017}
J.~A. Mingo and R.~Speicher.
\newblock \emph{Free Probability and Random Matrices}.
\newblock Schwinger-Verlag, 2017.

\bibitem[Nayak(2018)]{nayak2018}
Soumyashant Nayak.
\newblock The \mbox{H}adamard determinant inequality -- \mbox{E}xtensions to
  operators on a \mbox{H}ilbert space.
\newblock \emph{J. Funct. Anal.}, 274:\penalty0 2978--3002, 2018.

\bibitem[Nica and Speicher(1997{\natexlab{a}})]{nica_speicher97}
Alexandru Nica and Roland Speicher.
\newblock \mbox{$R$}-diagonal pairs: A common approach to \mbox{H}aar unitaries
  and circular elements.
\newblock In Dan-Virgil Voiculescu, editor, \emph{Free Probability Theory},
  volume~12 of \emph{Fields Institute Communications}, pages 149--188. American
  Mathematical Society, 1997{\natexlab{a}}.

\bibitem[Nica and Speicher(1997{\natexlab{b}})]{nica_speicher97b}
Alexandru Nica and Roland Speicher.
\newblock A ``\mbox{F}ourier transform" for multiplicative functions on
  non-crossing partitions.
\newblock \emph{Journal of Algebraic Combinatorics}, 6:\penalty0 141--160,
  1997{\natexlab{b}}.

\bibitem[Nica and Speicher(2006)]{nica_speicher06}
Alexandru Nica and Roland Speicher.
\newblock \emph{Lectures on the Combinatorics of Free Probability}, volume 335
  of \emph{London Mathematical Society Lecture Note Series}.
\newblock Cambridge University Press, 2006.

\bibitem[Rao and Speicher(2007)]{rajrao_speicher2007}
N.~Raj Rao and Roland Speicher.
\newblock Multiplication of free random variables and the \mbox{$S$}-transform:
  the case of vanishing mean.
\newblock \emph{Electronic Communications in Probability}, 12:\penalty0
  248--258, 2007.

\bibitem[Rota(1964)]{rota64}
Gian-Carlo Rota.
\newblock On the foundations of combinatorial theory i. theory of \mbox{M}obius
  functions.
\newblock \emph{Zeit. Fur Wahrscheinlichkeitstheorie und Verw. Gebiete},
  2:\penalty0 340--368, 1964.

\bibitem[Ryan(1998)]{ryan98}
Oyvind Ryan.
\newblock On the construction of free random variables.
\newblock \emph{Journal of Functional Analysis}, 154:\penalty0 291--322, 1998.

\bibitem[Speed(1983)]{speed83}
T.~P. Speed.
\newblock Cumulants and partition lattices.
\newblock \emph{Australian Journal of Statistics}, 25:\penalty0 378--388, 1983.

\bibitem[Speicher(1998)]{speicher98}
R.~Speicher.
\newblock \emph{Combinatorial Theory of the Free Product with Amalgamation and
  Operator-Valued Free Probability Theory}, volume 627 of \emph{Memoirs of
  American Mathematical Society}.
\newblock A.M.S. Providence, RI, 1998.

\bibitem[Speicher(2019)]{speicher2019}
Roland Speicher.
\newblock Non-commutative distributions. lecture notes.
\newblock https://arxiv.org/abs/2009.03589, 2019.

\bibitem[Stanley(2012)]{stanley2012}
Richard~P. Stanley.
\newblock \emph{Enumerative Combinatorics, volume 1}.
\newblock Cambridge University Press, second edition, 2012.

\bibitem[Stanley(2015)]{stanley2015}
Richard~P. Stanley.
\newblock \emph{Catalan Numbers}.
\newblock Cambridge University Press, 2015.

\bibitem[Voiculescu et~al.(1992)Voiculescu, Dykema, and
  Nica]{voiculescu_dykema_nica92}
D.~Voiculescu, K.~Dykema, and A.~Nica.
\newblock \emph{Free Random Variables}.
\newblock A.M.S. Providence, RI, 1992.
\newblock CRM Monograph series, No.1.

\bibitem[Voiculescu(1983)]{voiculescu83}
Dan Voiculescu.
\newblock Symmetries of some reduced free product $\mbox{C}^*$-algebras.
\newblock In \emph{Lecture Notes in Mathematics}, volume 1132, pages 556--588.
  Springer-Verlag, New York, 1983.

\bibitem[Voiculescu(1986)]{voiculescu86}
Dan Voiculescu.
\newblock Addition of certain non-commuting random variables.
\newblock \emph{Journal of Functional Analysis}, 66:\penalty0 323--346, 1986.

\bibitem[Voiculescu(1987)]{voiculescu87}
Dan Voiculescu.
\newblock Multiplication of certain non-commuting random variables.
\newblock \emph{Journal of Operator Theory}, 18:\penalty0 223--235, 1987.

\bibitem[Voiculescu(1991)]{voiculescu91}
Dan Voiculescu.
\newblock Limit laws for random matrices and free products.
\newblock \emph{Inventiones mathematicae}, 104:\penalty0 201--220, 1991.

\bibitem[Voiculescu(1993)]{voiculescu93}
Dan Voiculescu.
\newblock The analogues of entropy and of \mbox{F}isher\mbox{'}s information
  measure in free probability theory \mbox{I}.
\newblock \emph{Communications in Mathematical Physics}, 155:\penalty0 71--92,
  1993.

\bibitem[Voiculescu(2000)]{voiculescu2000}
Dan Voiculescu.
\newblock The coalgebra of the free difference quotient and free probability.
\newblock \emph{International Mathematics Research Notices}, pages 79--106,
  2000.

\bibitem[Voiculescu(2002)]{voiculescu2002}
Dan Voiculescu.
\newblock Analytic subordination consequences of free \mbox{M}arkovianity.
\newblock \emph{Indiana University Mathematics Journal}, 51:\penalty0
  1161--1166, 2002.

\bibitem[Weingarten(1978)]{weingarten78}
Don Weingarten.
\newblock Asymptotic behavior of group integrals in the limit of infinite rank.
\newblock \emph{Journal of Mathematical Physics}, 19:\penalty0 999--1001, 1978.

\end{thebibliography}

\printindex
\end{document}